\documentclass[opre,nonblindrev]{informs4}
\OneAndAHalfSpacedXI
\usepackage[linktocpage,colorlinks,linkcolor=blue,anchorcolor=blue,citecolor=blue,urlcolor=blue,pagebackref]{hyperref}
\hypersetup{}
\renewcommand*{\backref}[1]{\ifx#1\relax \else ~ \fi}
\renewcommand*{\backrefalt}[4]{%
  \ifcase #1 \footnotesize{(Not cited.)}%
  \or        \footnotesize{(Cited on page~#2.)}%
  \else      \footnotesize{(Cited on pages~#2.)}%
  \fi
}

\usepackage{amsmath,amssymb,amsfonts}

\usepackage{natbib}
 \bibpunct[, ]{(}{)}{,}{a}{}{,}%
 \def\bibfont{\small}%

\usepackage{anyfontsize}
\usepackage{endnotes}
\usepackage{rotating}
\usepackage{fancyvrb}
\usepackage{mathtools}
\usepackage{algorithm}
\usepackage{algorithmic}
\usepackage{multirow}

\def\E{{\mathbb E}}
\def\P{{\mathbb P}}
\def\R{{\mathbb R}}

\def\ca{{\mathcal A}}

\def\ci{{\mathcal I}}

\def\cm{{\mathcal M}}
\def\cn{{\mathcal N}}
\def\co{{\mathcal O}}

\def\cs{{\mathcal S}}

\def\cx{{\mathcal X}}
\def\cy{{\mathcal Y}}
\def\cz{{\mathcal Z}}

\TheoremsNumberedThrough     
\ECRepeatTheorems
\JOURNAL{Operations Research}

\EquationsNumberedThrough    

\MANUSCRIPTNO{IJDS-0001-1922.65}

\begin{document}



\RUNAUTHOR{Chen, Hu, and Zhao}

\RUNTITLE{Landscape of Finite-Horizon MDPs with General States and Actions}

\TITLE{Landscape of Policy Optimization for Finite Horizon MDPs with General State and Action}


\ARTICLEAUTHORS{%
\AUTHOR{Xin Chen}
\AFF{H. Milton Stewart School of Industrial and Systems Engineering, Georgia Institute of Technology}

\AUTHOR{Yifan Hu}
\AFF{Department of Statistics, Rutgers University
}

\AUTHOR{Minda Zhao}
\AFF{H. Milton Stewart School of Industrial and Systems Engineering, Georgia Institute of Technology}
}

\ABSTRACT{Policy gradient methods are widely used in reinforcement learning. Yet, the nonconvexity of policy optimization poses significant challenges in understanding the global convergence of policy gradient methods. For a class of finite-horizon Markov Decision Processes (MDPs) with general state and action spaces, we identify a set of structural properties to establish a benign nonconvex landscape, the Polyak-{\L}ojasiewicz-Kurdyka (P{\L}K) condition of the policy optimization. Leveraging the P{\L}K condition, policy gradient methods converge to the globally optimal policy with a non-asymptotic rate despite nonconvexity. Our results apply to various control and operations models, including entropy-regularized tabular MDPs, Linear Quadratic Regulator problems, and both stochastic inventory models and stochastic cash balance problems with strongly convex costs. In these models, stochastic policy gradient methods obtain an $\epsilon$-optimal policy using a sample size of $\tilde{\co}(\epsilon^{-1})$ and polynomial in terms of the planning horizon. To the best of our knowledge, we provide the first sample-complexity guarantees for multi-period inventory systems with Markov-modulated demand and for stochastic cash balance problems. We complement the theory with numerical experiments showing that policy gradient methods outperform several benchmark algorithms from the literature across these operations models.
}


\maketitle

\section{Introduction}
\label{section: introduction}
Reinforcement Learning (RL) has achieved remarkable success in various real-world applications, including the game of Go \citep{silver2016mastering} and robotics \citep{hwangbo2019learning}. An important class of algorithms for solving these RL problems is policy gradient methods, which search over a parameterized policy space by applying (stochastic) gradient methods to the total expected cost of a Markov Decision Process (MDP). This approach is particularly attractive when the parameterized policy class contains (or closely approximates) the optimal policy. For example, \citet{kunnumkal2008using} apply stochastic gradient methods to the classical inventory model by optimizing over the base-stock policy class, which contains the optimal policy. More broadly, policy gradient methods have been studied in multi-echelon inventory systems \citep{glasserman1995sensitivity}, control problems \citep{hu2023toward}, and many other domains.

Despite wide applicability, our understanding of the global and non-asymptotic convergence behavior of policy gradient methods remains limited, largely because the underlying policy optimization problem is generally nonconvex \citep{agarwal2021theory, bhandari2024global}. To help bridge this gap, we study the nonconvex landscape of the policy optimization and establish the Polyak-{\L}ojasiewicz-Kurdyka (P{\L}K) condition \citep{polyak1963gradient, lojasiewicz1963topological, kurdyka1998gradients}. Informally, the P{\L}K condition states that the norm of the gradient dominates the suboptimality gap. It is a relaxation of strong convexity while maintaining a key property that any point satisfying the first-order necessary optimality condition \citep{nocedal1999numerical} is globally optimal. Since policy gradient methods are designed to seek such stationary points, the P{\L}K condition provides a natural route to establishing global convergence guarantees for policy gradient methods on nonconvex MDP problems.

Existing results that establish P{\L}K conditions for policy optimization are often developed for special MDP classes (see, for example, \citealp{fazel2018global}). Consequently, general structural conditions under which the policy optimization satisfies the P{\L}K condition are less developed, and it remains unclear how to extend existing analysis to general MDPs with analogous structures. To address this issue, we identify a set of structural properties that ensure the P{\L}K condition for the policy optimization in finite-horizon MDPs with general state and action spaces, assuming the finite-dimensional parameterized policy class contains an optimal policy. In particular, we show that the P{\L}K condition holds when (i) the policy objective has bounded gradients, (ii) the expected optimal Q-value functions satisfy the P{\L}K condition, and (iii) sequential decomposition inequalities hold. Roughly speaking, sequential decomposition inequalities control the differences between the policy gradients for the current policy and the policy whose parameter at a later period is replaced by its optimal counterpart. The resulting change is controlled by the suboptimality gap of the corresponding expected optimal Q-value function at that period.

The proposed structural conditions provide a general framework to establish the P{\L}K condition for policy optimization. To illustrate this, we verify conditions (i)-(iii) for a range of control and operations models with optimal policy classes. The examples include (1) entropy-regularized tabular MDPs with stochastic policies, (2) linear quadratic regulator (LQR) problems with linear policies, (3) multi-period inventory systems with Markov-modulated demand and strongly convex costs under state-dependent base-stock policies, and (4) stochastic cash balance problems with strongly convex costs under two-sided base-stock policies. A common thread across these models is that the associated dynamic programming recursions exhibit a (hidden) convexity structure, which is crucial for verifying the P{\L}K condition.

Leveraging the P{\L}K condition, we establish linear convergence of the exact policy gradient method to an $\epsilon$-optimal policy for entropy-regularized tabular MDPs and LQR problems, which aligns with existing results \citep{bhandari2024global, hambly2021policy}. For inventory systems with Markov-modulated demand and strongly convex costs, as well as stochastic cash balance problems with strongly convex costs, we further prove an $\tilde{\mathcal{O}}(\epsilon^{-1})$ sample complexity for stochastic policy gradient methods to obtain an $\epsilon$-optimal policy, yielding the first sample-complexity results in the literature. In all cases, the resulting bounds scale polynomially on the time horizon. Finally, we complement our theory with numerical experiments showing that policy gradient methods outperform several benchmark algorithms from the literature across these operations models.

\subsection{Highlights of Contributions}
Our work contributes to \textit{optimization}, \textit{operations}, and \textit{reinforcement learning}.

\begin{itemize}
    \item From an \textit{optimization} perspective, we complement the literature by identifying a concrete class of problems that provably satisfy the P{\L}K condition. In contrast, much of the prior literature takes the P{\L}K condition as a standing assumption to derive global convergence guarantees for first-order methods \citep{attouch2013convergence, lewis2025complexity}. 
    \item From an \textit{operations} perspective, to the best of our knowledge, we provide the first sample complexity guarantees for inventory systems with Markov-modulated demand and for stochastic cash balance problems. Our bounds scale polynomially with the time horizon. In particular, for the classical multi-period inventory control problem, our results improve over the exponential dependence implied by the framework of \citet{huh2014online}. This improvement is enabled by a new technical result (Lemma~\ref{lemma: sequence}), which could be of independent interest.
    \item From a \textit{reinforcement learning} perspective, we identify structural conditions for general MDPs that ensure the P{\L}K condition for policy optimization. Under the P{\L}K condition, policy gradient methods converge to the globally optimal policy with a non-asymptotic rate. In contrast, existing global convergence guarantees for policy gradient methods typically focus on special settings, most notably LQR problems \citep{fazel2018global, hambly2021policy} and tabular MDPs \citep{agarwal2021theory, lan2023policy}.
\end{itemize}

Relatedly, \citet[Theorem 2]{bhandari2024global} identify different structural conditions under which policy gradient methods converge to the globally optimal policy. In particular, they establish the gradient dominance condition for policy optimization in infinite-horizon discounted MDPs, which in turn yields non-asymptotic convergence guarantees for policy gradient methods. Their result requires two assumptions: gradient dominance for all Q-value functions and closure under policy improvement. Although their framework can extend to finite-horizon MDPs, these assumptions do not hold for the inventory model studied in this paper. In contrast, our conditions are tailored to structural properties that arise in a broad class of control and operations problems, including inventory and cash balance models. Thus, our results establish the nonconvex landscape and the first non-asymptotic complexity for these problems. 

Beyond a broader model coverage, we provide a more quantitative characterization of the policy optimization landscape by establishing the P{\L}K condition. \citet[Theorem 3]{bhandari2024global} presents structural conditions for finite-horizon MDPs under which first-order stationary points are globally optimal, leading to an asymptotic convergence rate for policy gradient methods. Yet, it remains unclear if asymptotic convergence translates to an exponential or a polynomial dependence on the planning horizon. In fact, they ``leave the study of a gradient dominance condition for finite-horizon problems as future work." Our work bridges this gap by establishing the P{\L}K condition for policy optimization. The established P{\L}K condition enables us to achieve a non-asymptotic convergence rate for policy gradient methods, yielding a sample complexity of $\tilde{\mathcal{O}}(\epsilon^{-1})$ with polynomial dependence on the time horizon.

\subsection{Related Literature}
This study intersects with three streams of literature: (i) nonconvex landscape conditions that ensure global convergence for algorithms, (ii) global optimality guarantees for policy gradient methods, and (iii) data-driven operations management.

\subsubsection*{Nonconvex Landscape Conditions}
In nonconvex optimization, several landscape conditions guarantee convergence to global optimality for algorithms, including the hidden convexity \citep{stern1995indefinite, ben1996hidden}, i.e., the problem admits a convex reformulation, the Polyak-{\L}ojasiewicz (P{\L}) condition \citep{polyak1963gradient, lojasiewicz1963topological}, the Polyak-{\L}ojasiewicz-Kurdyka (P{\L}K) condition \citep{kurdyka1998gradients, bolte2007lojasiewicz, bento2025convergence}, and others as summarized in \citet{karimi2016linear}.

Hidden convexity has emerged in numerous modern applications, such as policy optimization in convex RL \citep{zhang2020variational, sun2021learning}, supply chain and revenue management \citep{feng2018supply, chen2018preservation, chen2019stochastic, chen2023network, miao2025network}. In the optimization community, \citet{stern1995indefinite} and \citet{ben1996hidden} studied hidden convexity in quadratic programming. Since then, several works have developed tools to identify hidden convexity \citep{ben2011hidden, ben2014hidden} and studied the landscape of hidden convex functions \citep{levin2024effect}. Subsequently, research focuses on algorithm design and non-asymptotic convergence guarantees under hidden convexity. For instance, \cite{chen2024efficient} addressed a specific problem arising from network revenue management. Building on a convex reformulation presented by \citet{feng2018supply}, \cite{chen2024efficient} developed Mirror Stochastic Gradient Descent which achieves global convergence with a sample complexity of order $\tilde{\co}(\epsilon^{-2})$. Moreover, \cite{fatkhullin2025stochastic} extended the results to general-purpose stochastic optimization under hidden convexity, demonstrating an $\co(\epsilon^{-3})$ sample complexity for projected Stochastic Gradient Descent (SGD).

Recently, several works have investigated the Polyak-{\L}ojasiewicz (P{\L}) condition, first introduced by \citet{polyak1963gradient} and \citet{lojasiewicz1963topological}. The P{\L} condition ensures a linear convergence rate for Gradient Descent \citep{karimi2016linear} and an $\co(\epsilon^{-1})$ sample complexity for SGD \citep{hu2024multi} to achieve an $\epsilon$-optimal solution. Additionally, \citet{karimi2016linear} and \citet{liao2024error} explored the connection between the P{\L} condition and other nonconvex landscape conditions, such as the error bounds \citep{luo1993error}, the quadratic growth \citep{anitescu2000degenerate}, the weak strong convexity \citep{necoara2019linear}, and the proximal P{\L} condition (gradient dominance). These conditions ensure that any first-order stationary point is globally optimal, leading to a global convergence of first-order methods despite nonconvexity.

Our study builds on the Polyak-{\L}ojasiewicz-Kurdyka (P{\L}K) condition \citep{kurdyka1998gradients, bolte2007lojasiewicz}, which extends the P{\L} condition to handle more general settings, particularly for non-smooth functions and constrained optimization problems. Assuming the optimization problem satisfies the P{\L}K condition, \citet{attouch2009convergence} demonstrated the global convergence of the classic proximal-point algorithm. Furthermore, \citet{attouch2013convergence} proved a convergence result for descent methods satisfying several conditions, which includes the forward-backward splitting algorithm.

\subsubsection*{Global Convergence for Policy Gradient Methods} 
Several works have analyzed the global convergence of policy gradient methods for MDPs with finite state and action spaces. Among these, \citet{agarwal2021theory} provided a comprehensive analysis of policy gradient methods solving discounted infinite-horizon MDPs. They demonstrated iteration complexities of $\co(\epsilon^{-2})$ for exact policy gradient methods and $\co(\epsilon^{-1})$ for exact Natural Policy Gradient (NPG) methods with softmax parameterization. \citet{cen2022fast} studied entropy-regularized NPG methods in conjunction with softmax parameterization. They established a linear convergence rate for exact entropy-regularized NPG methods solving entropy-regularized MDPs and further proved an $\tilde{\co}(\epsilon^{-2})$ sample complexity for approximate entropy-regularized NPG methods. \citet{lan2023policy} introduced stochastic Policy Mirror Descent (PMD) approaches and demonstrated $\co(\epsilon^{-1})$ (resp., $\co(\epsilon^{-2})$) sample complexity for solving RL problems with strongly convex (resp., convex) regularizers. Following these works, \citet{klein2024beyond} investigated the policy gradient methods with a softmax parameterization for finite-horizon MDPs. They established a weak P{\L} condition of the policy gradient objective function and proved an $\co(\epsilon^{-1})$ iteration complexity for exact policy gradient methods. For more references, we refer the readers to \citet{mei2020global}, \citet{xiao2022convergence}, and \citet{fatkhullin2023stochastic}.

When solving MDPs with general state and action spaces, it is impractical to extend existing results directly as the complexity for policy gradient methods solving tabular MDPs depends on the cardinalities of state and action sets \citep{lan2023policy, klein2024beyond}. To address this issue, many works impose additional assumptions to ensure the global convergence of policy gradient methods. For instance, when dealing with discounted infinite-horizon MDPs, \citet{ju2022policy} demonstrated a linear convergence rate of PMD when the advantage function is convex in action and the regularizer is strongly convex. \citet[Theorem 2]{bhandari2024global} identified certain structural characteristics shared by several discounted infinite-horizon control problems, which guarantee that the policy gradient objective function satisfies the so-called $(c,\mu)$-gradient dominance condition \citep[Definition 2]{bhandari2024global}, a condition equivalent to the P{\L}K condition under some mild assumptions \citep{karimi2016linear}. They established iteration complexities of $\co(\epsilon^{-2})$ and $\co(\log(\epsilon^{-1}))$ for exact policy gradient methods to achieve an $\epsilon$-optimal policy under the $(c,0)$-gradient dominance condition and $(c,\mu)$-gradient dominance condition with $\mu > 0$, respectively.

For finite-horizon MDPs, \citet[Theorem 3]{bhandari2024global} established conditions for a class of finite-horizon MDPs under which first-order stationary points are globally optimal, leading to an asymptotic convergence rate for policy gradient methods \citep[Lemma 2]{bhandari2024global}. Yet, they ``leave the study of a gradient dominance condition for finite-horizon problems as future work". Our work bridges this gap by establishing the P{\L}K condition of the policy gradient optimization for such problems.

In addition to the results for general MDPs, \citet{fazel2018global} established the gradient domination of the policy gradient objective function for the infinite-horizon discounted LQR problem. They demonstrated a linear convergence rate for exact policy gradient methods to achieve an optimal policy. \citet{han2025policy} extended results to infinite-horizon nearly linear quadratic regulators, where the dynamic system integrates linear and nonlinear components. They established a linear convergence rate of policy gradient methods using the nonlinear policy class to achieve globally optimal policies. More closely related to our work, \citet{hambly2021policy} studied policy gradient methods for solving the finite-horizon LQR problem. Their complexity admits a polynomial dependence on the time horizon. However, the applicability of their analysis is limited due to the specific structures inherent to LQR problems, making it challenging to generalize to other applications.

\subsubsection*{Data-driven Operations Management} 
Most of the literature studying multi-period operations models relies on the Sample Average Approximation (SAA) approach, which constructs an empirical objective using samples and then solves the corresponding empirical problem \citep{kleywegt2002sample}. For instance, \citet{levi2007provably} established the sample complexity required for SAA to find an $\epsilon$-optimal base-stock policy of the multi-period inventory system. \citet{cheung2019sampling} applied the SAA method for multi-period capacitated stochastic inventory control problems and derived a sample complexity required to achieve a near-optimal expected cost. They further proposed a polynomial-time approximation scheme that also uses polynomially many samples to solve the empirical counterpart. \citet{qin2022data} investigated multi-period joint pricing and inventory control models and established the sample complexity for the SAA approach. They applied a sparsification technique and proposed a polynomial-time approximation algorithm for the empirical problem. \citet{zhang2025sampling} applied the SAA method for managing inventories in an infinite-horizon series system with multiple stages and derived sample complexities. \citet{xie2024vc} investigated uniform generalization errors of the SAA approach for different inventory policy classes and established a sample complexity independent of the time horizon length for learning a base-stock policy which incurs an averaged cost no more than the optimal averaged cost plus $\epsilon$. 

In addition to SAA methods, some works apply value-based methods to solve multi-period operations problems. For example, \citet{qin2023sailing} introduced a variance-reduced value iteration algorithm for multi-period stochastic inventory control with independent demands, establishing matching upper and lower bounds on the sample complexity. \citet{gong2023bandits} developed online Q-learning methods for stochastic inventory models with cyclic demands. They considered two scenarios: the episodic model where inventory is discarded at the end of each cycle and the non-discarding case. We remark that the SAA approach and value-based methods usually rely on solving the empirical counterpart through dynamic programming, whereas policy gradient methods directly optimize a single objective.

Similar to policy gradient methods, some studies apply stochastic gradient methods to solve data-driven operations models. For instance, \citet{kunnumkal2008using} proposed a biased stochastic gradient method to solve finite-horizon inventory systems with independent demands and established an asymptotic convergence rate to achieve optimal base-stock levels. \citet{huh2014online} applied the same biased stochastic gradient method for a class of multistage stochastic optimization problems, where the objective satisfies a generalized convex condition called the sequentially convex condition. 

Though our sequential decomposition inequalities are inspired by one of the conditions in sequential convexity~\citet{huh2014online}, our work is fundamentally different from theirs. First, \citet{huh2014online} (and \citet{kunnumkal2008using}) applied biased gradient methods to minimize cost-to-go functions, which differs from standard policy gradient methods that optimize the objective function directly as studied in the literature and this work. Second, the sequential convexity proposed in \citet{huh2014online} cannot be directly applied to the inventory system with Markov-modulated demand and stochastic cash balance problems. For example, the minimizers of cost-to-go functions differ from the optimal policy parameters of stochastic cash balance problems. On the other hand, our framework is readily applicable to the operations applications in \citet{huh2014online}, including inventory control, capacity allocation, and lifetime buy decision problems. Third, \citet{huh2014online}, leveraging the sequential convexity, proved that biased stochastic gradient methods achieve a sample complexity with an exponential dependence on the planning horizon. In contrast, we show that policy gradient optimization satisfies the P{\L}K condition and policy gradient methods admit a sample complexity that scales polynomially with the planning horizon.

\subsection{Organization}
The rest of this paper is structured as follows. Section~\ref{section: problem formulation} presents the MDP formulation and the policy optimization problem. Section~\ref{section: landscape} outlines the definition and properties of the P{\L}K condition and identifies structural conditions required to ensure the P{\L}K condition of the policy optimization. In Sections~\ref{section: tabular MDPs} and~\ref{section: LQR}, we verify that the policy optimization for entropy-regularized tabular MDPs and LQR, respectively, satisfy the P{\L}K condition. Sections~\ref{section: inventory system} and \ref{section: cash balance} establish the P{\L}K condition for the policy optimization of the inventory system with Markov-modulated demand and the stochastic cash balance problem, respectively. Utilizing the P{\L}K condition, we provide the first sample complexity results for these settings in the literature. Section~\ref{sec: numerical exp} presents numerical experiments benchmarking policy gradient methods on operations models against established algorithms from the literature. The results demonstrate that policy gradient methods achieve strong solution quality and remain computationally efficient.

\subsection{Notations and Definitions}
We use the following notations throughout the paper. Let $\bar{\R} = \R \cup \{+\infty\}$. $\mathbb{N}_+$ denotes the set of all natural numbers. For $n\in\mathbb{N}_+$, denote $[n]$ as the set $\{1, \dots, n\}$. $\|x\|_2 = \sqrt{\sum_{i=1}^n|x_i|^2}$ denotes the $l_2$-norm of a vector $x\in\R^n$. We use $\|A\|_2$ to denote the spectral norm of a matrix $A\in\R^{m \times n}$, the largest singular value of $A$. Let $\|A\|_F = \sqrt{\sum_{i=1}^m \sum_{j=1}^n |a_{ij}|^2}$ denote the Frobenius norm of a matrix $A\in\R^{m \times n}$. We use $r_\sigma(A)$ to denote the spectral radius of a square matrix $A\in\R^{n \times n}$, which is the largest modulus of the eigenvalues of $A$. Let $e$ denote $\exp(1)$. We use $\lfloor x \rfloor$ to denote the greatest integer less than or equal to $x$. For $D\in\R^n$, define $D_{[j, k]} \coloneqq \sum_{i=j}^kD_i$ for $1\le j \le k \le n$. We say a point $x\in\cx$ satisfies the first-order necessary optimality condition of the optimization problem $\min_{x\in\cx} f(x)$ if $\langle \nabla f(x), x' - x\rangle \ge 0, \forall x'\in\cx$ for a differentiable function $f$. A point $\bar{x}$ is an $\epsilon$-optimal solution of $\min_{x} f(x)$ if $f(\bar{x}) - \min_x f(x) \le \epsilon$. We use $\co(\cdot)$ to denote the order in terms of $\epsilon^{-1}$ and $\tilde{\co}(\cdot)$ to denote the order hiding the logarithmic dependency on $\epsilon^{-1}$.

\section{Problem Formulation}
\label{section: problem formulation}
We specify a finite-horizon MDP $\cm = (\mathcal{S}, \mathcal{A}, P, C, T, \rho)$ defined in \citet{puterman2014markov}: the time horizon $T$; the state space $\cs=\cs_1\cup\dots\cup\cs_T$, where $\cs_t\subseteq \R^m$ is the feasible region for a state $s$ at period $t$; the action space $\ca = \cup_{s\in\cs}\ca_s$, where $\ca_s\subseteq \R^n$ is the set of feasible actions for state $s \in \cs \subseteq \R^m$; the transition kernel $P: \cs \times \ca \times [T] \to \cs$, where $P(s'|s, a, t)$ is the probability density function (or probability mass function in the discrete setting) of transitioning into $s'$ when taking action $a$ in state $s$ at period $t$; the cost function $C: \cs \times \ca \times [T] \to \R$, where $C(s, a, t)$ is the immediate cost after taking action $a$ in state $s$ at period $t$; and the initial state distribution $\rho$. For simplicity, we use $P_t(\cdot|s,a) \coloneqq P(\cdot|s,a,t)$, and $C_t(s,a) \coloneqq C(s,a,t)$ for all $s \in \cs, a \in \ca, t\in[T]$. The agent starts at state $s_1 \in \cs_1$, which follows the initial state distribution $\rho$. At period $t$, the agent first observes the current state $s_t \in \cs_t$ and then takes an action $a_t\in \ca_{s_t}$. Afterwards, it receives an immediate cost $C_t(s_t,a_t)$ and proceeds to the next period with state $s_{t+1}\sim P_t(\cdot | s_t,a_t)$. 

A non-stationary policy $\pi:\cs\times[T]\to\ca$ is a function that maps the current state $s$ to a feasible action $a$ at period $t$, e.g., $a = \pi(s, t)$. Similarly, we use $\pi_t(\cdot)$ to denote the policy at period $t$ and $\pi_t(s) \coloneqq \pi(s, t)$ for all $s\in\cs, t\in[T]$. Let $\Pi$ denote the set of feasible policies and $\Pi_t$ denote the set of feasible policies at period $t$. For any $\pi\in\Pi$, the total expected cost starting from state $s$ is
\begin{equation*}
    \begin{aligned}
        J^\pi(s) = \E \left[ \sum_{t=1}^T C_t \left( s_t, \pi_t(s_t) \right) \Biggl| s_1 = s, \pi \right].
    \end{aligned}
\end{equation*}
We take the expectation over a Markovian sequence $(s_1, \dots, s_T)$, where $s_1$ is the initial state and $s_{t+1}\sim P_t(\cdot|s_t, \pi_t(s_t))$ for all $t=1,\dots,T-1$. A policy $\pi^*$ is optimal if it minimizes the total expected cost $J(\pi)$ with the initial distribution $\rho$:
\begin{equation*}
    \begin{aligned}
        J(\pi) = \E_{s\sim \rho} \left[ J^\pi(s) \right] = \E \left[ \sum_{t=1}^T C_t \left( s_t, \pi_t(s_t) \right) \Biggl| s_1 \sim \rho, \pi \right].
    \end{aligned}
\end{equation*}

\subsection{Bellman Equation}
We introduce several terminologies commonly used in the literature of MDPs. Let $\pi$ be a given policy. We define $\rho_t(\cdot|\pi)$ as the cumulative distribution function of $s_t$ incurred by policy $\pi$ starting with the initial distribution $\rho$. By definition, we have $\rho_1(\cdot|\pi) = \rho$. Furthermore, we define the value function $V^{\pi}_t:\cs\to\R$, which represents the total expected cost at time $t$ starting with the initial state $s$ and policy $\pi$:
\begin{equation*}
    \begin{aligned}
        V^{\pi}_t(s) = \E \left[ \sum_{k=t}^T C_k \left( s_k, \pi_k(s_k) \right) \Biggl| s_t = s, \pi \right].
    \end{aligned}
\end{equation*}
In the same manner, we define the function $Q^{\pi}_t:\cs \times \ca \to \R$ as the action-value (or Q-value) function:
\begin{equation*}
    \begin{aligned}
        Q^{\pi}_t(s, a) = C_t(s, a) + \E \left[ \sum_{k=t+1}^T C_k \left( s_k, \pi_k(s_k) \right) \Biggl| s_t = s, a_t = a, \pi \right].
    \end{aligned}
\end{equation*}
By definition, the value function $V^{\pi}$ and action-value function $Q^{\pi}$ have the following relationships:
\begin{equation}
\label{bellman equation}
    \left\{    
    \begin{aligned}
        V^{\pi}_t(s) & = Q^{\pi}_t \left( s, \pi_t(s) \right),\\
        Q^{\pi}_t(s, a) & = C_t(s, a) + \E \left[ V^{\pi}_{t+1}(s') | s'\sim P_t(\cdot | s, a) \right],
    \end{aligned}
    \right.
\end{equation}
for all $s\in\cs, a\in\ca, t\in[T]$ and the boundary condition is $V^{\pi}_{T+1}(\cdot) = 0$ for all $\pi\in\Pi$. These are commonly known as the Bellman equations in the literature \citep{bellman1952theory}. By the principle of optimality \citep{puterman2014markov}, an optimal policy $\pi^*$ solves the following Bellman equations:
\begin{equation}
\label{bellman optimality equation}
    \left\{
    \begin{aligned}
        V^*_t(s) & = \min_{\pi_t\in\Pi_t} Q^*_t \left( s, \pi_t(s) \right),\\
        Q^*_t(s, a) & = C_t(s, a) + \E \left[ V^*_{t+1}(s') | s'\sim P_t(\cdot | s, a) \right],
    \end{aligned}
    \right.
\end{equation}
for all $s\in\cs, a\in\ca, t\in[T]$ and the boundary condition is $V^*_{T+1}(\cdot) = 0$. Here $V_t^*$ and $Q_t^*$ denote the value function and the Q-value function corresponding with the optimal policy $\pi^*$, respectively.

\subsection{Policy Gradient Formulation}
\label{subsection: policy gradient formulation}
Policy gradient methods apply first-order algorithms to minimize the total expected cost $J(\pi)$. Note that general policy optimization falls into functional optimization as we search over the function class $\Pi$, which is computationally intractable. To avoid functional optimization, it is common to parameterize the policy through finite-dimensional parameters $\theta = (\theta_1,\dots,\theta_T)$ \citep{sutton1999policy}. At time $t$, the parameterized policy is $\pi_t(\cdot|\theta_t)$ and $\theta_t$ belongs to a convex and compact set $\Theta_t \subseteq \R^d$. The feasible region of $\theta$ is a product set $\Theta = \Theta_1\times\dots\times\Theta_T$, which is also convex and compact. In such a case, the parameterized policy class is $\Pi_\Theta = \{ \pi(s, t|\theta): \cs \times [T] \times \Theta \to \ca \} \subseteq \Pi$ and we use $\pi_\theta$ to denote $\pi(\cdot|\theta)$ for simplicity.

Given parameterized policy $\pi_\theta$, we represent the total expected cost by $l(\theta) \coloneqq J(\pi_\theta)$, called the policy gradient objective function. We define a policy $\pi_{\theta}$ to be $\epsilon$-optimal if $\theta$ is an $\epsilon$-optimal solution of $l(\theta)$. Let $\theta^*$ denote one of the minimizers of $\min_{\theta\in\Theta} l(\theta)$ and $\pi_{\theta^*}$ as the corresponding policy. Throughout the paper, we assume that the parameterized policy class $\Pi_\Theta$ contains the optimal policy $\pi^*$. This often occurs when the optimal policy class is known, e.g., affine policies in the Linear Quadratic Regulator (LQR) problem and base-stock policies in the multi-period inventory control model. Hence, $V_t^*$ (resp., $Q_t^*$) and $V_t^{\pi_{\theta^*}}$ (resp., $Q_t^{\pi_{\theta^*}}$) are identical.

\section{Landscape Characterization}
\label{section: landscape}
In this section, we first formally define the P{\L}K condition and discuss its properties. Next, we present the non-asymptotic convergence rate of algorithms for optimization problems satisfying the P{\L}K condition. Lastly, we provide a set of verifiable assumptions for the policy gradient optimization in Section~\ref{subsection: policy gradient formulation} to satisfy the P{\L}K condition, enabling a non-asymptotic convergence rate for policy gradient methods to achieve globally optimal solutions.

\subsection{Definition and Properties of P{\L}K Condition}
\label{subsection: definition}
To characterize the landscape of constrained smooth optimization problems, we adopt a particular form of the P{\L}K condition \citep[Appendix G]{karimi2016linear}. For more general formulations that accommodate nonsmooth objectives, we refer readers to \citep[Definition 2.4]{attouch2013convergence}.

\begin{definition}[P{\L}K condition]
\label{def: KL condition}
Let $\mathcal{X}\subseteq\mathbb{R}^n$ be a convex and compact set, and let $f:\mathbb{R}^n\to\mathbb{R}$ be differentiable. Define the optimal value $f^* \coloneqq \min_{x\in\mathcal{X}} f(x)$. We say that $f$ satisfies the P{\L}K condition on $\mathcal{X}$ if there exists $\mu>0$ such that
\begin{equation*}
    f(x)-f^* \le \frac{1}{2\mu}\min_{g\in \partial \delta_{\mathcal{X}}(x)} \left\| \nabla f(x)+g \right\|_2^2, \qquad \forall x\in\mathcal{X},
\end{equation*}
where $\mu$ is the P{\L}K constant and $\delta_{\mathcal{X}}$ is the indicator function of $\mathcal{X}$:
\begin{equation*}
    \delta_{\mathcal{X}}(x)=
    \begin{cases}
    0, & x\in \mathcal{X},\\
    +\infty, & x\notin \mathcal{X}.
    \end{cases}
\end{equation*}
Here $\partial \delta_{\mathcal{X}}(x)$ denotes the (convex) subdifferential of $\delta_{\cx}$ at $x$, which is the normal cone of $\cx$ at $x$.
\end{definition}

\begin{remark}
    When the decision variable $x$ is a matrix, e.g., parameters in the LQR problem (Section~\ref{section: LQR}), we replace the $l_2$-norm with the Frobenius norm.
\end{remark}

When $\cx = \R^n$, the P{\L}K Condition reduces to the P{\L} condition. Similar to the P{\L} condition, the P{\L}K condition is a relaxed condition of the strong convexity (Corollary~\ref{corollary: strong convexity implies KL} in Appendix~\ref{subsection: nonconvex conditions}). It is well known that strongly convex functions exclude all suboptimal stationary or local optimal points. Optimization problems with the P{\L}K condition exhibit the same structural property.

\begin{proposition}[\citealp{karimi2016linear}]
\label{proposition: KL no local optimal}
    Consider a convex and compact set $\cx\subseteq\R^n$. Suppose a function $f:\cx\to\R$ satisfies the P{\L}K condition with a P{\L}K constant $\mu>0$ over $\cx$. Then, any point satisfying the first-order necessary optimality condition of the optimization problem $\min_{x\in\cx} f(x)$ is globally optimal.
\end{proposition}

\subsection{Convergence Rate under P{\L}K Condition}
\label{section: alg}
This subsection presents the global convergence results of the projected (stochastic) gradient descent for solving the stochastic optimization problem over a convex and compact set $\mathcal{X}$ under the P{\L}K condition:
\begin{equation}
\label{stochastic opt}
    \min_{x\in\mathcal{X}} \ f(x) \coloneqq \mathbb{E}_{\xi\sim\mathbb{P}(\xi)} \left[ F(x,\xi) \right],
\end{equation}
where $x$ is the decision variable and $\xi$ is the random variable with a cumulative distribution function $\mathbb{P}(\xi)$. Furthermore, we analyze the algorithm's sample complexity, which refers to the number of samples required to obtain an $\epsilon$-optimal solution.

\begin{assumption}
\label{assumption: smoothness}
    For the function $f:\mathcal{X}\to\mathbb{R}$, we assume that it is differentiable and is $L$-smooth, i.e., its gradient is $L$-Lipschitz:
    \begin{equation*}
        \left \| \nabla f(x) - \nabla f(y) \right \|_2 \le L \|x - y\|_2, \quad \forall x,y\in\mathcal{X}.
    \end{equation*}
\end{assumption}

The smoothness assumption is standard in the stochastic gradient methods literature \citep{nemirovski2009robust, ghadimi2013stochastic}. In our applications, the smoothness holds under mild regularity conditions. We provide a detailed verification in Lemma~\ref{lemma: lipschitz_pg invt control} in Appendix~\ref{subsection: sample complexity} for the multi-period inventory control with Markov-modulated demand.

In what follows, we consider two types of gradient oracle. In the first, an oracle provides the exact gradient $\nabla f(x)=\nabla \mathbb{E}_{\xi\sim \mathbb{P}(\xi)}[F(x,\xi)]$. Because computing $\nabla f(x)$ is often impractical in stochastic optimization, we also study a more realistic oracle in which the algorithm can access only a stochastic gradient estimator constructed from samples drawn from $\mathbb{P}(\xi)$.

\begin{assumption}[Gradient Oracle]
\label{assumption: gradient information}
    We consider one of the following oracles:
    \begin{henumerate}
        \item The exact gradient $\nabla f(x)$ is available.
        \item \label{assumption: stochastic gradient} A stochastic gradient $\nabla F(x, \xi)$ is available such that, conditional on $x$, there exists $\sigma > 0$,
        \begin{equation}
        \label{property:bias_variance}
            \begin{aligned}
                &\mathbb{E}_{\xi} \left[ \nabla F(x, \xi) \right] = \nabla f(x),\\
                &\mathbb{E}_{\xi} \left[ \left\| \nabla F(x, \xi) - \nabla f(x) \right\|_2^2 \right] \le \sigma^2.
            \end{aligned}
        \end{equation}
    \end{henumerate}
\end{assumption}

More specifically, at iteration $k$, let $\{\xi_{k}^{(i)}\}_{i=1}^N$ be i.i.d. fresh samples drawn from $\mathbb{P}(\xi)$ (independent across $k$), we consider the following mini-batch estimator:
\begin{equation*}
    \nabla \hat{f}(x_k) \coloneqq \frac{1}{N}\sum_{i=1}^N \nabla F(x_k, \xi_k^{(i)}).
\end{equation*}
Setting $N=1$ recovers the standard single-sample stochastic gradient estimator. Under Assumption~\ref{assumption: gradient information}.\ref{assumption: stochastic gradient}, the mini-batch estimator $\nabla \hat{f}(x_k)$ is unbiased since
\begin{equation*}
    \mathbb{E}_{\{\xi_k^{(i)}\}_{i=1}^N} \left[ \nabla \hat{f}(x_k) \right] = \mathbb{E}_{\{\xi_k^{(i)}\}_{i=1}^N} \left[ \frac{1}{N}\sum_{i=1}^N \nabla F(x_k, \xi_k^{(i)}) \right] = \nabla f(x_k).
\end{equation*}
Furthermore, since $\{\xi_k^{(i)}\}_{i=1}^N$ are i.i.d. conditional on $x_k$, its variance is reduced by a factor $N$:
\begin{equation}\label{mini-batch estimator}
    \begin{aligned}
        \mathbb{E}_{\{\xi_k^{(i)}\}_{i=1}^N} \left[ \left\| \nabla \hat{f}(x_k) - \nabla f(x_k) \right\|_2^2 \right] &= \mathbb{E}_{\{\xi_k^{(i)}\}_{i=1}^N} \left[ \left\| \frac{1}{N}\sum_{i=1}^N \left( F(x_k, \xi_k^{(i)}) - \nabla f(x_k) \right) \right\|_2^2 \right]\\
        &= \frac{1}{N} \mathbb{E}_{\{\xi_k^{(i)}\}_{i=1}^N} \left[ \left\| F(x_k, \xi_k^{(i)}) - \nabla f(x_k) \right\|_2^2 \right] \le \frac{\sigma^2}{N}.
    \end{aligned}
\end{equation}

Extensive research has studied the convergence of first-order methods for optimization problems satisfying the P{\L}K condition. \citet{attouch2013convergence} presented a general framework for analyzing the convergence of a class of descent methods, in which the projected gradient descent is a special case. Following the same framework, one can extend the convergence results to the stochastic setting. We state the results for completeness and leave the proof in Appendix~\ref{subsection: PGD KL}.
\begin{lemma}
\label{lemma: PGD KL}
    Consider an optimization problem $\min_{x\in\cx} f(x)$ over a convex and compact set $\cx$. Assume $f$ satisfies the P{\L}K condition on $\cx$ with a parameter $\mu > 0$, and Assumptions~\ref{assumption: smoothness} and \ref{assumption: gradient information} hold. Denote $f^* = \min_{x\in\mathcal{X}}f(x)$.
    \begin{henumerate}
        \item \label{alg: DPGD KL} \citep{attouch2013convergence} The sequence of projected gradient descent $x_{k+1} = \text{Proj}_\cx (x_k - \gamma_k \nabla f(x_k) )$ with stepsizes $\gamma_k = \frac{1}{L}$ achieves a linear convergence rate, i.e.,
        \[
        f(x_{k}) - f^* \le \left(1 - \frac{\mu}{4L + \mu} \right)^k \left( f(x_0) - f^* \right).
        \]
        \item \label{alg: SPGD KL} The sequence of projected stochastic gradient descent $x_{k+1} = \text{Proj}_\cx (x_k - \gamma_k \nabla \hat{f}(x_k) )$ with the mini-batch estimator $\nabla \hat{f}(x_k)$ defined in (\ref{mini-batch estimator}) and stepsizes $\gamma_k = \frac{1}{L}$ admits a sublinear convergence rate, i.e.,
        \[
        \mathbb{E}[f(x_k)] - f^* \le \left( 1 - \frac{\mu}{16L + \mu} \right)^k\left( \E[f(x_0)] - f^* \right) + \frac{17\sigma^2}{\mu N}.
        \]
    \end{henumerate}
\end{lemma}
\begin{remark}\label{remark: PGD KL}
    From Lemma~\ref{lemma: PGD KL}.\ref{alg: SPGD KL}, we need to set $N = \co(\epsilon^{-1})$ and $k = \co(\log(\epsilon^{-1}))$ to obtain an $\epsilon$-optimal solution. As a result, the sample complexity of projected stochastic gradient descent is $\tilde{\co}(\epsilon^{-1})$.
\end{remark}

\subsection{P{\L}K Condition in Policy Gradient Formulation}
\label{subsection: KL in PG}
Leveraging the P{\L}K condition, we can establish the global convergence of first-order methods for nonconvex smooth optimization problems. However, verifying the P{\L}K condition for the policy gradient optimization is challenging. To address this difficulty, we develop a general framework to validate the P{\L}K condition for a class of MDPs in the following theorem. To maintain consistency in notation, we present the theorem using the same terminology in Section~\ref{section: problem formulation}.

\begin{theorem}
\label{theorem: main result}
    Consider a Markov Decision Process $\cm = (\cs, \ca, P, C, T, \rho)$ and a policy class $\Pi_\Theta$ with a convex and compact set $\Theta$. Suppose the following conditions hold.
    \begin{henumerate}
        \item \textbf{(Bounded Gradients)} For any $t\in[T]$, the expected Q-value function 
        \begin{equation*}
            \theta_t \mapsto \E_{s_t \sim \rho_t(\cdot|\pi_\theta)} \left[ Q^{\pi_\theta}_t \left(s_t, \pi_t(s_t|\theta_t) \right) \right]
        \end{equation*}
        is continuously differentiable on $\Theta_t$ with the $2$-norm of its gradient upper bounded by $G$.
        \item \textbf{(P{\L}K Condition of Expected Optimal Q-value Functions)} For any $t\in[T]$, the expected optimal Q-value function
        \begin{equation*}
            \theta_t \mapsto \E_{s_t \sim \rho_t(\cdot|\pi_\theta)} \left[ Q^{\pi_{\theta^*}}_t \left(s_t, \pi_t(s_t|\theta_t) \right) \right]
        \end{equation*}
        satisfies the P{\L}K condition with a P{\L}K constant $\mu_Q$ on $\Theta_t$.
        \item \textbf{(Sequential Decomposition Inequality)} There exists $M_g > 0$ such that for any $\theta\in\Theta$ and $1\le t < k \le T$,
        \begin{equation}
        \label{sequential decomposition inequality}
            \begin{aligned}
                & \ \left\| \nabla_{\theta_t} l(\theta_1, \dots, \theta_{k-1}, \theta_k, \theta_{k+1}^*, \dots, \theta_T^*) - \nabla_{\theta_t} l \left( \theta_1, \dots, \theta_{k-1}, \theta_k^*, \theta_{k+1}^* \dots, \theta_T^* \right) \right\|_2\\
                \le & \ M_g \left( \mathbb{E}_{s_k \sim \rho_k(\cdot|\pi_\theta)} \left[ Q_k^{\pi_{\theta^*}} \left( s_k, \pi_k(s_k|\theta_k) \right)\right] - \mathbb{E}_{s_k \sim \rho_k(\cdot|\pi_\theta)} \left[Q_k^{\pi_{\theta^*}} \left( s_k, \pi_k(s_k|\theta^*_k) \right) \right] \right).
            \end{aligned}
        \end{equation}
    \end{henumerate}
    Then the policy gradient objective function $l(\theta)$ satisfies the P{\L}K condition on $\Theta$. Furthermore, the corresponding P{\L}K constant is $\mu_l = \frac{\mu_Q^3}{eM_g^2G^2T^2}$, i.e.,
    \begin{equation*}
        l(\theta) - l(\theta^*) \le \frac{eM_g^2G^2T^2}{2\mu_Q^3} \min_{g \in \partial \delta_{\Theta}(\theta)} \left\|\nabla l(\theta) + g \right\|^2_2, \quad \forall \theta\in\Theta.
    \end{equation*}
\end{theorem}

\begin{remark}
    The parameters $G$, $\mu_Q$, and $M_g$ may depend on the time horizon $T$. The following sections show these parameters exhibit polynomial dependence on $T$ for different applications.
\end{remark}

Theorem~\ref{theorem: main result} provides structural conditions under which the policy optimization satisfies the P{\L}K condition. Table~\ref{table: comparison} summarizes differences between Theorem~\ref{theorem: main result} and \citet{bhandari2024global}. First, Theorem~\ref{theorem: main result} applies to the inventory and cash balance models, which will be verified in Sections~\ref{section: inventory system} and~\ref{section: cash balance}. In contrast, \citet[Theorem~2]{bhandari2024global} does not cover these operations models. Second, relative to \citet[Theorem~3]{bhandari2024global}, Theorem~\ref{theorem: main result} delivers a sharper convergence rate for policy gradient methods, which is polynomially in $T$. While the structural conditions in \citet[Theorem~3]{bhandari2024global} are more general than ours, Theorem~\ref{theorem: main result} works for all applications studied in \citet{bhandari2024global}.

\vspace{-10pt}
\begin{table}[htbp]
    \scriptsize
    \centering
    \caption{Differences between Theorem~\ref{theorem: main result} and \citet{bhandari2024global}.}
    \renewcommand{\arraystretch}{1.2}
    \begin{tabular}{cccc}
        \hline\hline
        & \begin{tabular}[c]{@{}c@{}} Theorem 2 in \\ \citet{bhandari2024global} \end{tabular} & \begin{tabular}[c]{@{}c@{}} Theorem 3 in \\ \citet{bhandari2024global} \end{tabular} & Our Work \\ \hline
            
        \begin{tabular}[c]{@{}c@{}} Problem \\ Setting \end{tabular} & Infinite-Horizon MDPs & Finite-Horizon MDPs & Finite-Horizon MDPs \\ \hline
            
        \begin{tabular}[c]{@{}c@{}} Sufficient \\ Conditions \end{tabular} & \begin{tabular}[c]{@{}c@{}} Gradient dominance of all \\ expected Q-value functions; \\ Closure under policy improvement \end{tabular} & \begin{tabular}[c]{@{}c@{}} No spurious local optimal points for \\ expected optimal Q-value functions; \end{tabular} & \begin{tabular}[c]{@{}c@{}} P{\L}K condition of expected \\ optimal Q-value functions; \\ Sequential decomposition inequality \end{tabular}\\ \hline
            
        \begin{tabular}[c]{@{}c@{}} Landscape of \\ Objective \end{tabular} & Gradient Dominance & No local optimal points & P{\L}K condition \\ \hline

        Convergence Rate & Non-asymptotic & Asymptotic & Non-asymptotic \\ \hline

        Dependence in $T$ & NA & Unknown & Polynomial in $T$ \\ \hline
            
        Applications & \begin{tabular}[c]{@{}c@{}} Inventory models ($\times$) \\ Cash balance problem ($\times$) \end{tabular} & \begin{tabular}[c]{@{}c@{}} Inventory models ($\checkmark$) \\ Cash balance problem ($\checkmark$) \end{tabular} & \begin{tabular}[c]{@{}c@{}} Inventory models ($\checkmark$) \\ Cash balance problem ($\checkmark$) \end{tabular} \\
        \hline\hline
    \end{tabular}
    \label{table: comparison}
\end{table}
\vspace{-10pt}

We provide some intuitions as to why the conditions in Theorem~\ref{theorem: main result} hold for a broad class of MDPs. The bounded gradients condition is standard and is likely to hold for many applications. For the P{\L}K condition of expected optimal Q-value functions, one can verify it using convex cost-to-go functions and strongly convex costs. We discuss more intuitions in the following sections for different cases. The less intuitive condition is the sequential decomposition inequality. A weaker result using standard assumptions can be derived. Suppose $l(\theta)$ is $S_l$-smooth, and expected optimal Q-value functions satisfy the P{\L}K condition:
\begin{equation*}
    \begin{aligned}
        & \ \left\| \nabla_{\theta_t} l(\theta_1, \dots, \theta_k, \theta_{k+1}^*, \dots, \theta_T^*) - \nabla_{\theta_t} l \left( \theta_1, \dots, \theta_{k-1}, \theta_k^*, \dots, \theta_T^* \right) \right\|_2 \le S_l \| \theta_k - \theta_k^* \|_2\\
        \le & \ S_l \sqrt{\frac{2}{\mu_Q}} \sqrt{\mathbb{E}_{s_k \sim \rho_k(\cdot|\pi_\theta)} \left[ Q_k^{\pi_{\theta^*}} \left( s_k, \pi_k(s_k|\theta_k) \right)\right] - \mathbb{E}_{s_k \sim \rho_k(\cdot|\pi_\theta)} \left[Q_k^{\pi_{\theta^*}} \left( s_k, \pi_k(s_k|\theta^*_k) \right) \right]}.
    \end{aligned}
\end{equation*}
The last inequality holds as the P{\L}K condition implies the quadratic growth condition \citep{karimi2016linear}. In our analysis, this weaker condition leads to a suboptimal characterization of the P{\L}K constant with an exponential dependence on the time horizon $T$ (see discussions in Appendix~\ref{subsection: KL with weak assumptions}). To remove such an exponential dependence, we instead use the stronger sequential decomposition inequality. Interestingly, we demonstrate that sequential decomposition inequalities indeed hold for the control and operations models analyzed in the following sections.

\textbf{Proof Sketch:} To understand the high-level idea of Theorem~\ref{theorem: main result}, we illustrate the role of each condition. First, under the continuous differentiability of expected Q-value functions, Deterministic Policy Gradient Theorem \citep{silver2014deterministic} provides the expression of $\nabla l(\theta)$:
\begin{equation}
\label{eq: deterministic pg}
    \nabla_{\theta_t} l(\theta) = \nabla_{\theta_t} \E_{s_t \sim \rho_t(\cdot|\pi_\theta)} \left[ Q^{\pi_\theta}_t \left(s_t, \pi_t(s_t|\theta_t) \right) \right].
\end{equation}
Based on the policy gradient formulation, our goal is to establish the relationship between the suboptimality gap $l(\theta) - l(\theta^*)$ and $\nabla l(\theta)$. However, $l(\theta) = J(\pi_\theta) = \E [ \sum_{t=1}^T C_t ( s_t, \pi_t(s_t|\theta_t) ) | s_1 \sim \rho, \pi_\theta ]$ has a nested formulation and is difficult for one to check the P{\L}K condition directly. Thanks to the Performance Difference Lemma \citep{kakade2002approximately}, we have
\begin{equation*}
    l(\theta) - l(\theta^*) = \sum_{t=1}^T \left( \mathbb{E}_{s_t \sim \rho_t(\cdot|\pi_\theta)} \left[ Q_t^{\pi_{\theta^*}} \left( s_t, \pi_t(s_t|\theta_t) \right)\right] - \mathbb{E}_{s_t \sim \rho_t(\cdot|\pi_\theta)} \left[Q_t^{\pi_{\theta^*}} \left( s_t, \pi_t(s_t|\theta^*_t) \right) \right] \right).
\end{equation*}

The suboptimality gap $l(\theta) - l(\theta^*)$ can be decomposed as the differences of expected optimal Q-value functions under two different single-stage policies at each period $t$. Therefore, we only need to check the P{\L}K condition for expected optimal Q-value functions across different applications. Leveraging the P{\L}K condition of expected optimal Q-value functions, we have
\begin{equation}
\label{ineq: different gradient}
    l(\theta) - l(\theta^*) \le \sum_{t=1}^T \frac{1}{2\mu_Q} \min_{g_t \in \partial \delta_{\Theta_t}(\theta_t)} \left\| \nabla_{\theta_t} \mathbb{E}_{s_t \sim \rho_t(\cdot|\pi_\theta)} \left[ Q_t^{\pi_{\theta^*}} \left( s_t, \pi_t(s_t|\theta_t) \right)\right] + g_t \right\|_2^2.
\end{equation}

The right-hand side of (\ref{ineq: different gradient}) differs from the gradient formulation $\nabla_{\theta_t}l(\theta)$ presented in (\ref{eq: deterministic pg}). To establish the P{\L}K condition of $l(\theta)$, we need to prove a gradient mismatch inequality:
\begin{equation}
\label{bounded gradient mismatch}
    \begin{aligned}
        & \ \sum_{t=1}^T  \min_{g_t \in \partial \delta_{\Theta_t}(\theta_t)} \left\| \nabla_{\theta_t} \E_{s_t \sim \rho_t(\cdot|\pi_\theta)} \left[ Q^{\pi_{\theta^*}}_t \left(s_t, \pi_t(s_t|\theta_t) \right) \right] + g_t \right\|_2^2\\
        \le & \ M \sum_{t=1}^T\min_{g_t \in \partial \delta_{\Theta_t}(\theta_t)}\left\|\nabla_{\theta_t} \E_{s_t \sim \rho_t(\cdot|\pi_\theta)} \left[ Q^{\pi_\theta}_t \left(s_t, \pi_t(s_t|\theta_t) \right) \right] + g_t \right\|^2_2.
    \end{aligned}
\end{equation}

This inequality captures the relationship of Q-value functions' gradients under two policies. To prove it, we note that
\begin{equation*}
    \begin{aligned}
        & \ \left\| \nabla_{\theta_t} \E_{s_t \sim \rho_t(\cdot|\pi_\theta)} \left[ Q^{\pi_{\theta^*}}_t \left(s_t, \pi_t(s_t|\theta_t) \right) \right] + g_t \right\|_2 - \left\|\nabla_{\theta_t} \E_{s_t \sim \rho_t(\cdot|\pi_\theta)} \left[ Q^{\pi_\theta}_t \left(s_t, \pi_t(s_t|\theta_t) \right) \right] + g_t \right\|_2\\
        \le & \ \left\| \nabla_{\theta_t} \E_{s_t \sim \rho_t(\cdot|\pi_\theta)} \left[ Q^{\pi_{\theta^*}}_t \left(s_t, \pi_t(s_t|\theta_t) \right) \right] - \nabla_{\theta_t} \E_{s_t \sim \rho_t(\cdot|\pi_\theta)} \left[ Q^{\pi_\theta}_t \left(s_t, \pi_t(s_t|\theta_t) \right) \right] \right\|_2\\
        = & \ \left\| \nabla_{\theta_t} l(\theta_1, \dots, \theta_t, \theta_{t+1}^*, \dots, \theta_T^*) - \nabla_{\theta_t} l \left( \theta_1, \dots, \theta_t, \theta_{t+1}, \dots, \theta_T \right) \right\|_2.
    \end{aligned}
\end{equation*}

Applying the sequential decomposition inequality and the P{\L}K condition of expected optimal Q-value functions, we end up with the following inequalities:
\begin{equation*}
    \begin{aligned}
        & \ \left\| \nabla_{\theta_t} l(\theta_1, \dots, \theta_t, \theta_{t+1}^*, \dots, \theta_T^*) - \nabla_{\theta_t} l \left( \theta_1, \dots, \theta_t, \theta_{t+1}, \dots, \theta_T \right) \right\|_2\\
        \le & \ \sum_{k=t+1}^T \left\| \nabla_{\theta_t} l(\theta_1, \dots, \theta_k, \theta_{k+1}^*, \dots, \theta_T^*) - \nabla_{\theta_t} l \left( \theta_1, \dots, \theta_{k-1}, \theta_k^*, \dots, \theta_T^* \right) \right\|_2\\
        \le & \ \sum_{k=t+1}^T M_g \left( \mathbb{E}_{s_k \sim \rho_k(\cdot|\pi_\theta)} \left[ Q_k^{\pi_{\theta^*}} \left( s_k, \pi_k(s_k|\theta_k) \right)\right] - \mathbb{E}_{s_k \sim \rho_k(\cdot|\pi_\theta)} \left[Q_k^{\pi_{\theta^*}} \left( s_k, \pi_k(s_k|\theta^*_k) \right) \right] \right)\\
        \le & \ \sum_{k=t+1}^T \frac{M_g}{2\mu_Q} \min_{g_k \in \partial \delta_{\Theta_k}(\theta_k)}\left\|\nabla_{\theta_k} \E_{s_k \sim \rho_k(\cdot|\pi_\theta)} \left[ Q^{\pi_{\theta^*}}_k \left(s_k, \pi_k(s_k|\theta_k) \right) \right] + g_k \right\|^2_2.
    \end{aligned}
\end{equation*}

The following technical lemma is essential to proving the gradient mismatch inequality.
\begin{lemma}
    \label{lemma: sequence}
    Assume that the nonnegative sequences $\{X_t\}_{t=1}^T$ and $\{Y_t\}_{t=1}^T$ satisfy
    \begin{equation}
        \label{ineq: bounded sequential gradient difference}
        |X_t - Y_t| \le M_g \sum_{k=t+1}^T X_k^2,
    \end{equation}
    with some positive constant $M_g$. If $X_T = Y_T$ and $X_t, Y_t \le G$ for all $t=1,\dots,T$, then
    \begin{equation*}
        \sum_{t=1}^TX_t^2 \le \max\{e, 4e M_g^2 G^2 T^2\} \sum_{t=1}^T Y_t^2.
    \end{equation*}
\end{lemma}

Without loss of generality, we focus on the case with $e<4eM_g^2G^2T^2$. Employing Lemma~\ref{lemma: sequence}, we prove the gradient mismatch inequality (\ref{bounded gradient mismatch}), thereby demonstrating how (\ref{eq: deterministic pg}) relates to (\ref{ineq: different gradient}). This completes the proof of the P{\L}K condition. We refer readers to Appendix~\ref{appendix: KL in PG} for a more rigorous proof.

Theorem~\ref{theorem: main result} allows us to verify the P{\L}K condition for policy gradient optimization problems by verifying several structural conditions. In the following sections, we demonstrate the P{\L}K condition using Theorem~\ref{theorem: main result} for several control and operations models, e.g., the entropy-regularized tabular MDPs, the LQR problem, the inventory systems with Markov-modulated demand, and the stochastic cash balance problem. Leveraging the P{\L}K condition, we demonstrate the global convergence of exact and stochastic policy gradient methods for solving these problems and provide their sample complexities to achieve $\epsilon$-optimal solutions.

\section{Entropy-Regularized Tabular MDPs}
\label{section: tabular MDPs}
Tabular MDP is one of the most popular models in Reinforcement Learning (RL) and many papers focus on this setting \citep{agarwal2021theory, lan2023policy, klein2024beyond}. This section considers a finite-horizon version of this problem with an entropy regularization to the per-period cost functions, a variant of infinite-horizon regularized tabular MDPs analyzed in \citet{bhandari2024global}. The regularization smooths the objective function such that optimal policies are stochastic. See \citet{geist2019theory} for a more detailed discussion about the regularized MDPs.

\subsection{Problem Formulation}
Consider an MDP $\cm = (\cs,\ca,P,C,T,\rho)$ with a finite state space $\cs = \{1, \dots, m\}$. We assume a finite set of actions $\cn = \{1,\dots,n\}$ to choose and take $\ca = \Delta(\cn)$ as the set of probability distributions over these actions. For the tabular setting, it is natural to work directly with a randomized policy $\pi_t(s_t|\theta_t) = \theta_{t}(s_t, \cdot)\in\ca$, instead of using any parameterization as the functional optimization reduces to a finite-dimensional optimization problem. Therefore, the per-period costs and transition functions of $\pi_\theta$ are
\begin{equation*}
    C_t \left( s_t, \pi_t(s_t | \theta_t) \right) = \sum_{i \in \cn} \theta_{t}(s_t, i) C_t (s_t, i), \quad P_t \left( s_{t+1}|s_t, \pi_t(s_t | \theta_t) \right) =  \sum_{i \in \cn} \theta_{t}(s_t, i) P_t (s_{t+1}|s_t, i).
\end{equation*}
We assume that per-period costs are nonnegative and uniformly bounded by $\bar{C}$ for all $s_t\in\cs$, $i\in\cn$, and $t\in[T]$. We further assume a uniformly positive transition kernel, i.e., $P_t(s_{t+1}\mid s_t,i)\ge \underline{p}>0$ for all $s_t,s_{t+1}\in\cs$, $i\in\cn$, and $t\in[T]$. Finally, we assume the initial-state distribution satisfies $\rho(s_1)\ge \underline{\rho}>0$ for all $s_1\in\cs$. Together, these assumptions ensure that every state is reachable with positive probability.

Define $\mathcal{R}(p) \coloneqq D_{\text{KL}}(U||p) = \sum_{i=1}^n\frac{1}{n}\log(\frac{1/n}{p_i})$ as the Kullback-Leibler (KL) divergence between a uniform distribution and $p\in\Delta(\cn)$. $\mathcal{R}(p)$ is a strongly convex function with $\text{dom}(\mathcal{R}) = \{p\in\Delta(\ca):\mathcal{R}(p) < +\infty\} = \{p \in \Delta(\ca): \min_{i} p_i > 0 \}$. The entropy-regularized per-period cost of $\pi$ is
\begin{equation*}
    C^r_t\left( s_t, \pi_t(s_t|\theta_t) \right) = C_t \left( s_t, \pi_t(s_t|\theta_t) \right) + \lambda\mathcal{R}\left(\pi_t(s_t|\theta_t)\right).
\end{equation*}

Since $\theta_t(s_t,\cdot)$ is a probability vector, it must satisfy the constraints $\sum_{i=1}^n \theta_t(s_t,i)=1$ and $\theta_t(s_t,i)\ge 0$ for all $s_t\in\cs$ and $i\in\cn$. With the regularization term, the policy has strictly positive components, i.e., $\theta_t(s_t,i)>0$ for all $s_t\in\cs$ and $i\in\cn$. Moreover, by analyzing properties of the optimal solution $\theta_t^*$, we can further restrict the feasible set to $\Theta_t = \{\theta_t \in\R^{m\times n}: \sum_{i=1}^n \theta_t(s_t, i) = 1, \theta_t(s_t,i) \ge \underline{\theta}, \forall s_t\in\cs, \forall i \in \cn\}$, where $\underline{\theta}\coloneqq \lambda/(n\bar{C}T+n\lambda)$. The detailed argument is deferred to Appendix~\ref{subsection: feasible region}.

\subsection{P{\L}K condition of Policy Gradient Objectives}\label{subsection: tabular MDPs PLK}
Let $\pi^*$ denote the optimal policy that minimizes the total expected cost and $Q^*$ denote the corresponding Q-value function. For any $\pi_\theta\in\Pi_\Theta$, the expected Q-value function satisfies the Bellman equation (\ref{bellman equation}):
\begin{equation*}
    \begin{aligned}
        & \ \E_{s_t \sim \rho_t(\cdot|\pi_\theta)} \left[ Q^{\pi_\theta}_t \left(s_t, \pi_t(s_t|\theta_t) \right) \right]\\
        = & \ \sum_{s_t \in \cs} \rho_t(s_t|\pi_\theta) \Bigl( C_t \left( s_t, \pi_t(s_t|\theta_t) \right) + \lambda\mathcal{R} \left( \pi_t(s_t|\theta_t) \right) + \sum_{s_{t+1}\in\cs}P_{t}\left(s_{t+1}|s_t,\pi_t(s_t|\theta_t)\right) V_{t+1}^{\pi_\theta}(s_{t+1}) \Bigr)\\
        = & \ \sum_{s_t\in\cs} \rho_t(s_t|\pi_\theta) \biggl( \underbrace{\lambda\mathcal{R} \left( \theta_t(s_t, \cdot) \right)}_{(\text{I})} + \underbrace{\sum_{i\in\cn} \theta_t(s_t, i) \Bigl( C_t(s_t, i) + \sum_{s_{t+1}\in\cs}P_{t} (s_{t+1}|s_t, i) V_{t+1}^{\pi_\theta}(s_{t+1}) \Bigr)}_{(\text{II})} \biggr).
    \end{aligned}
\end{equation*}
Term (I) is strongly convex in $\theta_t$ \citep{bhandari2024global} and (II) is linear in $\theta_t$. The continuous differentiability condition holds as the linear function and entropy regularization are smooth. If we replace $V_{t+1}^{\pi_\theta}$ with $V_{t+1}^*$ in (II), we get the expression of expected optimal Q-value function at period $t$, which is strongly convex in $\theta_t$. From Corollary~\ref{corollary: strong convexity implies KL}, expected optimal Q-value functions satisfy the P{\L}K condition.

\begin{lemma}
\label{lemma: tabular MDPs single stage KL}
    The expected Q-value function $\E_{s_t \sim \rho_t(\cdot|\pi_\theta)} \left[ Q^{\pi_\theta}_t \left(s_t, \pi_t(s_t|\theta_t) \right) \right]$ is continuously differentiable on $\Theta_t$. Furthermore, the expected optimal Q-value function $\E_{s_t \sim \rho_t(\cdot|\pi_\theta)} \left[ Q^*_t \left(s_t, \pi_t(s_t|\theta_t) \right) \right]$ satisfies the P{\L}K condition with constant $\lambda\min\{\underline{\rho}, \underline{p}\} / n$ over $\Theta_t$.
\end{lemma}

To verify the bounded gradients and sequential decomposition inequality, we apply the Policy Gradient Theorem \citep[Theorem 1]{sutton1999policy} and get
\begin{equation*}
    \begin{aligned}
        \nabla_{\theta_t(s_t, i)} l(\theta) &= \nabla_{\theta_t(s_t, i)} \E_{s_t \sim \rho_t(\cdot|\pi_\theta)} \left[ Q^{\pi_\theta}_t \left(s_t, \pi_t(s_t|\theta_t) \right) \right]\\
        &= \rho_t(s_t|\pi_\theta) \Bigl( \underbrace{-\frac{\lambda}{n\theta_t(s_t, i)}}_{\text{(I)}} + \underbrace{C_t(s_t, i)}_{\text{(II)}} + \sum_{s_{t+1}\in\cs}P_{t}(s_{t+1}|s_t, i) V_{t+1}^{\pi_\theta}(s_{t+1}) \Bigr).
    \end{aligned}
\end{equation*}

The absolute value of (I) is upper bounded by $\lambda / (n\underline{\theta})$. (II) is uniformly bounded by $\bar{C}$. The regularized per-period costs are bounded as well. Therefore, $l(\theta)$ has bounded gradients by mathematical induction.

\begin{lemma}
\label{lemma: tabular MDPs bounded gradient}
    The policy gradient objective function $l(\theta)$ has bounded gradients
    \begin{equation*}
        \left\| \nabla_{\theta_t}l(\theta) \right\|_F \le \sqrt{n} \bigl[ 2T \bar{C} + \lambda + \lambda T \log((T\bar{C}+\lambda) / \lambda) \bigr], \quad \forall t\in[T].
    \end{equation*}
\end{lemma}

According to Theorem~\ref{theorem: main result}, the last condition we need to verify is the sequential decomposition inequality. This inequality characterizes how far the gradients under the two policies
\begin{equation*}
    \pi_\alpha \coloneqq \left( \pi_1(\cdot|\theta_1), \dots, \pi_{k-1}(\cdot|\theta_{k-1}), \pi_k(\cdot|\theta_k), \pi_{k+1}(\cdot|\theta_{k+1}^*), \dots, \pi_T(\cdot|\theta_T^*) \right)
\end{equation*}
and 
\begin{equation*}
    \pi_\beta \coloneqq \left( \pi_1(\cdot|\theta_1), \dots, \pi_{k-1}(\cdot|\theta_{k-1}), \pi_k(\cdot|\theta_k^*), \pi_{k+1}(\cdot|\theta_{k+1}^*), \dots, \pi_T(\cdot|\theta_T^*) \right).
\end{equation*}

The structure of tabular MDPs naturally builds the connection between the difference in gradients and the difference in optimal Q-value functions for $t < k$:
\begin{equation*}
    \begin{aligned}
        & \ \left\| \nabla_{\theta_t} l(\theta_1, \dots, \theta_{k-1}, \theta_k, \theta_{k+1}^*, \dots, \theta_T^*) - \nabla_{\theta_t} l( \theta_1, \dots, \theta_{k-1}, \theta_k^*, \theta_{k+1}^*, \dots, \theta_T^* ) \right\|_F\\
        \le & \ \sum_{s_t\in\cs, i_t\in\cn} \left| \nabla_{\theta_t(s_t, i_t)} l(\theta_1, \dots, \theta_{k-1}, \theta_k, \theta_{k+1}^*, \dots, \theta_T^*) - \nabla_{\theta_t(s_t, i_t)} l( \theta_1, \dots, \theta_{k-1}, \theta_k^*, \theta_{k+1}^*, \dots, \theta_T^* ) \right|\\
        = & \ \sum_{s_t \in \cs, i_t \in \cn} \left| \rho_t(s_t | \pi_\theta) \sum_{s_{t+1} \in \cs} P_t(s_{t+1} | s_t, i_t) \left( Q_{t+1}^{\pi_\alpha} \left(s_{t+1}, \pi_{t+1}(s_{t+1}|\theta_{t+1}) \right) - Q_{t+1}^{\pi_\beta} \left(s_{t+1}, \pi_{t+1}(s_{t+1}|\theta_{t+1}) \right) \right) \right|\\
        & \dots\\
        = & \ \sum_{s_t \in \cs, i_t\in\cn} \Biggl| \rho_t(s_t | \pi_\theta) \sum_{s_{t+1} \in \cs} P_t(s_{t+1} | s_t, i_t) \sum_{i_{t+1}\in\cn} \theta_{t+1}(s_{t+1}, i_{t+1}) \sum_{s_{t+2} \in \cs} P_{t+1}(s_{t+2} | s_{t+1}, i_{t+1}) \dots\\
        & \quad \quad \quad  \sum_{i_{k-1}\in\cn} \theta_{k-1}(s_{k-1}, i_{k-1}) \sum_{s_k\in\cs} P_{k-1}(s_k|s_{k-1}, i_{k-1}) \underbrace{ \left( Q_k^* \left( s_k, \pi_k(s_k|\theta_k) \right) - Q_k^* \left( s_k, \pi_k(s_k|\theta_k^*) \right) \right)}_{\text{(I)}} \Biggr|.
    \end{aligned}
\end{equation*}

Term (I) is exactly the difference in optimal Q-value functions at period $k > t$. From the assumption that $\theta_t(s, i) \ge \underline{\theta}$ for any $s\in\cs$, $i\in\cn$, and $t\in[T]$, the sequential decomposition inequality holds. The detailed proof is deferred to Appendix~\ref{subsection: tabular MDPs sequential decomposition inequality}.

\begin{lemma}
\label{lemma: tabular MDPs sequential decomposition}
    Sequential decomposition inequalities hold with $M_g = 1/\underline{\theta}$.
\end{lemma}

We have checked all the required conditions in Theorem~\ref{theorem: main result}. The following theorem establishes the P{\L}K condition of the policy gradient optimization problem for entropy-regularized tabular MDPs.
\begin{theorem}
\label{theorem: tabular MDPs multistage KL}
    Consider the entropy-regularized tabular MDPs. The policy gradient optimization problem satisfies the P{\L}K condition with the corresponding P{\L}K constant $\mu_l = \frac{\lambda^3 \min\{\underline{\rho}, \underline{p}\}^3 \underline{\theta}^2}{en^4T^2 \left[ 2T \bar{C} + \lambda + \lambda T \log((T\bar{C}+\lambda) / \lambda) \right]^2}$:
    \begin{equation*}
        l(\theta) - l(\theta^*) \le \frac{1}{2\mu_l} \min_{g \in \partial \delta_{\Theta}(\theta)} \left\|\nabla l(\theta) + g \right\|^2_F, \quad \forall \pi\in\Pi.
    \end{equation*}
\end{theorem}
\proof{Proof of Theorem~\ref{theorem: tabular MDPs multistage KL}} Plugging Lemma~\ref{lemma: tabular MDPs single stage KL}, \ref{lemma: tabular MDPs bounded gradient}, and \ref{lemma: tabular MDPs sequential decomposition} into Theorem~\ref{theorem: main result} yields the result. \Halmos

Leveraging the P{\L}K condition, one can establish a linear convergence rate for exact policy gradient methods and $\tilde{\co}(\epsilon^{-1})$ sample complexity for stochastic policy gradient methods to achieve an optimal policy by Lemma~\ref{lemma: PGD KL}. This is essentially the same as the result in \citet{bhandari2024global}, which demonstrated the gradient domination of $l(\theta)$, implying a linear convergence rate for exact policy gradient methods.

\section{Linear Quadratic Regulator}
\label{section: LQR}
The Linear Quadratic Regulator (LQR) is one of the fundamental problems in the optimal control theory. It seeks an optimal control for a linear dynamic system, in which the state's dynamic is a linear function of the current state and action while incurring a quadratic cost. We present the problem following most of the terminologies in \citet{fazel2018global} and \citet{hambly2021policy}, while keeping some differences to maintain consistency in Section~\ref{section: problem formulation}.

\subsection{Problem Formulation}
Consider an MDP with $\cm=(\cs, \ca, P, C, T, \rho)$. We aim to solve the following optimization problem over a finite time horizon $T$:
\begin{equation}
\label{optimization: LQR objective}
    \min_{ \{a_t\}_{t=0}^{T-1} } \E \left[ \sum_{t=0}^{T-1} (s_t^\top Q_t s_t + a_t^\top R_t a_t) + s_T^\top Q_T s_T \Biggl| s_0 \sim \rho\right],
\end{equation}
such that for all $t=0,\dots,T-1$,
\begin{equation}
\label{equation: LQR transition kernel}
    s_{t+1} = A s_t + B a_t + w_t.
\end{equation}

Here $s_t\in\R^m$ is the state of the system with an initial distribution $\rho$, $a_t\in\R^n$ is the action at period $t$, and $\{w_t\}_{t=0}^{T-1}$ are independent and identical distributed random variables with zero mean that are independent from the initial distribution. Dynamic (\ref{equation: LQR transition kernel}) captures the transition kernel $P_t$ and the cost function is
\begin{equation*}
    C_t(s_t, a_t) = s_t^\top Q_t s_t + a_t^\top R_t a_t, \quad \forall t=0,\dots,T-1.
\end{equation*}
Our analysis can deal with time-dependent parameters in (\ref{equation: LQR transition kernel}), i.e., $A_t$ and $B_t$. For simplicity, we assume that these parameters are time-independent. To ensure the problem is well-defined, we make the following assumptions:
\begin{assumption}
\label{assumption: LQR}
    Assume that the following assumptions hold.
    \begin{henumerate}
        \item \textbf{(Cost Parameters)} Assume that $Q_t \in \R^{m \times m}$ and $R_t \in \R^{n \times n}$ are positive definite matrices for all $t=0,\dots,T-1$. Furthermore, define $\underline{\sigma}_Q$ and $\underline{\sigma}_R$ as the smallest eigenvalue of $\{Q_t\}_{t=0}^{T-1}$ and $\{R_t\}_{t=0}^{T-1}$ respectively:
        \begin{equation*}
            \begin{aligned}
                \underline{\sigma}_Q &= \min_{t=0,\dots,T-1} \sigma_{\text{min}}(Q_t) > 0,\\
                \underline{\sigma}_R &= \min_{t=0,\dots,T-1} \sigma_{\text{min}}(R_t) > 0.
            \end{aligned}
        \end{equation*}
        \item \textbf{(Randomness)} Assume that the second moments of $s_0$ and $\{w_t\}_{t=0}^{T-1}$ are finite. Furthermore, we assume that $\E[s_0s_0^\top]$ and $\E[w_tw_t^\top]$ are positive definite matrices for all $t=0,\dots,T-1$.
    \end{henumerate}
\end{assumption}
Similarly, define $\underline{\sigma}_X$ as the smallest eigenvalue of $\E[s_ts_t^\top]$:
\begin{equation*}
    \underline{\sigma}_X = \min_{t=0,\dots,T} \sigma_{\text{min}}\left( \E_{s_t \sim \rho_t(\cdot| \pi_\theta) }[s_ts_t^\top] \right).
\end{equation*}
Then, we have the following result that shows the well-definedness of the state covariance matrix.
\begin{lemma}[\citealp{hambly2021policy}, Lemma 3.2]
\label{lemma: positive definite state covariance matrix}
    Suppose that Assumption~\ref{assumption: LQR} holds, the second moment of $s_t$ is positive definite for any $t=0,\dots, T$ under any policy $\pi\in\Pi$. Therefore, we have $\underline{\sigma}_X > 0$.
\end{lemma}

This lemma is essential for the landscape characterization and the global convergence of policy gradient methods. In this setting, it is well-known that the linear policy $\pi_t(s_t|\theta_t) = \theta_t s_t$ is optimal for some unknown parameters $\theta_t \in \R^{n \times m}$ \citep{bertsekas1995dynamic}. The policy gradient objective function using linear policies is
\begin{equation*}
    l(\theta) = \E \left[ \sum_{t=0}^{T-1} \left( s_t^\top Q_t s_t + (\theta_t s_t)^\top R_t (\theta_t s_t) \right) + s_T^\top Q_T s_T \Biggl| s_0 \sim \rho\right]
\end{equation*}
with $s_{t+1} = (A + B\theta_t) s_t + w_t$. When the linear system is unstable, complexity in \citet{hambly2021policy} has an exponential dependence on $T$. To stabilize the system, one needs to restrict to the set $\{\theta: r_\sigma(A + B \theta_t) \le 1, \ \forall t=0, \dots, T-1\}$. However, this set is unbounded and non-convex, making the analysis difficult. 

In the following, we only consider the landscape of $l(\theta)$ within a convex and compact set $\Theta$ such that $r_\sigma(A+B\theta_t) \le \|A+B\theta_t\|_2 \le 1$ for all $\theta_t \in \Theta_t$ and $0\le t \le T-1$. Furthermore, we denote $\bar{\sigma}_\Theta$ as the largest spectral norm of $\theta_t\in\Theta_t$ for all $0\le t \le T-1$, i.e.,
\begin{equation*}
    \bar{\sigma}_\Theta = \max_{t=0,\dots,T-1} \left\{ \max\{ \|\theta_t\|_2 : \theta_t \in \Theta_t\} \right\}.
\end{equation*}
For the positive definite matrices $Q_t, R_t, \E[s_0s_0^\top], \E[w_t w_t^\top]$, we define $\bar{\sigma}_Q, \bar{\sigma}_R, \bar{\sigma}_X, \bar{\sigma}_W$ as the upper bound on their eigenvalues for all $t=0,\dots,T-1$, respectively. In addition, we assume $r_\sigma(Q_T) \le \bar{\sigma}_Q$ as well.

\subsection{P{\L}K Condition of Policy Gradient Objectives}
To establish the P{\L}K condition of $l(\theta)$, we aim to verify all the conditions in Theorem~\ref{theorem: main result}. The continuous differentiability and the P{\L}K condition of expected optimal Q-value functions come from the definition of Q-value functions. Given any policy $\pi_\theta \in \Pi_\Theta$, Q-value functions satisfy (\ref{bellman equation}):
\begin{equation*}
    \begin{aligned}
        Q_t^{\pi_\theta} \left( s_t, \pi_t(s_t|\theta_t) \right) = \underbrace{s_t^\top Q_t s_t + s_t^\top \theta_t^\top R_t \theta_t s_t}_{(\text{I})} + \underbrace{\E_{w_t} \left[ V_{t+1}^{\pi_\theta} \left( (A + B\theta_t)s_t + w_t \right) \right]}_{(\text{II})}.
    \end{aligned}
\end{equation*}

(I) is a quadratic function with a positive definite matrix $R_t$, and is therefore continuously differentiable. For (II), the value function is continuously differentiable by mathematical induction. Since the composition of a continuously differentiable function and a linear function is continuously differentiable, we have the continuous differentiability of (II). If we plug $\pi^*$ into (II), we get an explicit expression of the optimal Q-value function $Q_t^*$. \citet{bertsekas1995dynamic} demonstrated the convexity of $V_t^*$ by mathematical induction. Therefore, (II) is a convex function of $\theta_t$, which implies the strong convexity (and the P{\L}K condition) of the optimal Q-value function combined with the strong convexity of (I).

\begin{lemma}
\label{lemma: LQR single stage KL}
    Suppose that Assumption~\ref{assumption: LQR} holds. The expected Q-value function is continuously differentiable on $\Theta_t$. Furthermore, the expected optimal Q-value function $\E_{s_t \sim \rho_t(\cdot|\pi_\theta)} \left[ Q^*_t \left(s_t, \pi_t(s_t|\theta_t) \right) \right]$ satisfies the P{\L}K condition on $\Theta_t$ with P{\L}K constant $2\underline{\sigma}_X\underline{\sigma}_R$.
\end{lemma}

To validate other conditions in Theorem~\ref{theorem: main result}, we need an explicit expression of the policy gradient $\nabla l(\theta)$. \citet{hambly2021policy} established the formulation of the Q-value function and the policy gradient $\nabla l(\theta)$ through a recursive form. First, let us define $P_t(\theta)$ by the following recursive equations:
\begin{equation}
\label{eq: LQR P update}
    P_t(\theta) \coloneqq Q_t + \theta_t^\top R_t \theta_t + (A + B \theta_t)^\top P_{t+1}(\theta) (A + B \theta_t), \quad \forall t=0,\dots, T-1,
\end{equation}
The boundary condition is $P_T(\theta) = Q_T$. In addition, define
\begin{equation*}
    L_t(\theta) \coloneqq L_{t+1}(\theta) + \E[w_t^\top P_{t+1}(\theta) w_t], \quad \forall t=0,\dots,T-1, 
\end{equation*}
with $L_T(\theta) = 0$, and
\begin{equation}
\label{eq: LQR E update}
    E_t(\theta) \coloneqq \left( R_t + B^\top P_{t+1}(\theta) B \right) \theta_t + B^\top P_{t+1}(\theta) A, \quad \forall t=0,\dots,T-1.
\end{equation}
The subsequent proposition presents an explicit formulation of the Q-value function and $\nabla l(\theta)$.
\begin{proposition}[\citealp{hambly2021policy}, Lemma~3.5]
\label{proposition: LQR policy gradient}
    The Q-value function has an explicit expression:
    \begin{equation*}
        \E_{s_t\sim\rho_t(\cdot|\pi_\theta)} \left[ Q_t^{\pi_\theta} \left( s_t, \pi_t(s_t|\theta_t) \right) \right] = \E_{s_t\sim\rho_t(\cdot|\pi_\theta)} \left[ s_t^\top P_t(\theta) s_t \right] + L_t(\theta).
    \end{equation*}
    Furthermore, the policy gradient objective function $l(\theta)$ has the following gradient form:
    \begin{equation*}
        \nabla_{\theta_t} l(\theta) = 2E_t(\theta) \E_{s_t\sim\rho_t(\cdot|\pi_\theta)}[s_t s_t^\top].
    \end{equation*}
\end{proposition}

From the expression of the policy gradient $\nabla l(\theta)$, we verify the bounded gradients condition by showing the boundedness of $P_t$. It utilizes the stability of the linear system and the compactness of the feasible region $\Theta_t$. The following lemma establishes a formal result.

\begin{lemma}
\label{lemma: LQR bounded gradient}
    Suppose that Assumption~\ref{assumption: LQR} holds. The policy gradient objective function has bounded gradients, i.e., $\|\nabla_{\theta_t} l(\theta)\|_F \le G$ for any $0\le t \le T-1$. Furthermore, $G$ is polynomial in the model parameters $(m, n, T, \bar{\sigma}_Q, \bar{\sigma}_R, \bar{\sigma}_\Theta, \bar{\sigma}_X, \bar{\sigma}_W, \|B\|_2)$.
\end{lemma}

Sequential decomposition inequalities are more complicated to verify. However, the structure of the LQR problem helps to construct the relationship between the difference in gradients and the difference in optimal Q-value functions. To see this, define $\Pi_{[j_1:j_2]} \coloneqq (A + B\theta_{j_2})(A + B\theta_{j_2 - 1}) \dots (A + B\theta_{j_1 + 1})(A + B\theta_{j_1})$ for $j_1 \le j_2$. Recall the gradient formulation in Proposition~\ref{proposition: LQR policy gradient}, for any $1\le t < k \le T$:
\begin{equation*}
    \begin{aligned}
        & \nabla_{\theta_t} l(\theta_{[0:k-1]}, \theta_k, \theta^*_{[k+1:T-1]}) - \nabla_{\theta_t} l(\theta_{[0:k-1]}, \theta^*_k, \theta^*_{[k+1:T-1]})\\
        = \ & 2B^\top \left( P_{t+1}(\theta_{[0:k-1]}, \theta_k, \theta^*_{[k+1:T-1]}) - P_{t+1}(\theta_{[0:k-1]}, \theta^*_k, \theta^*_{[k+1:T-1]}) \right) (A + B\theta_t) \E_{s_t\sim\rho_t(\cdot|\pi_\theta)}[s_ts_t^\top]\\
        = \ & 2B^\top (A + B\theta_{t+1})^\top \bigl( P_{t+2}(\theta_{[0:k-1]}, \theta_k, \theta^*_{[k+1:T-1]})\\
        & \quad \quad \quad \quad \quad \quad \quad - P_{t+2}(\theta_{[0:k-1]}, \theta^*_k, \theta^*_{[k+1:T-1]}) \bigr) (A + B\theta_{t+1}) (A + B\theta_t) \E_{s_t\sim\rho_t(\cdot|\pi_\theta)}[s_ts_t^\top]\\
        &\dots\\
        = \ & 2B^\top \Pi_{[t+1:k-1]}^\top \left( P_k(\theta_{[0:k-1]}, \theta_k, \theta^*_{[k+1:T-1]}) - P_k(\theta_{[0:k-1]}, \theta^*_k, \theta^*_{[k+1:T-1]}) \right) \Pi_{[t:k-1]} \E_{s_t\sim\rho_t(\cdot|\pi_\theta)}[s_ts_t^\top].
    \end{aligned}
\end{equation*}
The second equation uses the update (\ref{eq: LQR P update}). Next, we proceed to the difference in optimal Q-value functions. Utilizing the explicit expression of the Q-value function in Proposition~\ref{proposition: LQR policy gradient}, we conclude that
\begin{equation*}
    \begin{aligned}
        & \E_{s_k\sim\rho_k(\cdot|\pi_\theta)} \left[ Q_k^* \left( s_k, \pi_k(s_k|\theta_k) \right) \right] - \E_{s_k\sim\rho_k(\cdot|\pi_\theta)} \left[ Q_k^* \left( s_k, \pi_k(s_k|\theta_k^*) \right) \right]\\
        = \ & \E_{s_k\sim\rho_k(\cdot|\pi_\theta)} \left[ s_k^\top \left( P_k(\theta_{[0:k-1]}, \theta_k, \theta^*_{[k+1:T-1]}) - P_k(\theta_{[0:k-1]}, \theta^*_k, \theta^*_{[k+1:T-1]}) \right) s_k \right]\\
        = \ & \text{Tr} \left( \left(P_k(\theta_{[0:k-1]}, \theta_k, \theta^*_{[k+1:T-1]}) - P_k(\theta_{[0:k-1]}, \theta^*_k, \theta^*_{[k+1:T-1]}) \right) \E_{s_k\sim\rho_k(\cdot|\pi_\theta)}[s_ks_k^\top] \right).
    \end{aligned}
\end{equation*}

Both the difference in gradients of $l(\theta)$ and the difference in expected optimal Q-value functions have the component $P_k(\theta_{[0:k-1]}, \theta_k, \theta^*_{[k+1:T-1]}) - P_k(\theta_{[0:k-1]}, \theta^*_k, \theta^*_{[k+1:T-1]})$. Leveraging this, we can prove the sequential decomposition inequality under some mild assumptions.

\begin{lemma} \label{lemma: LQR sequential decomposition}
    Suppose that Assumption~\ref{assumption: LQR} holds. The \textit{Sequential Decomposition} condition holds with $M_g > 0$. Furthermore, $M_g$ is polynomial in the model parameters $(m, n, T, \bar{\sigma}_Q, \bar{\sigma}_R, \bar{\sigma}_\Theta, \bar{\sigma}_X, \bar{\sigma}_W, \|B\|_2)$.
\end{lemma}

With the continuous differentiability condition (Lemma~\ref{lemma: LQR single stage KL}), P{\L}K condition of optimal Q-value functions (Lemma~\ref{lemma: LQR single stage KL}), bounded gradient condition (Lemma~\ref{lemma: LQR bounded gradient}), and sequential decomposition inequality (Lemma~\ref{lemma: LQR sequential decomposition}), we are ready to demonstrate the P{\L}K condition of the policy gradient optimization problem.

\begin{theorem}
\label{theorem: LQR multistage KL}
    Consider the LQR problem. Suppose that Assumption~\ref{assumption: LQR} holds. The policy gradient objective function $l(\theta)$ satisfies the P{\L}K condition on $\Theta$:
    \begin{equation*}
        l(\theta) - l(\theta^*) \le \frac{1}{2\mu_l} \min_{g \in \partial \delta_{\Theta}(\theta)} \sum_{t=0}^{T-1} \|\nabla_{\theta_t} l(\theta) + g_t\|_F^2, \quad \forall \theta \in \Theta,
    \end{equation*}
    where $g = (g_0,\dots,g_{T-1})$. In addition, the reciprocal of P{\L}K constant $\mu_l$ is polynomial in the model parameters $(m, n, T, \bar{\sigma}_Q, \bar{\sigma}_R, \underline{\sigma}_R^{-1}, \bar{\sigma}_\Theta, \bar{\sigma}_X, \underline{\sigma}_X^{-1}, \bar{\sigma}_W, \|B\|_2)$.
\end{theorem}

\proof{Proof of Theorem~\ref{theorem: LQR multistage KL}} Plugging Lemma~\ref{lemma: LQR single stage KL}, \ref{lemma: LQR bounded gradient}, and \ref{lemma: LQR sequential decomposition} into Theorem~\ref{theorem: main result} yields the result. \Halmos

Following the proof of \citet[Theorem~3.3]{hambly2021policy}, we can establish a linear convergence rate of exact policy gradient methods. It is worth noting that the P{\L}K constant in Theorem~\ref{theorem: LQR multistage KL} is different from that in \citet[Lemma~3.6]{hambly2021policy}. They fully explored the special structure of the LQR problem and used an important fact that expected Q-value functions induced by \textit{any} policy $\pi_\theta$ is quadratic. To illustrate the general applicability of our unified framework, we instead utilize the P{\L}K condition of expected \textit{optimal} Q-value functions, as this is a less restrictive condition shared with other operations problems. 

\section{Inventory Models}
\label{section: inventory system}
This section demonstrates how to validate the assumptions in Theorem~\ref{theorem: main result} to establish the P{\L}K condition for the policy gradient optimization problem of the multi-period inventory system with Markov-modulated demand, where unsatisfied demands are backlogged. One can extend the result to the lost sales model and derive a sample complexity with the same order.

Many works in inventory control assume that random demands are independent across time. However, this assumption is often unrealistic in the real world, and the underlying demand process might be correlated, i.e., economic conditions and seasons affect the random demands of different periods. To capture correlations, some literature models the underlying demand process by an exogenous discrete-time, discrete-state Markov chain \citep{song1993inventory}. We briefly state the problem formulation and validate the P{\L}K condition of the policy gradient optimization problem using a state-dependent base-stock policy class.

\subsection{Problem Formulation}

Consider a MDP framework with $\cm = (\cs, \ca, P, C, T, \rho)$. At the beginning of period $t$, the decision-maker observes state $s_t = (x_t, i_t)$ and determines the order quantity $a_t \ge 0$, where $x_t$ represents the current inventory level and $i_t$ is the state of the world. The inventory level $x_1$ follows an initial distribution $\rho$. The replenishment immediately raises the inventory level to $y_t = x_t + a_t \ge x_t$. Subsequently, the decision-maker observes a random demand $D_t$ whose distribution depends on the current state of the world $i_t$. Finally, the inventory level at the beginning of the next period follows the linear transition kernel $x_{t+1} = y_t - D_t = x_t + a_t - D_t$, where the negative inventory level represents backlogged demands.

Suppose the exogenous Markov chain has a finite state space $\ci$. In state $i_t\in\ci$, the random demand $D_t$ follows a cumulative distribution function $P_{D}(\cdot|i_t)$. State $i$ moves to the next state $j$ with probability $p(j|i) \in [0,1]$ and $\sum_{j\in\ci}p(j|i) = 1$ for any $i\in\ci$. As a finite-state time-homogeneous Markov Chain has at least one stationary distribution, let us pick $\nu\in\R^{|\ci|}$ as one of its stationary distributions. Assume that initial state $i_1$ follows the stationary distribution $\nu$, e.g., $i_t \sim \nu$ for any $t\in[T]$.

When the on-hand inventory level exceeds the realized demand $D_t$, it incurs a holding cost of $h_t \ge 0$ per unit. Otherwise, insufficient inventory causes a backlogging cost of $b_t \ge 0$ per unit. Let $L_t$ denote the expected holding and backlogging cost for period $t$, which is a convex function of the order-up-to level $y_t$:
\begin{equation*}
    \begin{aligned}
        L_t(y_t|i_t) = \mathbb{E}_{D_t \sim P_D(\cdot | i_t)}\left[ h_t \left( y_t - D_t \right)^+ + b_t \left( y_t - D_t \right)^- \right].
    \end{aligned}
\end{equation*}

For simplicity, we omit the ordering cost. Our analysis can accommodate positive linear ordering costs, which can be subsumed in holding and backlogging costs. The per-period cost in the MDP framework is $C_t(s_t, a_t) = L_t(x_t + a_t|i_t)$. We aim to identify an ordering policy that minimizes the total expected cost over $T$ periods:
\begin{equation}
\label{optimization: markov modulated demand}
    \min_{\{a_t\}_{t=1}^T} \E \left[ \sum_{t=1}^T L_t(x_t + a_t|i_t) \Biggl| x_1 \sim \rho, i_1 \sim \nu \right],
\end{equation}
where $x_{t+1} = x_t + a_t - D_t$, $D_t \sim P_D(\cdot|i_t)$, and $i_{t+1}\sim p(\cdot|i_t)$. One can use dynamic programming to reformulate (\ref{optimization: markov modulated demand}). Let $V_t^*$ denote the cost-to-go function, which starts at the beginning of period $t$. It satisfies the Bellman optimality equation (\ref{bellman optimality equation}):
\begin{equation*}
    \begin{aligned}
        V_t^*(x_t, i_t) &= \min_{a_t \ge 0} \left\{ L_t \left( x_t + a_t|i_t \right) + \E_{i_{t+1}\sim p(\cdot|i_t), D_t \sim P_D(\cdot | i_t)} \left[ V^*_{t+1} \left( x_t + a_t - D_t, i_{t+1} \right) \right] \right\}\\
        &= \min_{y_t \ge x_t} \left\{ L_t \left( y_t | i_t \right) + \E_{i_{t+1}\sim p(\cdot|i_t), D_t \sim P_D(\cdot | i_t)} \left[ V^*_{t+1} \left( y_t - D_t, i_{t+1} \right) \right] \right\},
    \end{aligned}
\end{equation*}
with $V_{T+1}^*(\cdot) = 0$. Let us define
\begin{equation}
\label{equation: markov-modulated function f}
    f_t(y_t | i_t) \coloneqq L_t \left( y_t | i_t \right) + \E_{i_{t+1}\sim p(\cdot|i_t), D_t \sim P_D(\cdot | i_t)} \left[ V^*_{t+1} \left( y_t - D_t, i_{t+1} \right) \right].
\end{equation}

\citet{song1993inventory} demonstrated the convexity of $V_t^*(x_t, i_t)$ with respect to $x_t$ and the convexity of $f_t(y_t | i_t)$ with respect to $y_t$ for any $i_t\in\ci$ by mathematical induction. Given the convexity of $f_t(y_t | i_t)$, they proved that a state-dependent base-stock policy is optimal. 

Thus, in the following, we focus on the state-dependent base-stock policy class. Specifically, for each state of the exogenous Markov Chain and period $t$, we define the corresponding base-stock level $\theta_{t,i} \in \R$. Let $\theta_t \in \R^{|\ci|}$ to denote the vector of all the base-stock levels in period $t$ with $\theta_{t, i}$ as its $i$-th component. Given the current inventory level $x_t$ and the state of the world $i_t$, the decision-maker orders $\pi_t(x_t, i_t | \theta_t) = ( \theta_{t, i_t} - x_t )^+$. For the optimal base-stock levels $\theta_t^*$, its $i$-th component $\theta^*_{t, i}$ is a minimizer of $f_t(\cdot|i)$.

Given the state-dependent base-stock policy class, we can express the policy gradient objective function of all the base-stock levels $\theta = (\theta_1,\dots,\theta_T) \in \R^{T \times |\ci|}$ by
\begin{equation*}
    l(\theta) = \E \left[ \sum_{t=1}^T L_t \left( x_t \vee \theta_{t, i_t} | i_t \right) \Biggl| x_1\sim\rho, i_1\sim\nu \right],
\end{equation*}
with $x_{t+1} = x_t \vee \theta_{t,i_t} - D_t$, $D_t \sim P_D(\cdot|i_t)$, and $i_{t+1}\sim p(\cdot | i_t)$. Here $\vee$ denotes a component-wise max operator. In the remaining part of this section, we focus on the landscape of $l(\theta)$ over a compact convex set $\Theta = \Theta_1 \times \dots \times \Theta_T$ with $\Theta_t = \{ \theta_t: \theta_{t, i} \in [0, B], \forall i \in \ci\}$. This is reasonable in practice since one can treat $B \ge 0$ as the capacity of the warehouse. The policy class is $\Pi_\Theta \coloneqq \{ \pi_t(x,i|\theta) = x \vee \theta_{t,i} : \theta_{t, i} \in [0, B]\}$. We restrict the initial inventory level $x_1 \in (-\infty, B]$.

\subsection{Nonconvex Landscape}
In this subsection, we demonstrate the P{\L}K condition of the policy gradient objective function $l(\theta)$ by verifying the conditions in Theorem~\ref{theorem: main result}. To establish the P{\L}K condition, we rely on the following assumptions.
\begin{assumption}
    \label{assumption: inventory}
    Assume that all the following assumptions hold.
    \begin{henumerate}
        \item \label{assumption: inventory independence of rv} The initial inventory level $x_1$ is independent of the exogenous Markov chain. In addition, $x_1$ is independent of the random demands. 
        \item \label{assumption: inventory Lipschitz of rv} The initial cumulative distribution function of $x_1$ is $L_\rho$-Lipschitz continuous. Furthermore, the cumulative distribution function $P_D(\cdot|i)$ is $L_D$-Lipschitz continuous for all $i\in\ci$. 
        \item \label{assumption: inventory KL} There exists a positive constant $\alpha$ such that $\rho(0) \ge \alpha$ and $P_D(B|i) \le 1 - \alpha$ for any $i\in\ci$.
        \item \label{assumption: inventory SC} The probability density functions of random demands at each period are uniformly bounded below by $\mu_D > 0$ on $[0, B]$. 
    \end{henumerate}
\end{assumption}

Assumption~\ref{assumption: inventory}.\ref{assumption: inventory independence of rv} is standard in the setting of Markov-modulated demand \citep{song1993inventory, chen2001optimal}. It allows the correlation of random demands across different periods through the exogenous Markov chain, which is less restrictive than the standard independent assumption in the literature \citep{huh2014online, cheung2019sampling}. Assumption~\ref{assumption: inventory}.\ref{assumption: inventory Lipschitz of rv} ensures the continuous differentiability, which is valid for many commonly used distributions, e.g., uniform, exponential, and Erlang distributions. Assumption~\ref{assumption: inventory}.\ref{assumption: inventory KL} provides sufficient conditions to ensure that the base-stock level $\theta_{t,i_t}$ exceeds the on-hand inventory $x_t$ with positive probability, i.e., $\mathbb{P}(\theta_{t,i_t}\ge x_t)>0$. This guarantees that $\theta_{t,i_t}$ contributes nontrivially to $l(\theta)$ and avoids suboptimal stationary points. Assumption~\ref{assumption: inventory}.\ref{assumption: inventory SC} leads to strongly convex costs $L_t(\cdot|i_t)$ and holds for numerous distributions, such as uniform, exponential, and Erlang distributions. In Appendix~\ref{Appendix: Robustness Check}, we conduct additional numerical experiments to show the robust performance of policy gradient methods even when Assumption~\ref{assumption: inventory}.\ref{assumption: inventory Lipschitz of rv} and Assumption~\ref{assumption: inventory}.\ref{assumption: inventory KL} are violated.

For any policy $\pi_\theta \in \Pi_\Theta$, the expected Q-value functions satisfy the Bellman equation (\ref{bellman equation}):
\begin{equation}
\label{equation: inventory expected Q value function}
    \begin{aligned}
        & \ \E_{(x_t, i_t)\sim \rho_t(\cdot|\pi_\theta)} \left[ Q^{\pi_\theta}_t \left(x_t, i_t, \pi_t(x_t, i_t|\theta_t) \right) \right]\\
        = & \ \E_{(x_t, i_t)\sim \rho_t(\cdot|\pi_\theta)} \left[ L_t(x_t \vee \theta_{t, i_t} | i_t) + \E_{i_{t+1}\sim p(\cdot|i_t), D_t \sim P_D(\cdot | i_t)} \left[ V^{\pi_\theta}_{t+1} \left( x_t \vee \theta_{t, i_t} - D_t, i_{t+1} \right) \right] \right].
    \end{aligned}
\end{equation}
We express the expected Q-value function by $\E_{(x_t, i_t)\sim \rho_t(\cdot|\pi_\theta)} \left[ h(x_t \vee \theta_{t, i_t} | i_t) \right]$, where
\begin{equation*}
    h(y_{t, i_t} | i_t) = L_t(y | i_t) + \E_{i_{t+1}\sim p(\cdot|i_t), D_t \sim P_D(\cdot | i_t)} \left[ V^{\pi_\theta}_{t+1} \left( y - D_t, i_{t+1} \right) \right].
\end{equation*}
Clearly, $h(\cdot|i_t)$ is continuously differentiable by backward mathematical induction. The cumulative distribution function $\rho_t(\cdot|\pi_\theta)$ is continuous. Thus, the expected Q-value function is also continuously differentiable. Furthermore, replacing $\pi_\theta$ with ${\pi_{\theta^*}}$ in (\ref{equation: inventory expected Q value function}) gives the expression of expected optimal Q-value functions:
\begin{equation}
\label{inventory optimal Q}
    F_t(\theta_t) \coloneqq \E_{(x_t, i_t)\sim \rho_t(\cdot|\pi_\theta)} \left[ Q^{\pi_{\theta^*}}_t \left(x_t, i_t, \pi_t(x_t, i_t|\theta_t) \right) \right] = \E_{(x_t, i_t)\sim \rho_t(\cdot|\pi_\theta)} \left[ f_t(x_t \vee \theta_{t, i_t} | i_t) \right].
\end{equation}

It is not straightforward to verify the P{\L}K condition of $F_t(\theta_t)$, and we present a proof sketch to understand why this condition holds. From (\ref{inventory optimal Q}), the suboptimality gap of $F_t(\cdot)$ can be upper bounded using the suboptimality gap of $f_t(\cdot|i_t)$. Since $f_t(\cdot|i_t)$ is strongly convex, its suboptimality gap is dominated by its gradient norm. Applying (\ref{inventory optimal Q}) again, we establish the connection between the gradients of $F_t(\cdot)$ and $f_t(\cdot|i_t)$ and prove the P{\L}K condition of $F_t(\theta_t)$. For more rigorous proof, we refer readers to Appendix~\ref{appendix: inventory models}.

\begin{lemma}
\label{lemma: markov modulated demand single stage KL}
    Suppose that Assumption \ref{assumption: inventory} holds. The expected Q-value function is continuously differentiable on $\Theta_t$. Furthermore, the expected optimal Q-value function $F_t(\theta_t)$ satisfies the P{\L}K condition on $\Theta_t$ with P{\L}K constant $\mu_Q = \min_{t\in[T]}\{h_t+b_t\}\mu_D\alpha^2\min_{i\in\ci}\{\nu_i\}$.
\end{lemma}

Similar to previous applications, verifying remaining conditions requires a gradient expression of the policy gradient objective function. Let $V_t^{\pi_\theta}$ denote the value function, which starts at the beginning of period $t$ and follows policy $\pi_\theta$. It satisfies the Bellman equation (\ref{bellman equation}):
\begin{equation*}
    \begin{aligned}
        V_t^{\pi_\theta}(x_t, i_t) &= L_t \left( x_t \vee \theta_{t, i_t} |i_t \right) + \E_{i_{t+1}\sim p(\cdot|i_t), D_t \sim P_D(\cdot | i_t)} \left[ V^{\pi_\theta}_{t+1} \left( x_t \vee \theta_{t, i_t} - D_t, i_{t+1} \right) \right],
    \end{aligned}
\end{equation*}
with $V_{T+1}^{\pi_\theta}(\cdot) = 0$. The subsequent proposition presents the result in a recursive form.

\begin{proposition}[Policy Gradient Expression]
\label{proposition: Markov modulated demand gradient formulation}
    For any $\theta\in \Theta$ and $t\in[T]$, the partial derivatives of value functions satisfy the following recursive form:
    \begin{equation}
    \label{markov modulated demand: recursive gradient}
        \nabla_{x} V_t^{\pi_\theta} (x_t, i_t) = \mathbf{1} \left( x_t \ge \theta_{t,i_t} \right) \times \Bigl( L'_t(x_t | i_t) + \sum_{i_{t+1}\in\ci} p(i_{t+1} | i_t) \E_{D_t\sim P_D(\cdot|i_t)} \left[ \nabla_x V^{\pi_\theta}_{t+1} (x_t - D_t, i_{t+1}) \right] \Bigr),
    \end{equation}
    where $\nabla_x V_{T+1}^{\pi_\theta}(\cdot, \cdot) = 0$ and $\nabla_x V_t^{\pi_\theta}(x_t, i_t)$ represent the partial derivative of $x_t$. Additionally, the policy gradient objective function $l(\theta)$ has the following gradient form for any $t\in[T]$ and $ i\in\ci$:
    \begin{equation*}
        \begin{aligned}
            & \frac{\partial}{\partial \theta_{t,i}} l(\theta)\\
            = &\E_{(x_t,i_t) \sim \rho_t(\cdot|\pi_\theta)} \biggl[ \mathbf{1} \left( i_t = i, \theta_{t,i} \ge x_t \right) \times \Bigl( L'_t \left( \theta_{t, i} | i \right) + \sum_{i_{t+1}\in\ci} p(i_{t+1}|i) \E_{D_t \sim P_D(\cdot|i)} \left[ \nabla_x V_{t+1}^{\pi_\theta} \left( \theta_{t,i} - D_t, i_{t+1} \right) \right] \Bigr) \biggr].
        \end{aligned}
    \end{equation*}
\end{proposition}

Leveraging Proposition~\ref{proposition: Markov modulated demand gradient formulation}, we can verify the smoothness (Assumption~\ref{assumption: smoothness}). We leave the detailed analysis in Appendix~\ref{subsection: sample complexity}. The following two lemmas verify the bounded gradient condition and the sequential decomposition inequality.

\begin{lemma}
\label{lemma: markov modulated demand bounded gradient}
    Suppose that Assumption \ref{assumption: inventory} holds. The policy gradient objective function $l(\theta)$ has bounded gradients for any $\theta\in\Theta$ and $t\in[T]$:
    \begin{equation*}
        \left\| \nabla_{\theta_t} l(\theta) \right\|_2 \le \max_{t\in[T]} \left\{ \max\{h_t, b_t\} \right\} T.
    \end{equation*}
\end{lemma}

\begin{lemma}
\label{lemma: markov modulated demand bounded sequence gradient difference}
    Suppose that Assumption~\ref{assumption: inventory} holds. The sequential decomposition inequality holds for any $\theta\in\Theta$ and $1\le t < k \le T$, i.e., 
    \begin{equation*}
        \begin{aligned}
            \left\| \nabla_{\theta_t} l(\theta_{[1:k-1]}, \theta_{k}, \theta_{[k+1:T]}^*) - \nabla_{\theta_t} l(\theta_{[1:k-1]}, \theta_{k}^*, \theta_{[k+1:T]}^*) \right\| \le \frac{L_D}{\alpha} \left( F_k(\theta_k) - F_k(\theta_k^*) \right).
        \end{aligned}
    \end{equation*}
\end{lemma}

\begin{theorem}
\label{theorem: markov modulated demand multistage KL}
    Suppose that Assumption~\ref{assumption: inventory} holds. The policy gradient objective function $l(\theta)$ of the multi-period inventory system with Markov-modulated demand satisfies the P{\L}K condition on $\Theta$ with P{\L}K constant
    \begin{equation*}
        \mu_l = \frac{\min_{t\in[T]}\{h_t+b_t\}^3\mu_D^3 \alpha^8 \min_{i\in\ci}\{\nu_i\}^3}{e L_D^2 \max_{t\in[T]} \left\{ \max\{h_t, b_t\} \right\}^2 T^4}.
    \end{equation*}
    More specifically, we have
    \begin{equation*}
        l(\theta) - l(\theta^*) \le \frac{1}{2\mu_l} \min_{g \in \partial \delta_{\Theta}(\theta)} \|\nabla l(\theta) + g\|_2^2, \quad \forall \theta\in\Theta.
    \end{equation*}
\end{theorem}

\proof{Proof of Theorem~\ref{theorem: markov modulated demand multistage KL}} Plugging Lemmas~\ref{lemma: markov modulated demand single stage KL}, \ref{lemma: markov modulated demand bounded gradient}, and \ref{lemma: markov modulated demand bounded sequence gradient difference} into Theorem~\ref{theorem: main result} yields the result. \Halmos

\begin{remark}
    The P{\L}K constant in Theorem~\ref{theorem: markov modulated demand multistage KL} has a dependence on $\min_{i\in\ci}\{\nu_i\}$. Since $\sum_{i\in\ci} \nu_i = 1$, the smallest probability is at most the order of $1/|\ci|$.
\end{remark}

Leveraging the P{\L}K condition, we establish the global convergence of stochastic policy gradient methods by Lemma~\ref{lemma: PGD KL}. The sample complexity required for achieving an $\epsilon$-optimal state-dependent base-stock policy is $\tilde{\co} \left( \epsilon^{-1} \text{poly}(T) \right)$. This is the first sample complexity result for the multi-period inventory system with Markov-modulated demand in the literature. As a byproduct, we improve the sample complexity for stochastic gradient methods to solve the inventory system with independent demands, which is a special case with no exogenous Markov chain ($|\ci| = 1$). Our sample complexity admits a polynomial dependence on the time horizon, representing a significant improvement compared to the exponential dependence in \citet{huh2014online} for a biased stochastic gradient method. We remark that \citet{huh2014online} assumes the convexity of cost-go-functions, whereas Theorem~\ref{theorem: main result} assumes the P{\L}K Condition of expected optimal Q-value functions (which is a relaxation of strong convexity). While the analysis in \citet{huh2014online} can be extended to the case with strongly convex cost-to-go functions, it remains unclear whether the exponential dependence on $T$ can be improved.

\section{Stochastic Cash Balance Problem}\label{section: cash balance}
The cash balance problem originally refers to a cost minimization problem when a firm has to decide how much cash to hold to meet the transaction requirements over a finite planning horizon. It can also model inventory management with rented equipment \citep{whisler1967stochastic, chen2009new}. In this section, we will briefly describe the problem formulation under the inventory setting and validate the P{\L}K condition of the objectives for policy gradient methods using a two-sided base-stock policy class.

\subsection{Problem Formulation}
The problem formulation of the stochastic cash balance problem is similar to the multi-period inventory system. We use the same notation in Section~\ref{section: inventory system}. Unlike the classic inventory system where the decision-maker can only raise the inventory level, the stochastic cash balance problem allows the decision-maker to reduce the inventory level. Let $a_t = y_t - x_t$ denote the ordering ($a_t > 0$) or return ($a_t < 0$) quantity. The transaction cost is a piecewise linear function: 
\begin{equation*}
    c(y, x) = 
    \left\{
    \begin{aligned}
        & k(y - x), \quad &y \ge x,\\
        & q(x - y), \quad &y < x,
    \end{aligned}
    \right.
\end{equation*}
with $k + q \ge 0$. The assumption for $k + q \ge 0$ implies that the unit refund can not exceed the unit ordering cost. Therefore, the transaction cost function is jointly convex in $(x, y)$. 

Additionally, the random demand $D_t$ can be positive or negative. Negative demand means the decision-maker receives more returns from customers than their purchase. For simplicity, we assume that demands among different periods are independent and identically distributed. One can extend our results to the setting when random demands are not identically distributed. Let $F_D$ represent the cumulative distribution function of random demands. All other settings are the same as the inventory model in Section~\ref{section: inventory system}.

For convenience, we recall some useful notations. Let $L_t$ denote the expected cost for the period $t$ as a function of the inventory level $y_t$ after the ordering and return decisions:
\begin{equation*}
    \begin{aligned}
        L_t(y_t) = \mathbb{E}_{D_t}\left[h_t(y_t - D_t)^+ + b_t(y_t - D_t)^- \right].
    \end{aligned}
\end{equation*}
Then the per-period cost in the MDP framework is $C_t(s_t, a_t) = c(s_t + a_t, s_t) + L_t(s_t + a_t)$. The objective of the decision-maker is to minimize the total expected costs over the finite-horizon $T$:
\begin{equation}\label{optimization: cash balance objective}
    \min_{\{a_t\}_{t=1}^T} \E \left[ \sum_{t=1}^{T} \left( c(s_t + a_t, s_t) + L_t(s_t + a_t) \right) \Biggl| s_1 \sim \rho \right],
\end{equation}
with $s_{t+1} = s_t + a_t - D_t$. Like the classic stochastic inventory system, one can present the optimization problem (\ref{optimization: cash balance objective}) by a dynamic program. Let $V_t^*$ be the cost-to-go function which starts at the beginning of period $t$, it satisfies the Bellman equation (\ref{bellman optimality equation}):
\begin{equation*}
    \begin{aligned}
        V_t^*(s_t) = \min_{y_t} \left\{ c(y_t, s_t) + L_t(y_t) + \E_{D_t} \left[ V_{t+1}^*(y_t - D_t) \right] \right\},
    \end{aligned}
\end{equation*}
with $V_{T+1}^*(\cdot) = 0$. By mathematical induction, we can show the convexity of cost-to-go functions since the transaction cost function is jointly convex in $(x,y)$ and $L_t(y_t)$ is a convex function. Define $f_t(x) = L_t(x) + \E_{D_t} \left[ V_{t+1}^*(x - D_t) \right]$, we can rewrite the dynamic programming recursion:
\begin{equation*}
    V_t^*(s_t) = \min \left\{ \min_{y_t \ge s_t} \left\{ k(y_t - s_t) + f_t(y_t) \right\}, \min_{y_t \le s_t} \left\{ q(s_t - y_t) + f_t(y_t) \right\} \right\}.
\end{equation*}

\citet{whisler1967stochastic} and \citet{eppen1969cash} studied the stochastic cash balance problem and proved the optimality of the two-sided base-stock policy. That is, at period $t$ there exists two parameters $\underline{\theta}_t$ and $\bar{\theta}_t$ with $\underline{\theta}_t \le \bar{\theta}_t$, such that the optimal inventory level $y_t(s_t)$ satisfies:
\begin{equation*}
    y_t(s_t) = 
    \left\{
    \begin{aligned}
        & \underline{\theta}_t, \quad & s_t \le \underline{\theta}_t,\\
        & s_t, \quad & \underline{\theta}_t < s_t < \bar{\theta}_t,\\
        & \bar{\theta}_t, \quad & s_t \ge \bar{\theta}_t.
    \end{aligned}
    \right.
\end{equation*}

Based on the convexity of the cost-to-go functions, if we let
\begin{equation*}
    \left\{
    \begin{aligned}
        \underline{f}_t(y_t) &= ky_t + f_t(y_t) = k y_t + L_t(y_t) + \E_{D_t} \left[ V_{t+1}^* (y_t - D_t) \right], \\
        \bar{f}_t(y_t) &= -qy_t + f_t(y_t) = -q y_t + L_t(y_t) + \E_{D_t} \left[ V_{t+1}^* (y_t - D_t) \right],
    \end{aligned}
    \right.
\end{equation*}
then both $\underline{f}_t(y_t)$ and $\bar{f}_t(y_t)$ are convex functions with minimizers $\underline{\theta}_t^*$ and $\bar{\theta}_t^*$. It is well known that the optimal parameters for the two-sided base-stock policy are $\underline{\theta}_t^*$ and $\bar{\theta}_t^*$ \citep{whisler1967stochastic, eppen1969cash}.

Let $\theta_t = (\underline{\theta}_t, \bar{\theta}_t)$ denote the parameters at period $t$. Then we can rewrite the two-sided base-stock policy by  $a_t = \pi_t(s_t | \theta_t) = (s_t \vee \underline{\theta}_t) \wedge \bar{\theta}_t - s_t$. Thus the policy gradient objective function is
\begin{equation*}
    l(\theta) = \E \left[ \sum_{t=1}^{T} \left( c \left( (s_t \vee \underline{\theta}_t) \wedge \bar{\theta}_t , s_t \right) + L_t \left( (s_t \vee \underline{\theta}_t) \wedge \bar{\theta}_t \right) \right) \Biggl| s_1 \sim \rho \right],
\end{equation*}
with $s_{t+1} = (s_t \vee \underline{\theta}_t) \wedge \bar{\theta}_t - D_t$. In the remaining part of this section, we only analyze the nonconvex landscape of $l(\theta)$ over a convex set $\Theta = \Theta_1\times,\dots,\times\Theta_T$ with $\Theta_t = \{ (\underline{\theta}_t, \bar{\theta_t}): \underline{B} \le \underline{\theta}_t \le \bar{\theta_t} \le \bar{B} \}$. This is reasonable in practice since one can treat $\bar{B}$ as the capacity of the warehouse, and the lower bound $\underline{B}$ ensures the firms will never hold too many backlogged demands.

\subsection{Nonconvex Landscape}
Like Section~\ref{section: inventory system}, we establish the P{\L}K condition for the policy gradient objective function by verifying the conditions in Theorem~\ref{theorem: main result}. Before showing the main results, we make some assumptions.

\begin{assumption}
\label{assumption: cash balance}
    Assume that all the following assumptions hold.
    \begin{henumerate}
        \item \label{assumption: cash balance independence of rv} The initial state $s_1$ and random demands $D_t$ in different periods $t$ are independent of each other. 
        \item \label{assumption: cash balance Lipschitz of rv} The initial distribution $\rho$ is $L_\rho$-Lipschitz continuous. Furthermore, the cumulative distribution function $F_D$ of random demands is $L_D$-Lipschitz continuous. 
        \item \label{assumption: cash balance KL} There exists $\alpha>0$ such that, for all $t\in[T]$, we have $\rho(\underline{B})\ge \alpha$, $\rho(\bar{B}) \le 1-\alpha$, $\mathbb{P}(D_t\ge \bar{B} - \underline{B})\ge \alpha$, and $\mathbb{P}(D_t \le \underline{B} - \bar{B})\ge \alpha$.
        \item \label{assumption: cash balance SC} The probability density function of the random demand $D_t$ is bounded below by $\mu_D > 0$ on $[\underline{B}, \bar{B}]$. 
    \end{henumerate}
\end{assumption}

Assumption~\ref{assumption: cash balance} plays the same role as Assumption~\ref{assumption: inventory}. In particular, Assumption~\ref{assumption: cash balance}.\ref{assumption: cash balance KL} provides sufficient conditions ensuring that the state $s_t$ crosses the base-stock levels with positive probability, i.e., $\mathbb{P}(s_t \ge \bar{\theta}_t)>0$ and $\mathbb{P}(s_t \le \underline{\theta}_t)>0$. These conditions ensure that the parameters $(\underline{\theta}_t,\bar{\theta}_t)$ affect the objective $l(\theta)$ in a nontrivial way, thereby ruling out suboptimal stationary points caused by inactive parameters.

For any policy $\pi_\theta \in \Pi_\Theta$, the expected Q-value functions satisfy the Bellman equation (\ref{bellman equation}):
\begin{equation*}
\label{bellman equation: cash balance}
    \begin{aligned}
        & \ \E_{s_t \sim \rho_t(\cdot|\pi_\theta)} \left[ Q_t^{\pi_\theta} \left( s_t, \pi_t(s_t | \theta_t) \right) \right]\\
        = & \ \E_{s_t \sim \rho_t(\cdot|\pi_\theta)} \left[ c \left( (s_t \vee \underline{\theta}_t) \wedge \bar{\theta}_t, s_t \right) + L_t \left( (s_t \vee \underline{\theta}_t) \wedge \bar{\theta}_t \right) + \E_{D_t} \left[ V_{t+1}^{\pi_\theta} \left( (s_t \vee \underline{\theta}_t) \wedge \bar{\theta}_t - D_t \right) \right] \right].
    \end{aligned}
\end{equation*}

Same as Section~\ref{section: inventory system}, the continuous differentiability of the expected Q-value function comes from the continuity of the cumulative distribution function of $s_t$, which holds under Assumption~\ref{assumption: cash balance}.2. Recall the definition $f_t(x) = L_t(x) + \E_{D_t} \left[ V_{t+1}^{\pi_{\theta^*}}(x - D_t) \right]$. We get the expression of expected optimal Q-value functions by replacing $\pi_\theta$ with $\pi^*$:
\begin{equation*}
    F_t(\theta_t) \coloneqq \E_{s_t \sim \rho_t(\cdot|\pi_\theta)} \left[ Q_t^{\pi_{\theta^*}} \left( s_t, \pi_t(s_t|\theta_t) \right) \right] = \E_{s_t \sim \rho_t(\cdot|\pi_\theta)} \left[ c \left( (s_t \vee \underline{\theta}_t) \wedge \bar{\theta}_t, s_t \right) +  f_t \left( (s_t \vee \underline{\theta}_t) \wedge \bar{\theta}_t \right) \right].
\end{equation*}

The structure of the expected optimal Q-value function is close to that for the inventory system in Section~\ref{section: inventory system}. The following lemma establishes the P{\L}K condition for the expected optimal Q-value function by a similar proof.

\begin{lemma}
\label{lemma: cash balance single stage KL}
    Suppose that Assumption~\ref{assumption: cash balance} holds. The expected Q-value function is continuously differentiable on $\Theta_t$ for any $t\in[T]$. In addition, the expected optimal Q-value function $F_t(\theta_t)$ satisfies the P{\L}K condition on $\Theta_t$ with P{\L}K constant $\min_{t\in[T]}\{h_t+b_t\}\mu_D\alpha^2$ for any $t\in[T]$.
\end{lemma}

Again, we need an explicit expression of the gradient $\nabla l(\theta)$ to validate the remaining conditions. Let $V_t^{\pi_\theta}$ be the value function that starts at the beginning of period $t$. It satisfies the Bellman equation (\ref{bellman equation}):
\begin{equation*}
    \begin{aligned}
        V_t^{\pi_\theta}(s_t) &= c\left( \pi_t(s_t|\theta_t), s_t \right) + L_t \left(\pi_t(s_t|\theta_t)\right) + \E_{D_t} \left[ V_{t+1}^{\pi_\theta}\left(\pi_t(s_t|\theta_t) - D_t \right) \right]\\
        &= c \left( (s_t \vee \underline{\theta}_t) \wedge \bar{\theta}_t, s_t \right) + L_t \left( (s_t \vee \underline{\theta}_t) \wedge \bar{\theta}_t \right) + \E_{D_t} \left[ V_{t+1}^{\pi_\theta} \left( (s_t \vee \underline{\theta}_t) \wedge \bar{\theta}_t - D_t \right) \right],
    \end{aligned}
\end{equation*}
with $V_{T+1}^{\pi_\theta}(\cdot) = 0$. The subsequent proposition presents a gradient formulation for $l(\theta)$.

\begin{proposition}[Policy Gradient Expression]
\label{proposition: cash balance gradient formulation}
    For any $\theta\in \Theta$ and $t\in[T]$, the derivatives of value functions satisfy the following recursive form:
    \begin{equation}
    \label{cash balance: recursive gradient}
        (V_t^{\pi_\theta})' (s_t) = - k \mathbf{1}(s_t \le \underline{\theta}_t) + q \mathbf{1}(s_t \ge \bar{\theta}_t) + \mathbf{1}(\underline{\theta}_t < s_t < \bar{\theta}_t) \times \left( L'_t(s_t) + \E_{D_t}\left[ (V^{\pi_\theta}_{t+1})'(s_t - D_t) \right] \right)
    \end{equation}
    with $(V_{T+1}^{\pi_\theta})(\cdot) = 0$. Additionally, the policy gradient objective function $l(\theta)$ has the following gradient form for any $t\in[T]$:
    \begin{equation*}
        \left\{
        \begin{aligned}
            \frac{\partial}{\partial \underline{\theta}_t} l(\theta) &= \E_{s_t \sim \rho_t(\cdot|\pi_\theta)} \left[ \mathbf{1}(\underline{\theta}_t \ge s_t) \times \left( k + L'_t(\underline{\theta}_t) + \E_{D_t} \left[ (V_{t+1}^{\pi_\theta})' (\underline{\theta}_t - D_t) \right] \right) \right], \\
            \frac{\partial}{\partial \bar{\theta}_t} l(\theta) &= \E_{s_t \sim \rho_t(\cdot|\pi_\theta)} \left[ \mathbf{1}( \bar{\theta}_t \le s_t) \times \left( -q + L'_t(\bar{\theta}_t) + \E_{D_t} \left[ (V_{t+1}^{\pi_\theta})' (\bar{\theta}_t - D_t) \right] \right) \right].
        \end{aligned}
        \right.
    \end{equation*}
\end{proposition}

With the explicit expression of the gradient, we can verify the bounded gradient condition and the sequential decomposition inequality by the following two lemmas.

\begin{lemma}
\label{lemma: cash balance bounded gradient}
    Suppose that Assumption \ref{assumption: cash balance} holds. It follows that the policy gradient objective function $l(\theta)$ has bounded gradients for any $\theta\in\Theta$ and $t\in[T]$:
    \begin{equation*}
        \left\| \nabla_{\theta_t} l(\theta) \right\|_2 \le 2 \bigl( k + |q| + \max_{t\in[T]} \left\{ \max\{h_t, b_t\} \right\} \bigr) T.
    \end{equation*}
\end{lemma}

\begin{lemma}
\label{lemma: cash balance bounded sequence gradient difference}
    Suppose that Assumption~\ref{assumption: cash balance} holds. For any $\theta\in\Theta$ and $1\le t < k \le T$, sequential decomposition inequalities hold, i.e.,
    \begin{equation*}
        \begin{aligned}
            & \left\| \nabla_{\theta_t} l(\theta_{[1:k-1]}, \theta_{k}, \theta_{[k+1:T]}^*) - \nabla_{\theta_t} l(\theta_{[1:k-1]}, \theta_{k}^*, \theta_{[k+1:T]}^*) \right\|_2\\
            \le \ & \frac{2L_D}{\alpha} \left( \E_{s_k \sim \rho_k(\cdot|\pi_\theta)}\left[ Q_k^* \left( s_k, \pi_k(s_k|\theta_k) \right) \right] - \E_{s_k \sim \rho_k(\cdot|\pi_\theta)}\left[ Q_k^* \left( s_k, \pi_k(s_k|\theta_k^*) \right) \right] \right).
        \end{aligned}
    \end{equation*}
\end{lemma}

Equipped with P{\L}K condition of optimal Q-value functions (Lemma~\ref{lemma: cash balance single stage KL}), bounded gradient condition (Lemma~\ref{lemma: cash balance bounded gradient}), and sequential decomposition inequality (Lemma~\ref{lemma: cash balance bounded sequence gradient difference}), we can demonstrate the P{\L}K condition of the policy gradient objective function $l(\theta)$ by applying Theorem~\ref{theorem: main result}.

\begin{theorem}
    \label{theorem: cash balance multistage KL}
    Suppose that Assumption~\ref{assumption: cash balance} holds. The policy gradient objective function $l(\theta)$ of the stochastic cash balance problem satisfies the P{\L}K condition on $\Theta$ with P{\L}K constant
    \begin{equation*}
        \mu_l = \frac{\min_{t\in[T]}\{h_t+b_t\}^3\mu_D^3 \alpha^8}{16eL_D^2 \left( k + |q| + \max_{t\in[T]} \left\{ \max\{h_t, b_t\} \right\} \right)^2 T^4 }.
    \end{equation*}
\end{theorem}

\proof{Proof of Theorem~\ref{theorem: cash balance multistage KL}} Plugging Lemmas~\ref{lemma: cash balance single stage KL}, \ref{lemma: cash balance bounded gradient}, and \ref{lemma: cash balance bounded sequence gradient difference} into Theorem~\ref{theorem: main result} yields the result. \Halmos

Similar to the inventory system in Section~\ref{section: inventory system}, we establish an $\tilde{O} \left( \epsilon^{-1} \text{poly}(T) \right)$ sample complexity of stochastic policy gradient methods converging to globally optimal policies by Lemma~\ref{lemma: PGD KL}. To the best of our knowledge, this is the first sample complexity result for data-driven methods solving the stochastic cash balance problem. Additionally, one can check that P{\L}K condition holds for the stochastic cash balance problem with Markov-modulated demand as well.

\section{Numerical Experiments}\label{sec: numerical exp}
In this section, we present three numerical experiments to demonstrate that policy gradient (PG) methods achieve strong solution quality and remain computationally efficient. The experiment settings include (i) standard inventory models, (ii) inventory control with Markov-modulated demand as described in Section~\ref{section: inventory system}, and (iii) stochastic cash balance problems as described in Section~\ref{section: cash balance}. In setting (i), we compare PG with several established methods from the literature. In settings (ii) and (iii), the best-known approaches, to our knowledge, rely on sample average approximation and dynamic programming; we are not aware of prior work applying policy gradient methods to these problems. All the experiments are conducted on a MacBook Pro with an Apple M3 chip and 18 GB RAM.

\subsection{Inventory Models}
\label{sec:exp inv}
We first evaluate PG in inventory models without an exogenous Markov chain, which reduces to the standard inventory control model under a base-stock policy class. Several algorithms have been proposed for this benchmark problem. We compare PG against four representative approaches from the literature: \citet[hereafter KT2008]{kunnumkal2008using}, \citet[hereafter HR2014]{huh2014online}, \citet[hereafter CS2019]{cheung2019sampling}, and SAIL from \citet{qin2023sailing}.

Following the setup of \citet{kunnumkal2008using}, we consider two demand distributions, denoted UN and EX. Under UN, demand in period $t$ is uniformly distributed: $D_t \sim \mathrm{Unif}[l_t, u_t]$. Under EX, demand is exponentially distributed with rate $\lambda_t$, i.e., $D_t \sim \mathrm{Exp}(\lambda_t)$. For each $t \in [T]$, we draw $l_t$, $u_t$ independently from $\mathrm{Unif}[1,20]$ and $\lambda_t \sim \mathrm{Unif}[1,10]$. We set the per-unit backlogging cost to $b_t = 1$ and impose an upper bound $B=20$ on the base-stock levels. The initial inventory $x_1$ is drawn from $\mathrm{Unif}[-10,20]$.

Let $\theta^{\text{alg}}$ denote the solution returned by algorithm $\text{alg} \in \{\text{KT2008}, \text{HR2014}, \text{CS2019}, \text{SAIL}, \text{PG}\}$. For example, $\theta^{\text{PG}}$ is the final iterate produced by the policy gradient method. We assess each algorithm along two dimensions: (i) the suboptimality gap $l(\theta^{\text{alg}})-l(\theta^*)$ and (ii) runtime (in seconds). Here, $\theta^*$ and $l(\theta^*)$ are computed using a high-accuracy benchmark based on sample average approximation (SAA) with a very large number of samples combined with dynamic programming (DP) with a very large number of discretization grids. To reduce Monte Carlo variability arising from random demand realizations, we perform $25$ independent runs (with different random seeds) and report the average suboptimality gap and runtime. Among the benchmark methods, KT2008, HR2014, and PG are first-order algorithms. To ensure a fair comparison, we initialize the base-stock vector at $\theta^0_t = 0$ for all $t \in [T]$ and use the stepsize schedule $\gamma^k = 100/(40+k)$ at iteration $k$ for all first-order algorithms. We run each method for $10{,}000$ iterations and use $N=8$ independent demand trajectories per iteration to construct the gradient estimator. Thus, each first-order method uses $N \times 10{,}000 = 80{,}000$ trajectories in total, corresponding to $80{,}000 \times T$ demand samples overall. We evaluate performance at the final iterate. For CS2019 and SAIL, we draw $80{,}000$ demand samples per period, matching the per-period sample size used by the first-order methods.

\begin{table}[htbp]
    \small
    \centering
    \caption{Suboptimality gaps and runtimes of algorithms in different problem settings for inventory models.}
    \renewcommand{\arraystretch}{1.2}
    \begin{tabular}{c *{5}{cc}}
    \hline\hline
    \multirow{2}{*}{Problem Setting}
      & \multicolumn{2}{c}{KT2008}
      & \multicolumn{2}{c}{HR2014}
      & \multicolumn{2}{c}{CS2019}
      & \multicolumn{2}{c}{SAIL}
      & \multicolumn{2}{c}{PG} \\
    \cline{2-11}
      & gap & runtime
      & gap & runtime
      & gap & runtime
      & gap & runtime
      & gap & runtime \\
    \hline
    $(20, 0.1, \text{UN})$
      & 0.0018 & 4.0598
      & 0.0007 & 1.9491
      & 0.0001 & 3.3513
      & 0.0066 & 48.7723
      & 0.0007 & \textbf{0.8614} \\
    $(20, 0.25, \text{UN})$
      & 0.0042 & 4.0664
      & 0.0017 & 1.9430
      & 0.0002 & 17.3859
      & 0.0056 & 49.5625
      & 0.0015 & \textbf{0.8389} \\
    $(50, 0.1, \text{UN})$
      & 0.0088 & 21.1454
      & 0.0016 & 4.3729
      & 0.0002 & 8.1829
      & 0.0085 & 136.8843
      & 0.0016 & \textbf{1.9816} \\
    $(50, 0.25, \text{UN})$
      & 0.0206 & 21.1597
      & 0.0042 & 4.3746
      & 0.0004 & 43.7511
      & 0.0095 & 138.0567
      & 0.0040 & \textbf{2.0194} \\
    $(100, 0.1, \text{UN})$
      & 0.0196 & 77.4049
      & 0.0032 & 8.4272
      & 0.0004 & 15.2019
      & 0.0291 & 298.0215
      & 0.0030 & \textbf{3.6420} \\
    $(100, 0.25, \text{UN})$
      & 0.0521 & 78.5368
      & 0.0076 & 8.4660
      & 0.0008 & 79.5898
      & 0.0273 & 298.8256
      & 0.0080 & \textbf{3.8819} \\
    \hline
    $(20, 0.1, \text{EX})$
      & 0.0016 & 3.8650
      & 0.0007 & 2.1156
      & 0.0033 & 9.8270
      & 0.0017 & 59.0148
      & 0.0010 & \textbf{0.6647} \\
    $(20, 0.25, \text{EX})$
      & 0.0047 & 3.9257
      & 0.0016 & 2.0292
      & 0.0008 & 124.0198
      & 0.0018 & 64.3200
      & 0.0017 & \textbf{0.6766} \\
    $(50, 0.1, \text{EX})$
      & 0.2334 & 20.7678
      & 0.0021 & 4.8242
      & 0.0076 & 17.5380
      & 0.0031 & 161.2359
      & 0.0022 & \textbf{1.5313} \\
    $(50, 0.25, \text{EX})$
      & 0.0135 & 20.6941
      & 0.0044 & 4.6663
      & 0.0024 & 248.3746
      & 0.0040 & 162.1205
      & 0.0046 & \textbf{1.5556} \\
    
    $(100, 0.1, \text{EX})$
      & 0.3259 & 76.6661
      & 0.0042 & 9.2802
      & 0.0160 & 41.3954
      & 0.0063 & 342.6958
      & 0.0045 & \textbf{2.9345} \\
    $(100, 0.25, \text{EX})$
      & 0.0464 & 77.2410
      & 0.0093 & 9.2519
      & 0.0050 & 556.1555
      & 0.0081 & 341.2839
      & 0.0094 & \textbf{3.0652} \\
    \hline\hline
    \end{tabular}
    \label{table: inventory}
\end{table}

The computational results are summarized in Table~\ref{table: inventory}. The first column reports the problem setting, indexed by $(T,h,D) \in \{20,50,100\}\times\{0.1,0.25\}\times\{\text{UN},\text{EX}\}$, where $T$ is the planning horizon, $h$ is the per-unit holding cost, and $D$ specifies the demand distribution. The remaining ten columns are organized in pairs: for each algorithm, we report its average suboptimality gap and average runtime (in seconds). For example, Columns~2--3 correspond to KT2008 and report, respectively, the mean suboptimality gap and mean runtime, averaged over 25 independent runs.

Our computational results indicate that PG delivers strong performance in both solution quality and computational efficiency. First, PG consistently achieves smaller suboptimality gaps and shorter runtimes than KT2008, with the largest differences occurring at $T=100$. In these instances, KT2008 can produce gaps exceeding $0.3$, whereas PG typically attains gaps below $0.01$. This pattern is consistent with the fact that KT2008 employs a potentially biased gradient estimator, while PG uses an unbiased estimator. Second, CS2019 and SAIL can also reach suboptimality gaps below $0.01$, but their runtimes increase rapidly with the horizon length $T$. For example, when $T=100$, SAIL requires more than $300$ seconds to terminate, whereas PG completes in under $5$ seconds. Finally, HR2014 performs comparably to PG in our experiments, which is reasonable given that HR2014 can be viewed as a minor modification of PG. From a theoretical standpoint, however, \citet{huh2014online} establishes a convergence guarantee for HR2014 with exponential dependence on $T$, whereas our analysis yields a polynomial dependence on $T$ for PG. Although the resulting polynomial order may be large, the computational evidence suggests that PG scales favorably with $T$ in practice.

\subsection{Markov Modulated Demands}
\begin{table}[b]
    \small
    \centering
    \caption{Suboptimality gaps and runtimes of PG for inventory models with Markov-modulated demand.}
    \renewcommand{\arraystretch}{1.2}
    \begin{tabular}{cc *{3}{ccc}}
    \hline\hline
    \multirow{2}{*}{Problem Setting}
      & \multirow{2}{*}{OP}
      & \multirow{2}{*}{DP Runtime}
      & \multicolumn{3}{c}{Gap}
      & \multicolumn{3}{c}{PG Runtime}\\
    \cline{4-9}
      & & 
      & MIN & AVG & MAX
      & MIN & AVG & MAX \\
    \hline
    $(20, 4, 0.1, \text{UN})$ & 6.2158 & 1061.586
    & 0.0021 & 0.0031 & 0.0040
    & 3.5858 & 3.6537 & 3.7625\\
    $(20, 4, 0.25, \text{UN})$ & 13.5999 & 1048.975
    & 0.0057 & 0.0081 & 0.0110
    & 3.6214 & 3.6579 & 3.7161\\
    $(20, 7, 0.1, \text{UN})$ & 5.8975 & 1918.619
    & 0.0049 & 0.0062 & 0.0077
    & 4.7029 & 4.7574 & 4.8113\\
    $(20, 7, 0.25, \text{UN})$ & 12.9438 & 1900.071
    & 0.0082 & 0.0098 & 0.0125
    & 4.7381 & 4.8286 & 4.8850\\
    $(50, 4, 0.1, \text{UN})$ & 15.3875 & 2567.935
    & 0.0064 & 0.0079 & 0.0098
    & 8.9757 & 9.0557 & 9.3422\\
    $(50, 4, 0.25, \text{UN})$ & 33.6448 & 2578.570
    & 0.0153 & 0.0198 & 0.0235
    & 8.6873 & 8.8549 & 8.9486\\
    $(50, 7, 0.1, \text{UN})$ & 15.5202 & 4679.607
    & 0.0131 & 0.0146 & 0.0165
    & 11.5365 & 11.6705 & 11.8529\\
    $(50, 7, 0.25, \text{UN})$ & 33.9814 & 4680.607
    & 0.0187 & 0.0210 & 0.0259
    & 11.4649 & 11.6540 & 11.8162\\
    $(100, 4, 0.1, \text{UN})$ & 29.0929 & 6317.594
    & 0.0133 & 0.0158 & 0.0180
    & 17.4354 & 17.5590 & 17.6526\\
    $(100, 4, 0.25, \text{UN})$ & 63.4891 & 6309.872
    & 0.0333 & 0.0393 & 0.0441
    & 17.3796 & 17.5005 & 17.6396\\
    $(100, 7, 0.1, \text{UN})$ & 29.8144 & 11193.468
    & 0.0294 & 0.0322 & 0.0370
    & 22.6566 & 22.8394 & 23.3496\\
    $(100, 7, 0.25, \text{UN})$ & 65.1083 & 11361.589
    & 0.0419 & 0.0452 & 0.0493
    & 22.7001 & 22.8076 & 23.0601\\

    \hline
    $(20, 4, 0.1, \text{EX})$ & 26.5243 & 1042.334
    & 0.0046 & 0.0068 & 0.0097
    & 3.1229 & 3.2278 & 3.6493\\
    $(20, 4, 0.25, \text{EX})$ & 43.7490 & 1049.311
    & 0.0058 & 0.0080 & 0.0114
    & 3.1432 & 3.2203 & 3.5707\\
    $(20, 7, 0.1, \text{EX})$ & 28.1818 & 1905.240
    & 0.0106 & 0.0140 & 0.0193
    & 4.0285 & 4.0815 & 4.3046\\
    $(20, 7, 0.25, \text{EX})$ & 46.6560 & 1873.820
    & 0.0127 & 0.0172 & 0.0251
    & 4.0349 & 4.0737 & 4.1058\\
    $(50, 4, 0.1, \text{EX})$ & 67.0168 & 2593.539
    & 0.0130 & 0.0172 & 0.0200
    & 7.6383 & 7.7100 & 7.7887\\
    $(50, 4, 0.25, \text{EX})$ & 109.7501 & 2573.056
    & 0.0155 & 0.0214 & 0.0258
    & 7.5757 & 7.6479 & 7.7044\\
    $(50, 7, 0.1, \text{EX})$ & 69.3063 & 4652.499
    & 0.0303 & 0.0353 & 0.0416
    & 9.8209 & 9.9193 & 9.9782\\
    $(50, 7, 0.25, \text{EX})$ & 113.7005 & 4674.328
    & 0.0315 & 0.0414 & 0.0467
    & 9.9546 & 10.0723 & 10.4359\\
    $(100, 4, 0.1, \text{EX})$ & 138.9419 & 6233.653
    & 0.0289 & 0.0347 & 0.0407
    & 14.8920 & 15.0216 & 15.3145\\
    $(100, 4, 0.25, \text{EX})$ & 226.8542 & 6208.550
    & 0.0372 & 0.0422 & 0.0500
    & 15.4383 & 15.5220 & 15.6675\\
    $(100, 7, 0.1, \text{EX})$ & 140.7806 & 10919.911
    & 0.0588 & 0.0690 & 0.0795
    & 19.4040 & 19.5352 & 19.7449\\
    $(100, 7, 0.25, \text{EX})$ & 229.9596 & 11130.466
    & 0.0731 & 0.0848 & 0.1000
    & 19.4703 & 20.1159 & 20.9017\\
    \hline\hline
    \end{tabular}
    \label{table:exp markov}
\end{table}
We next evaluate the performance of PG for inventory models with Markov-modulated demand described in Section~\ref{section: inventory system}. For the exogenous Markov chain, we generate the transition probabilities by drawing $p(\cdot \mid i)$, independently for $i \in \mathcal{I}$, from a Dirichlet distribution with concentration parameter $(2,\ldots,2)$. As in Section~\ref{sec:exp inv}, we consider two demand families, UN and EX. For each period $t \in [T]$ and state $i \in \mathcal{I}$, under UN we have $D_{t,i} \sim \mathrm{Unif}[l_{t,i}, u_{t,i}]$, and under EX we have $D_{t,i} \sim \mathrm{Exp}(\lambda_{t,i})$. We draw $l_{t,i}$, $u_{t,i}$ independently from $\mathrm{Unif}[1,20]$ and $\lambda_{t,i} \sim \mathrm{Unif}[1,10]$ for all $t \in [T]$ and $i \in \mathcal{I}$. All other parameters follow Section~\ref{sec:exp inv}. In this subsection, we test the performance of PG under a constant stepsize. The stepsize is selected from the set $\{1, 0.5, 0.2, 0.1, 0.05, 0.02, 0.01\}$.

Table~\ref{table:exp markov} summarizes the computational results for PG. The first column denotes the problem setting with parameters $(T, |\ci|, h, D) \in \{20, 50, 100\} \times \{4, 7\} \times \{0.1, 0.25\} \times \{\text{UN}, \text{EX}\}$. We select the number of states for the exogenous Markov chain, $|\ci|\in\{4, 7\}$, to mimic seasonal ($|\ci| = 4$) and weekly ($|\ci| = 7$) patterns. The second column shows the optimal function values achieved by SAA and DP with extremely large numbers of samples and grids. The remaining six columns are divided into two parts. The first three columns report the minimum, average, and maximum suboptimality gaps achieved across $25$ runs, respectively. The last three columns record the runtime of policy gradient methods. Similar to Table~\ref{table: inventory}, policy gradient methods solve inventory models with Markov-modulated demand fast and accurately. All worst-case suboptimality gaps are bounded by $0.1$, and the longest time to terminate is less than $21$ seconds. Although in theory we have established the sample and computational complexities with a high-order polynomial dependence on $T$, our computational results demonstrate the scalability of policy gradient methods with respect to $T$. We further conduct a series of robustness checks to assess the stability of the policy gradient methods’ performance. Details are provided in Appendix~\ref{Appendix: Robustness Check}.

\subsection{Stochastic Cash Balance Problems}
At last, we test the performance of PG for solving stochastic cash balance problems discussed in Section~\ref{section: cash balance}. Unlike the standard inventory model, random demands in the stochastic cash balance problem can be negative. Therefore, we use two different distributions to sample the random demands. For each time period $t \in [T]$, the random demands either follow a uniform distribution (UN) over $[l_t, u_t]$ or a normal distribution (NR) with mean $\mu_t$ and standard deviation $\sigma_t$. For each parameter, we sample $l_t \sim \text{Unif}[-20, -10]$, $u_t \sim \text{Unif}[10, 20]$, $\mu_t \sim \text{Unif}[-20, 20]$, and $\sigma_t \sim \text{Unif}[10, 20]$. We set the per-unit ordering cost to $k = 0.5$ and the return cost to $q = -0.25$. All the other settings are the same as those in inventory models in Section~\ref{sec:exp inv}.

\begin{table}[htbp]
    \small
    \centering
    \caption{Suboptimality gaps and runtimes of PG for stochastic cash balance problems.}
    \renewcommand{\arraystretch}{1.2}
    \begin{tabular}{cc *{3}{ccc}}
    \hline\hline
    \multirow{2}{*}{Problem Setting}
      & \multirow{2}{*}{OP}
      & \multirow{2}{*}{DP Runtime}
      & \multicolumn{3}{c}{Gap}
      & \multicolumn{3}{c}{PG Runtime}\\
    \cline{4-9}
      & &
      & MIN & AVG & MAX
      & MIN & AVG & MAX \\
    \hline
    $(20, 0.1, \text{UN})$ & 49.4905 & 376.099
    & 0.0076 & 0.0123 & 0.0156
    & 2.2449 & 2.2744 & 2.3318\\
    $(20, 0.25, \text{UN})$ & 82.1427 & 375.528
    & 0.0031 & 0.0053 & 0.0097
    & 2.2741 & 2.3081 & 2.3944\\
    $(50, 0.1, \text{UN})$ & 114.6947 & 924.637
    & 0.0238 & 0.0275 & 0.0312
    & 5.3240 & 5.3912 & 5.5816\\
    $(50, 0.25, \text{UN})$ & 195.8018 & 926.638
    & 0.0101 & 0.0135 & 0.0183
    & 5.5156 & 5.6196 & 5.8634\\
    $(100, 0.1, \text{UN})$ & 222.9917 & 2082.380
    & 0.0511 & 0.0567 & 0.0668
    & 10.7756 & 10.9304 & 11.5080\\
    $(100, 0.25, \text{UN})$ & 386.7041 & 2074.342
    & 0.0212 & 0.0257 & 0.0309
    & 10.9924 & 11.2792 & 11.4449\\
    
    \hline
    $(20, 0.1, \text{NR})$ & 91.1015 & 315.547
    & 0.0119 & 0.0183 & 0.0309
    & 2.2061 & 2.2208 & 2.2485\\
    $(20, 0.25, \text{NR})$ & 135.5680 & 315.006
    & 0.0028 & 0.0063 & 0.0134
    & 2.2057 & 2.2486 & 2.3368\\
    $(50, 0.1, \text{NR})$ & 262.4016 & 802.203
    & 0.0228 & 0.0467 & 0.1057
    & 5.3854 & 5.4315 & 5.7128\\
    $(50, 0.25, \text{NR})$ & 372.0718 & 809.170
    & 0.0083 & 0.0139 & 0.0269
    & 5.3910 & 5.4331 & 5.4788\\
    $(100, 0.1, \text{NR})$ & 512.7532 & 1822.244
    & 0.0658 & 0.0814 & 0.0978
    & 10.7331 & 10.8022 & 10.8569\\
    $(100, 0.25, \text{NR})$ & 731.1269 & 1864.371
    & 0.0133 & 0.0232 & 0.0337
    & 10.8923 & 11.1434 & 11.7270\\
    \hline\hline
    \end{tabular}
    \label{table:exp cash}
\end{table}

Table~\ref{table:exp cash} summarizes our computational results. The first column reports the problem setting with parameters $(T, h, D) \in \{20, 50, 100\} \times \{0.1, 0.25\} \times \{\text{UN}, \text{NR}\}$. The remaining columns report the same meaning as in Table~\ref{table:exp markov}. Similarly, the computational results demonstrate the strong performance of policy gradient methods. Even in the worst case, policy gradient methods can obtain a solution within $12$ seconds, with a suboptimality gap of approximately $0.1$.

\section{Conclusion}
This work provides a framework with several structural conditions to establish the P{\L}K condition for policy gradient optimization problems of finite-horizon MDPs with general state and action spaces. Despite nonconvexity, the P{\L}K condition guarantees a linear convergence rate for exact policy gradient methods and an $\tilde{\co}(\epsilon^{-1})$ sample complexity for stochastic policy gradient methods. Our framework covers a broad range of control and operations models, including entropy-regularized tabular MDPs, LQR problems, multi-period inventory systems with Markov-modulated demand, and stochastic cash balance problems. Furthermore, we establish the first sample complexity solving the stochastic cash balance problem and multi-period inventory system with Markov-modulated demand allowing backorders, and an extension to the lost sales model. The complexity has a polynomial instead of an exponential dependence on the planning horizon.

Our work opens up several directions for future research. Firstly, in many of the applications discussed, one can further explore the structural properties to build a more precise characterization of the P{\L}K constant and reduce the dependence on $T$. Secondly, our results build upon the P{\L}K condition of expected optimal Q-value functions, requiring strongly convex per-period costs. It remains interesting to further generalize the results to general convex per-period costs for applications like inventory models. Adding regularization could be one potential solution, yet it might deteriorate the dependence on accuracy $\epsilon$. Finally, it is interesting to explore the applicability of our developed framework to other applications.

\ACKNOWLEDGMENT{Yifan Hu is supported by the Swiss National Science Foundation under NCCR Automation, grant agreement 51NF40\_180545.}

\bibliographystyle{informs2014}
\bibliography{myrefs}

\newpage

\renewcommand{\theHsection}{A\arabic{section}}
\begin{APPENDICES}
\ECHead{Online Appendices}

\section{Omitted Proofs in Section~\ref{section: landscape}}
\label{appendix: landscape}

\subsection{Strong Convexity, Gradient Dominance, and P{\L}K condition}\label{subsection: nonconvex conditions}
\begin{definition}[Gradient Dominance]
\label{definition: gradient_dominance}
    Consider a convex and compact set $\mathcal{X}\subseteq\mathbb{R}^d$ and a differentiable function $f$. Suppose that $f^*$ is the optimal objective value $f^* \coloneqq \min_{x\in\mathcal{X}} f(x)$. The function $f$ is said to be $(\alpha,\mu)$-gradient dominated over $\mathcal{X}$ if there exists constants $\alpha>0$ and $\mu\ge 0$ such that
    \begin{equation}
    \label{ineq: gradient dominance}
        f(x) - f^* \le \max_{x'\in\mathcal{X}} \left\{ \alpha\langle \nabla f(x), x - x'\rangle - \frac{\mu}{2}\|x - x'\|_2^2 \right\}, \quad\forall x\in\mathcal{X}.
    \end{equation}
\end{definition}
One important property is that the $(\alpha,0)$-gradient dominance (for $\mu > 0$) and $(\alpha,\mu)$-gradient dominance are relaxations of convexity and $\mu$-strong convexity, respectively. In particular, if $f$ is convex (or $\mu$-strongly convex), then $f$ is $(1,0)$-gradient dominated (or $(1,\mu)$-gradient dominated) by definition.

\begin{lemma}
    \label{lemma: gradient dominance implies KL}
    Consider a convex and compact set $\cx\subseteq\R^n$. If a function $f:\cx \to \R$ is $(\alpha, \mu)$-gradient dominated over $\cx$, then $f$ satisfies the P{\L}K condition with constant $\mu/\alpha^2$ over $\cx$.
\end{lemma}

Lemma~\ref{lemma: gradient dominance implies KL} is established in Appendix G of \citet{karimi2016linear}. We include its proof here for completeness and to keep the paper self-contained.

\proof{Proof of Lemma~\ref{lemma: gradient dominance implies KL}}
By definition,
\begin{equation}
\label{ineq: SC to gradient dominance}
    f(x) - f^* \le \max_{x'\in\cx} \left\{ \alpha \langle \nabla f(x), x - x' \rangle - \frac{\mu}{2} \|x - x'\|_2^2 \right\}.
\end{equation}
From the convexity of the indicator function $\delta_{\cx}(x)$, for any $s \in \partial \delta_\cx(x)$, we have
\begin{equation}
\label{ineq: convexity of indicator}
    \delta_\cx(x') - \delta_\cx(x) \ge \langle s, x' - x \rangle, \quad \forall x\in\cx, x'\in \R^n.
\end{equation}
Then for any $x\in\cx$, we can derive the following:
\begin{eqnarray*}
    f(x) - f^* & \overset{(a)}{\le} & \max_{x'\in\cx} \left\{ \alpha \langle \nabla f(x), x - x' \rangle - \frac{\mu}{2} \|x - x'\|_2^2 \right\}\\
    & \overset{(b)}{=} & \max_{x'\in\R^n} \left\{ \alpha \langle \nabla f(x), x - x' \rangle - \frac{\mu}{2} \|x - x'\|_2^2 + \alpha \left[ \delta_\cx(x) - \delta_\cx(x') \right] \right\}\\
    & \overset{(c)}{\le} & \max_{x'\in\R^n} \left\{ \alpha \langle \nabla f(x) + s, x - x' \rangle - \frac{\mu}{2} \|x - x'\|_2^2 \right\}\\
    & \overset{(d)}{=} & \frac{\alpha^2}{2\mu} \left\| \nabla f(x) + s \right\|_2^2.
\end{eqnarray*}
Here inequality (a) uses (\ref{ineq: SC to gradient dominance}), equation (b) moves constraints into the objective, inequality (c) uses (\ref{ineq: convexity of indicator}), and equation (d) solves the unconstrained quadratic optimization problem explicitly. Thus, taking the minimum value over all the subdifferentials $s\in\partial \delta_\cx(x)$, we have
\begin{equation*}
    f(x) - f^* \le \frac{\alpha^2}{2\mu} \min_{s\in\partial \delta_\cx(x)} \left \| \nabla f(x) + s \right \|_2^2, \quad \forall x\in\cx.
\end{equation*}
This completes the proof. \Halmos

\begin{corollary}
    \label{corollary: strong convexity implies KL}
    Consider a convex and compact set $\cx\subseteq\R^n$. If the function $f:\cx \to \R$ is $\mu$-strongly convex over $\cx$, then $f$ satisfies the P{\L}K condition over $\cx$ with constant $\mu$.
\end{corollary}

\proof{Proof of Corollary~\ref{corollary: strong convexity implies KL}}
From the $\mu$-strong convexity of $f$, we have
\begin{equation*}
    f(x) - f^* \le \langle \nabla f(x), x - x^*\rangle - \frac{\mu}{2}\|x - x^*\|_2^2 \le \max_{x'\in\cx} \left\{ \langle \nabla f(x), x - x' \rangle - \frac{\mu}{2} \|x - x'\|_2^2 \right\}
\end{equation*}
Here the last inequality uses the fact that $x^*\in\cx$. Then applying Lemma~\ref{lemma: gradient dominance implies KL} completes the proof. \Halmos

\subsection{No Suboptimal Stationary Points}
\proof{Proof of Proposition~\ref{proposition: KL no local optimal}}
Suppose $\bar{x}\in\cx$ satisfies the first-order necessary optimality condition. We have:
\begin{equation*}
    \begin{aligned}
        \langle \nabla f(\bar{x}), \bar{x} - x \rangle \le 0, \quad \forall x\in\mathcal{X}.
    \end{aligned}
\end{equation*}
Recall that the subdifferential of $\delta_\cx(\bar{x})$ is the normal cone of $\cx$ at $\bar{x}$, i.e.,
\begin{equation*}
    \partial \delta_{\cx}(\bar{x}) = \{g | \langle g, x - \bar{x} \rangle \le 0, \ \forall x\in\cx\}.
\end{equation*}
Therefore, we have $-\nabla f(\bar{x}) \in \partial \delta_\cx(\bar{x})$, which implies
\begin{equation*}
    f(\bar{x}) - f^* \le \frac{1}{2\mu} \min_{g \in \partial \delta_\cx(\bar{x})} \left \| \nabla f(\bar{x}) + g \right \|_2^2 \overset{(a)}{\le} \frac{1}{2\mu} \left \| \nabla f(\bar{x}) - \nabla f(\bar{x}) \right \|_2^2 = 0.
\end{equation*}
Here inequality (a) holds because $-\nabla f(\bar{x}) \in \partial \delta_\cx(\bar{x})$. Since $f(\bar{x}) \ge f^*$ by definition, it follows that $f(\bar{x}) = f^*$ and $\bar{x}$ is a global optimal point. \Halmos

\subsection{Convergence Rate under the P{\L}K Condition}\label{subsection: PGD KL}
\proof{Proof of Lemma~\ref{lemma: PGD KL}} 
\textbf{Projected Gradient Descent:} The non-asymptotic convergence result of projected gradient descent can be found in \citet{attouch2013convergence}. We include the proof for completeness and to keep the paper self-contained. From the optimality condition of $x_{k+1} = \argmin_{x\in\cx}\|x - ( x_k - \gamma_k \nabla f(x_k) ) \|_2^2$, we have
\begin{equation}\label{ineq: optimality projection-PGD}
    \|x_{k+1} - x_k\|_2^2 + \langle x_{k+1} - x_k, \gamma_k \nabla f(x_k) \rangle \le 0.
\end{equation}
From the smoothness of $f$, we conclude that
\begin{equation}
\label{ineq1- PGD KL}
    \begin{aligned}
        f(x_{k+1}) - f(x_k) &\le \langle \nabla f(x_k), x_{k+1} - x_k \rangle + \frac{L}{2}\|x_{k+1} - x_k\|_2^2 \le -\frac{L}{2}\|x_{k+1} - x_k\|_2^2.
    \end{aligned}
\end{equation}
The last inequality uses \eqref{ineq: optimality projection-PGD} and $\gamma_k = \frac{1}{L}$. From the optimality condition of the optimization problem that defines the projection operator,
\begin{equation*}
    x_k - \gamma_k \nabla f(x_k) - x_{k+1} \in \partial \delta_\cx(x_{k+1}).
\end{equation*}
Since the subdifferentials of $\delta_\cx$ form a normal cone, we have
\begin{equation}\label{normal cone}
    \frac{x_k - x_{k+1}}{\gamma_k} - \nabla f(x_k) \in \partial \delta_\cx(x_{k+1}).
\end{equation}
By the P{\L}K condition, we have
\begin{equation*}
    \begin{aligned}
        f(x_{k+1}) - f^* &\le \frac{1}{2\mu} \min_{g \in \partial \delta_\cx(x_{k+1})} \left \| \nabla f(x_{k+1}) + g \right \|_2^2\\
        &\overset{(a)}{\le} \frac{1}{2\mu} \left\| \frac{x_k - x_{k+1}}{\gamma_k} + \nabla f(x_{k+1}) - \nabla f(x_k) \right\|_2^2\\
        &\overset{(b)}{\le} \frac{1}{2\mu} \left( 2\left\| \frac{x_k - x_{k+1}}{\gamma_k} \right\|_2^2 + 2 \left\| \nabla f(x_{k+1}) - \nabla f(x_k) \right\|_2^2 \right)\\
        &\overset{(c)}{\le} \frac{2L^2}{\mu}\|x_{k+1} - x_k\|_2^2.
    \end{aligned}
\end{equation*}
Here, inequality (a) uses \eqref{normal cone}, inequality (b) comes from $\|a + b\|_2^2 \le 2\|a\|_2^2 + 2\|b\|_2^2$, and inequality (c) utilizes $\gamma_k = 1/L$. Then, applying \eqref{ineq1- PGD KL}, we have that
\begin{equation*}
    f(x_{k+1}) - f^* \le \frac{2L^2}{\mu}\|x_{k+1} - x_k\|_2^2 \le \frac{4L}{\mu} \left[ f(x_k) - f(x_{k+1}) \right].
\end{equation*}
Rearranging the terms, we have
\begin{equation*}
    \left( 1 + \frac{4L}{\mu} \right) \left[ f(x_{k+1}) - f^* \right] \le \frac{4L}{\mu} \left[ f(x_k) - f^* \right].
\end{equation*}
This implies that
\begin{equation*}
    f(x_{k}) - f^* \le \left(1 - \frac{\mu}{4L + \mu} \right)^k \left( f(x_0) - f^* \right).
\end{equation*}

\noindent \textbf{Projected Stochastic Gradient Descent:} The proof follows the framework in \citet{attouch2013convergence} with an extension to the stochastic setting. From the optimality condition of $x_{k+1} = \argmin_{x\in\cx}\|x - ( x_k - \gamma_k \nabla \hat{f}(x_k) ) \|_2^2$, we have
\begin{equation}
\label{ineq: optimality projection}
    \|x_{k+1} - x_k\|_2^2 + \langle x_{k+1} - x_k, \gamma_k \nabla \hat{f}(x_k) \rangle \le 0.
\end{equation}
From the smoothness of $f$, we conclude that
\begin{equation}
\label{ineq1: PGD KL}
    \begin{aligned}
        f(x_{k+1}) - f(x_k) &\le \langle \nabla f(x_k), x_{k+1} - x_k \rangle + \frac{L}{2}\|x_{k+1} - x_k\|_2^2\\
        &= \langle \nabla \hat{f}(x_k), x_{k+1} - x_k \rangle + \langle \nabla f(x_k) - \nabla \hat{f}(x_k), x_{k+1} - x_k \rangle + \frac{L}{2}\|x_{k+1} - x_k\|_2^2\\
        &\overset{(a)}{\le} \langle \nabla f(x_k) - \nabla \hat{f}(x_k), x_{k+1} - x_k \rangle - \frac{L}{2}\|x_{k+1} - x_k\|_2^2\\
        &\overset{(b)}{\le} \frac{1}{L}\|\nabla f(x_k) - \nabla \hat{f}(x_k)\|_2^2 + \frac{L}{4}\|x_{k+1} - x_k\|_2^2 - \frac{L}{2}\|x_{k+1} - x_k\|_2^2\\
        &= \frac{1}{L}\|\nabla f(x_k) - \nabla \hat{f}(x_k)\|_2^2 - \frac{L}{4}\|x_{k+1} - x_k\|_2^2.
    \end{aligned}
\end{equation}
Here inequality (a) uses (\ref{ineq: optimality projection}) and $\gamma_k = \frac{1}{L}$, and (b) uses the inequality $2\langle a, b \rangle \le \|a\|_2^2 + \|b\|^2_2$. From the optimality condition of the optimization problem that defines the projection operator,
\begin{equation*}
    x_k - \gamma_k \nabla \hat{f}(x_k) - x_{k+1} \in \partial \delta_\cx(x_{k+1}).
\end{equation*}
Since the subdifferentials of $\delta_\cx$ is a normal cone, we have
\begin{equation*}
    \frac{x_k - x_{k+1}}{\gamma_k} - \nabla \hat{f}(x_k) \in \partial \delta_\cx(x_{k+1}).
\end{equation*}
Thus, we get
\begin{equation*}
    \begin{aligned}
        \min_{g \in \partial \delta_\cx(x_{k+1})} \left \| \nabla f(x_{k+1}) + g \right \|_2 & \le \left\| \frac{x_k - x_{k+1}}{\gamma_k} + \nabla f(x_{k+1}) - \nabla \hat{f}(x_k) \right\|_2\\
        &\overset{(a)}{\le} \left\| \frac{x_k - x_{k+1}}{\gamma_k} \right\|_2 + \|\nabla f(x_{k+1}) - \nabla f(x_k)\|_2 + \|\nabla f(x_k) - \nabla \hat{f}(x_k)\|_2\\
        &\overset{(b)}{\le} 2L \|x_{k+1} - x_k\|_2 + \|\nabla f(x_k) - \nabla \hat{f}(x_k)\|_2.
    \end{aligned}
\end{equation*}
Inequality (a) uses the triangle inequality, and inequality (b) uses the smoothness of $f$ and $\gamma_k = \frac{1}{L}$. By the P{\L}K Condition, we have
\begin{equation*}
    \begin{aligned}
        f(x_{k+1}) - f^* &\le \frac{1}{2\mu} \min_{g \in \partial \delta_\cx(x_{k+1})} \left \| \nabla f(x_{k+1}) + g \right \|_2^2\\
        &\le \frac{1}{2\mu} \left( 2L \|x_{k+1} - x_k\|_2 + \|\nabla f(x_k) - \nabla \hat{f}(x_k)\|_2 \right)^2\\
        &\overset{(a)}{\le} \frac{1}{2\mu} \left( 8L^2 \|x_{k+1} - x_k\|_2^2 + 2 \|\nabla f(x_k) - \nabla \hat{f}(x_k)\|_2^2 \right)\\
        &\overset{(b)}{\le} \frac{16L}{\mu} \left(f(x_k) - f(x_{k+1}) \right) + \frac{17}{\mu} \|\nabla f(x_k) - \nabla \hat{f}(x_k)\|_2^2.
    \end{aligned}
\end{equation*}
Here (a) uses the inequality $(a + b)^2 \le 2a^2 + 2b^2$, and inequality (b) uses (\ref{ineq1: PGD KL}). Notice that $x_k = x_k(\xi_{[k-1]})$ is a function of the $\xi_{[k-1]} = (\xi_1,\dots,\xi_{k-1})$. For simplicity, we will use $\mathbb{E}[f(x_{k+1})]$ and $\mathbb{E}[f(x_k)]$ to denote $\mathbb{E}_{\xi_{[k]}}[f(x_{k+1})]$ and $\mathbb{E}_{\xi_{[k-1]}}[f(x_k)]$ respectively. Then, taking the expectation on both sides and using the assumption that $\E_{\xi_k} \|\nabla \hat{f}(x_k) - \nabla f(x_k)\|_2^2 \le \sigma^2 / N$, we have
\begin{equation*}
    \E[f(x_{k+1})] - f^* \le \frac{16L}{\mu}\left( \E[f(x_k)] - \E[f(x_{k+1})] \right) + \frac{17\sigma^2}{\mu N}.
\end{equation*}
Rearranging the terms, we have
\begin{equation*}
    \E[f(x_{k+1})] - f^* \le \left( 1 - \frac{\mu}{16L + \mu} \right)\left( \E[f(x_k)] - f^* \right) + \frac{17\sigma^2}{(16L + \mu) N}.
\end{equation*}
Taking the telescoping sum, we have
\begin{equation*}
    \begin{aligned}
        \E[f(x_{k})] - f^* &\le \left( 1 - \frac{\mu}{16L + \mu} \right)^k\left( \E[f(x_0)] - f^* \right) + \frac{17\sigma^2}{(16L + \mu) N} \sum_{l=0}^{k-1}\left( 1 - \frac{\mu}{16L + \mu} \right)^l\\
        &\le \left( 1 - \frac{\mu}{16L + \mu} \right)^k\left( \E[f(x_0)] - f^* \right) + \frac{17\sigma^2}{\mu N}.
    \end{aligned}
\end{equation*}
Here the last inequality holds because $\sum_{l=0}^{k-1} (1-q)^l \le 1/q$ for any $q\in(0,1)$. \Halmos

\subsection{P{\L}K Condition in Policy Gradient Formulation}
\label{appendix: KL in PG}
We first prove a technical lemma (Lemma~\ref{lemma: sequence}), which is useful for our main results in Theorem~\ref{theorem: main result}.

\proof{Proof of Lemma~\ref{lemma: sequence}} 
Let us define $u_t^2 \coloneqq \sum_{l=t}^T X_l^2$, $v_t^2 \coloneqq \sum_{l=t}^T Y_l^2$. Without loss of generality, we assume that $v_t > 0$ for any $t\in[T]$. Otherwise, there exists $t \in [T]$ such that $Y_l = 0$ for any $t \le l \le T$, which implies $X_l = Y_l = 0$ for any $t \le l \le T$ by (\ref{ineq: bounded sequential gradient difference}). Therefore, discarding these terms will not affect the final result. 

Furthermore, we define $f_t \coloneqq u_t^2 / v_t^2$ and a series $\{\delta_t\}_{t=1}^T$ where
\begin{equation*}
    \delta_{t} \coloneqq (1 + 2 Y_t M_g + M_g^2 \delta_{t+1} v_{t+1}^2) \delta_{t+1},
\end{equation*}
and $\delta_T = 1$. We use backward mathematical induction to show that $f_t \le \delta_t$ for any $t\in[T]$.

\textbf{Induction Base}: Since $f_T = u_T^2 / v_T^2 = X_T^2 / Y_T^2 = 1$, we have $f_T = \delta_T = 1$.

\textbf{Induction Step}: Suppose that $f_{t+1} \le \delta_{t+1}$ holds. We have
\begin{equation*}
    \begin{aligned}
        f_t & = \frac{u_t^2}{v_t^2} = \frac{u_{t+1}^2 + X_t^2}{v_{t+1}^2 + Y_t^2} \overset{(a)}{\le} \frac{u_{t+1}^2 + (M_g u_{t+1}^2 + Y_t)^2}{v_{t+1}^2 + Y_t^2} = \frac{u_{t+1}^2 + Y_t^2 + M_g^2 u_{t+1}^4 + 2 Y_t M_g u_{t+1}^2}{v_{t+1}^2 + Y_t^2}\\
        & \overset{(b)}{=} \frac{ (1 + 2Y_t M_g + M_g^2 f_{t+1} v_{t+1}^2) f_{t+1} v_{t+1}^2 + Y_t^2}{v_{t+1}^2 + Y_t^2}\\
        & \overset{(c)}{\le} (1 + 2 Y_t M_g + M_g^2 \delta_{t+1} v_{t+1}^2) \delta_{t+1}\\
        & = \delta_t.
    \end{aligned}
\end{equation*}

Here inequality (a) uses (\ref{ineq: bounded sequential gradient difference}), equation (b) comes from the definition $f_{t+1} = u_{t+1}^2 / v_{t+1}^2$, and inequality (c) utilizes the induction assumption $f_{t+1} \le \delta_{t+1}$ and the fact that $\delta_t$ is non-increasing in $t$ with $\delta_t \ge \delta_T = 1$. Therefore, by mathematical induction, we can conclude that $f_t \le \delta_t$ for any $t\in[T]$. By definition, we have the following inequalities,
\begin{equation*}
    \delta_t = \prod_{l=t}^{T-1}(1 + 2 Y_l M_g + M_g^2 \delta_{l+1} v_{l+1}^2) \overset{(a)}{\le} \prod_{l=t}^{T-1}(1 + 2 Y_l M_g + M_g^2 \delta_{t+1} v_t^2), \quad \forall t \in [T-1].
\end{equation*}
Inequality (a) holds because $\delta_t$ and $v_t^2$ are non-increasing in $t$. Consider the following optimization problem:
\begin{equation*}
    \begin{aligned}
        \max_{Y_{t}, \dots, Y_{T-1}} & \quad \prod_{l=t}^{T-1}(1 + 2 Y_l M_g + M_g^2 \delta_{t+1} v_{t}^2),\\
        \text{s.t.} & \quad \sum_{l=t}^{T-1} Y_l^2 \le v_{t}^2,\\
        & \quad Y_l \ge 0, \quad \quad \quad \forall t \le l \le T-1.
    \end{aligned}
\end{equation*}
Since all the terms are non-negative, the optimal solution satisfies the equation $\sum_{l=t}^{T-1} Y_l^2 = v_{t}^2$. Otherwise, we can increase one of $Y_l$ for a larger objective. By the KKT condition \citep{nocedal1999numerical}, the optimal solution has an explicit form $(Y_l^*)^2 = \frac{1}{T - t}v_{t}^2$ for any $t \le l \le T-1$. This implies
\begin{equation}
\label{ineq:delta recursive}
    \delta_{t} \le (1 + 2M_g \frac{1}{\sqrt{T-t}}v_{t} + M_g^2 \delta_{t+1}v_{t}^2)^{T-t}, \quad \forall t\in[T-1].
\end{equation}

Next, we consider the following two cases:
\begin{henumerate}
    \item Suppose that there exists $\tilde{t}\in[T - 1]$ such that $2M_g v_{\tilde{t}} \frac{1}{\sqrt{T-\tilde{t}}} + M_g^2 v_{\tilde{t}}^2 e > \frac{1}{T-\tilde{t}}$. Since $v_t \ge 0$ for all $t\in[T-1]$, we can conclude that
    \begin{equation*}
        v_{\tilde{t}} > \frac{1}{\sqrt{T-\tilde{t}}M_g} \times \frac{\sqrt{e + 1} - 1}{e}.
    \end{equation*}
    Since $v_t$ is non-increasing and $v_1 \ge v_{\tilde{t}}$, it holds that
    \begin{equation*}
        f_1 = \frac{u_1^2}{v_1^2} \le \frac{TG^2}{v_{\tilde{t}}^2} \le 4eM_g^2G^2T^2.
    \end{equation*}
    The first inequality holds as $u_1^2 = \sum_{t=1}^TX_t^2 \le \sum_{t=1}^TG^2 = TG^2$.
    \item Suppose that for any $t\in[T-1]$, we have $2M_gv_t \frac{1}{\sqrt{T-t}} + M_g^2 v_t^2 e \le \frac{1}{T-t}$. We use backward mathematical induction to show $\delta_t \le e$ for all $t\in[T]$.
    
    \textbf{Induction Base}: We have $\delta_T = 1 < e$.

    \textbf{Induction Step}: Suppose that $\delta_{t+1} \le e$ holds. From (\ref{ineq:delta recursive}), we have
    \begin{equation*}
        \begin{aligned}
            \delta_t \le (1 + 2M_g \frac{1}{\sqrt{T-t}}v_t + M_g^2 \delta_{t+1}v_{t}^2)^{T-t} \le (1 + 2M_g \frac{1}{\sqrt{T-t}}v_t + M_g^2 v_{t}^2 e)^{T-t} \le (1 + \frac{1}{T-t})^{T-t} \le e.
        \end{aligned}
    \end{equation*}
    
    Therefore, we have $\delta_t \le e, \forall t\in[T]$ by mathematical induction, which implies $f_1 \le \delta_1 \le e$.
\end{henumerate}

Combining these two cases, we conclude that $f_1 \le \max\{e, 4eM_g^2 G^2 T^2\}$. Since $f_1 = u_1^2 / v_1^2$, we have $\sum_{t=1}^TX_t^2 \le \max\{e, 4eM_g^2G^2T^2\} \sum_{t=1}^TY_t^2$. This completes the proof. \Halmos

The bound in Lemma~\ref{lemma: sequence} depends crucially on $M_g$, $G$, and $T$. This dependence could influence the KL constant in Theorem~\ref{theorem: main result}. In what follows, we demonstrate the tightness of the dependence on $M_g$, $G$, and $T$ in Lemma \ref{lemma: sequence}. Without loss of generality, we consider the case when $M_g>1$, $G>1$, and $M_g G \ge 4$.  

We first construct two sequences $\{X_t\}_{t=1}^T$ and $\{Y_t\}_{t=1}^T$ such that $f_1$ is order of $M_g^2 G^2 T^2$ up to some logarithmic factors. In particular, let us define $Z \coloneqq \log_2(M_g G)$ and
\begin{equation*}
    X_t =
    \begin{cases}
        \frac{2^{Z + 1 - t}}{M_g} & \quad t < \lfloor Z \rfloor, \\
        \frac{Z}{t M_g} & \quad t \ge \lfloor Z \rfloor,
    \end{cases}
    \quad\quad\quad
    Y_t =
    \begin{cases}
        0 & \quad t < T, \\
        \frac{Z}{t M_g} & \quad t = T.
    \end{cases}
\end{equation*}

With the sequences $\{X_t\}_{t=1}^T$ and $\{Y_t\}_{t=1}^T$,  we first verify the inequality \eqref{ineq: bounded sequential gradient difference}. Since $Y_t = 0$ for any $t<T$, the inequality (\ref{ineq: bounded sequential gradient difference}) simplifies to $X_t \le M_g \sum_{k=t+1}^T X_k^2$. Introducing a transformation $\bar{X}_t = M_g X_t$, it is sufficient to check that if the sequence $\{\bar{X}_t\}_{t=1}^T$ satisfies
\begin{equation}
\label{ineq: sequence lower bound}
    \bar{X}_t \le \sum_{k=t+1}^T \bar{X}_k^2,
\end{equation}
with $0 \le \bar{X}_t \le M_g G$ for any $t\in[T]$. To validate (\ref{ineq: sequence lower bound}), it is sufficient to verify that $\bar{X}_t \le \bar{X}_{t+1} + \bar{X}_{t+1}^2$ by mathematical induction. We consider the following three cases:
\begin{henumerate}
    \item For any $t \ge \lfloor Z \rfloor$, we have
    \begin{equation*}
        \bar{X}_{t+1} + \bar{X}_{t+1}^2 - \bar{X}_t = Z \left( \frac{1}{t + 1} - \frac{1}{t} + \frac{ Z }{(t+1)^2} \right) = Z \frac{Zt - t - 1} {t (t+1)^2} \overset{(a)}{\ge} Z \frac{t - 1}{t(t+1)^2} \overset{(b)}{\ge} 0.
    \end{equation*}
    Here inequality (a) holds because $Z = \log_2(M_g G) \ge 2$, and inequality (b) follows since $t \ge 1$.
    \item For $t = \lfloor Z \rfloor - 1$, we have
    \begin{equation*}
        \bar{X}_{t+1} + \bar{X}_{t+1}^2 - \bar{X}_t = \frac{Z}{\lfloor Z \rfloor} + \left( \frac{Z}{\lfloor Z \rfloor} \right)^2 - 2^{Z - \lfloor Z \rfloor} \ge 1 + 1 - 2 = 0.
    \end{equation*}
    The inequality applies $0 \le Z - \lfloor Z \rfloor \le 1$.
    \item  For any $t < \lfloor Z \rfloor - 1$, we obtain
    \begin{equation*}
            \bar{X}_{t+1} + \bar{X}_{t+1}^2 - \bar{X}_t \overset{(a)}{\ge} \bar{X}_{t+1}^2 - \bar{X}_t = (2^{Z - t})^2 - 2^{Z + 1 - t} \overset{(b)}{\ge} 0.
    \end{equation*}
    Here inequality (a) holds as $X_t \ge 0$ for any $t\in[T]$, and inequality (b) applies $Z \ge \lfloor Z \rfloor > t + 1$.
\end{henumerate}

Therefore, the sequences $\{X_t\}_{t=1}^T$ and $\{Y_t\}_{t=1}^T$ satisfy (\ref{ineq: bounded sequential gradient difference}). Additionally, we have
\begin{equation*}
    f_1 = \frac{u_1^2}{v_1^2} \ge \frac{X_1^2}{Y_T^2} = \frac{M_g^2 G^2 T^2}{Z^2} = \frac{M_g^2 G^2 T^2}{\log_2(M_g G)^2}.
\end{equation*}

This example establishes a lower bound of $f_1$, indicating that the dependence on $M_g$, $G$, and $T$ is tight up to some logarithmic factors and the bound in Lemma~\ref{lemma: sequence} is sharp.  Next, we prove the main result (Theorem~\ref{theorem: main result}).

\proof{Proof of Theorem~\ref{theorem: main result}}
The proof mainly follows the proof sketch in Section~\ref{subsection: KL in PG}. For readers' convenience, we divide the proof into several parts.

\textbf{Step 1: Gradient Mismatch Inequality.} For any $g_t \in \partial \delta_{\Theta_t}(\theta_t)$, we have
\begin{equation}
\label{multi-stage KL: step 1}
    \begin{aligned}
        & \ \left\| \nabla_{\theta_t} \E_{s_t \sim \rho_t(\cdot|\pi_\theta)} \left[ Q^{\pi_{\theta^*}}_t \left(s_t, \pi_t(s_t|\theta_t) \right) \right] + g_t\right\|_2 - \left\|\nabla_{\theta_t} \E_{s_t \sim \rho_t(\cdot|\pi_\theta)} \left[ Q^{\pi_\theta}_t \left(s_t, \pi_t(s_t|\theta_t) \right) \right] + g_t \right\|_2\\
        \le & \ \left\| \nabla_{\theta_t} \E_{s_t  \sim \rho_t(\cdot|\pi_\theta)} \left[ Q^{\pi_{\theta^*}}_t \left(s_t, \pi_t(s_t|\theta_t) \right) \right] + g_t - \nabla_{\theta_t} \E_{s_t \sim \rho_t(\cdot|\pi_\theta)} \left[ Q^{\pi_\theta}_t \left(s_t, \pi_t(s_t|\theta_t) \right) \right] - g_t \right\|_2\\
        \overset{(a)}{=} & \ \left\| \nabla_{\theta_t} l(\theta_{[1:t]}, \theta_{[t+1:T]}^*) - \nabla_{\theta_t} l(\theta_{[1:t]}, \theta_{[t+1:T]}) \right\|_2 \\
        \le & \ \sum_{k=t+1}^T \left\| \nabla_{\theta_t} l(\theta_1, \dots, \theta_k, \theta_{k+1}^*, \dots, \theta_T^*) - \nabla_{\theta_t} l \left( \theta_1, \dots, \theta_{k-1}, \theta_k^*, \dots, \theta_T^* \right) \right\|_2\\
        \overset{(b)}{\le} & \ \sum_{k=t+1}^T M_g \left( \mathbb{E}_{s_k \sim \rho_k(\cdot|\pi_\theta)} \left[ Q_k^{\pi_{\theta^*}} \left( s_k, \pi_k(s_k|\theta_k) \right)\right] - \mathbb{E}_{s_k \sim \rho_k(\cdot|\pi_\theta)} \left[Q_k^{\pi_{\theta^*}} \left( s_k, \pi_k(s_k|\theta^*_k) \right) \right] \right)\\
        \overset{(c)}{\le} & \ \sum_{k=t+1}^T \frac{M_g}{2\mu_Q} \min_{g_k \in \partial \delta_{\Theta_k}(\theta_k)}\left\|\nabla_{\theta_k} \E_{s_k \sim \rho_k(\cdot|\pi_\theta)} \left[ Q^{\pi_{\theta^*}}_k \left(s_k, \pi_k(s_k|\theta_k) \right) \right] + g_k \right\|^2_2.
    \end{aligned}
\end{equation}
Here $(a)$ comes from the deterministic policy gradient theorem \citep{silver2014deterministic}, $(b)$ relies on sequential decomposition inequalities, and $(c)$ uses the P{\L}K condition of expected optimal Q-value functions. Define
\begin{equation*}
    \begin{aligned}
        X_t & \coloneqq \min_{g_t \in \partial \delta_{\Theta_t}(\theta_t)}\left\| \nabla_{\theta_t} \E_{s_t \sim \rho_t(\cdot|\pi_\theta)} \left[ Q^{\pi_{\theta^*}}_t \left(s_t, \pi_t(s_t|\theta_t) \right) \right] + g_t \right\|_2,\\
        Y_t & \coloneqq \min_{g_t \in \partial \delta_{\Theta_t}(\theta_t)}\left\|\nabla_{\theta_t} \E_{s_t \sim \rho_t(\cdot|\pi_\theta)} \left[ Q^{\pi_\theta}_t \left(s_t, \pi_t(s_t|\theta_t) \right) \right] + g_t \right\|_2.
    \end{aligned}
\end{equation*}
We have
\begin{equation*}
    \begin{aligned}
        & \ X_t - Y_t\\
        \overset{(a)}{\le} & \ \max_{g_t \in \partial \delta_{\Theta_t}(\theta_t)} \left\{ \left\| \nabla_{\theta_t} \E_{s_t \sim \rho_t(\cdot|\pi_\theta)} \left[ Q^{\pi_{\theta^*}}_t \left(s_t, \pi_t(s_t|\theta_t) \right) \right] + g_t \right\|_2 - \left\|\nabla_{\theta_t} \E_{s_t \sim \rho_t(\cdot|\pi_\theta)} \left[ Q^{\pi_\theta}_t \left(s_t, \pi_t(s_t|\theta_t) \right) \right] + g_t \right\|_2 \right\}\\
        \overset{(b)}{\le} & \ \frac{M_g}{2\mu_Q} \sum_{k=t+1}^T \min_{g_k \in \partial \delta_{\Theta_k}(\theta_k)} \left\| \nabla_{\theta_k} \E_{s_k \sim \rho_k(\cdot|\pi_\theta)} \left[ Q^{\pi_{\theta^*}}_k \left(s_k, \pi_k(s_k|\theta_k) \right) \right] + g_k \right\|_2^2\\
        = & \ \frac{M_g}{2\mu_Q} \sum_{k=t+1}^T X_k^2.
    \end{aligned}
\end{equation*}
Inequality (a) utilizes the fact that $\min_{x\in\cx} f(x) - \min_{x\in\cx} h(x) \le \max_{x\in\cx}\{f(x) - h(x)\}$. Inequality (b) comes from (\ref{multi-stage KL: step 1}). Similarly, we have $Y_t - X_t \le \frac{M_g}{2\mu_Q} \sum_{k=t+1}^T X_k^2$. Given the fact that $X_T = Y_T$ and $X_t, Y_t \le G$ for all $t=1,\dots,T$, applying Lemma~\ref{lemma: sequence}, we have that the bounded gradient mismatch condition holds:
\begin{equation*}
    \begin{aligned}
        \sum_{t=1}^T X_t^2 \le \frac{eM_g^2G^2T^2}{\mu_Q^2}\sum_{t=1}^T Y_t^2.
    \end{aligned}
\end{equation*}

\textbf{Step 2: P{\L}K Condition of $l(\theta)$.} Although \citet{kakade2002approximately} only derived the Performance Difference Lemma for infinite-horizon discounted MDPs, one can use the same trick to get a similar result for finite-horizon MDPs. Therefore, we have
\begin{eqnarray*}
    l(\theta) - l(\theta^*) & = & \sum_{t=1}^T\mathbb{E}_{s_t \sim \rho_t(\cdot|\pi_\theta)} \left[ Q_t^{\pi_{\theta^*}} \left( s_t, \pi_t(s_t|\theta_t) \right) - V_t^{\pi_{\theta^*}}(s_t) \right]\\
    & = & \sum_{t=1}^T\mathbb{E}_{s_t \sim \rho_t(\cdot|\pi_\theta)} \left[ Q_t^{\pi_{\theta^*}} \left( s_t, \pi_t(s_t|\theta_t) \right) - Q_t^{\pi_{\theta^*}} \left( s_t, \pi_t(s_t|\theta^*_t) \right) \right]\\
    & \overset{(a)}{\le} & \sum_{t=1}^T  \frac{1}{2\mu_Q}\min_{g_t \in \partial \delta_{\Theta_t}(\theta_t)} \left\| \nabla_{\theta_t} \E_{s_t \sim \rho_t(\cdot|\pi_\theta)} \left[ Q^{\pi_{\theta^*}}_t \left(s_t, \pi_t(s_t|\theta_t) \right) \right] + g_t \right\|_2^2\\
    & \overset{(b)}{\le} & \frac{eM_g^2G^2T^2}{2\mu_Q^3} \sum_{t=1}^T \min_{g_t \in \partial \delta_{\Theta_t}(\theta_t)} \left\| \nabla_{\theta_t} \E_{s_t\sim \rho_t(\cdot|\pi_\theta)} \left[ Q^{\pi_\theta}_t \left(s_t, \pi_t(s_t|\theta_t) \right) \right] + g_t \right\|_2^2.
\end{eqnarray*}
Inequality (a) utilizes the P{\L}K condition of optimal Q-value functions, and inequality (b) comes from the gradient mismatch inequality. Employing the deterministic policy gradient theorem in \citet{silver2014deterministic},
\begin{equation}
\label{eq: deterministic policy gradient}
    \begin{aligned}
        \nabla_{\theta_t} l(\theta) = \nabla_{\theta_t} \E_{s_t \sim \rho_t(\cdot|\pi_\theta)} \left[ Q^{\pi_\theta}_t \left(s_t, \pi_t(s_t|\theta_t) \right) \right].
    \end{aligned}
\end{equation}
Therefore, we obtain
\begin{equation*}
    \begin{aligned}
        l(\theta) - l(\theta^*) & \le \frac{eM_g^2G^2T^2}{2\mu_Q^3} \sum_{t=1}^T\min_{g_t \in \partial \delta_{\Theta_t}(\theta_t)}\left\|\nabla_{\theta_t} \E_{s_t \sim \rho_t(\cdot|\pi_\theta)} \left[ Q^{\pi_\theta}_t \left(s_t, \pi_t(s_t|\theta_t) \right) \right] + g_t \right\|^2_2\\
        & \overset{(a)}= \frac{eM_g^2G^2T^2}{2\mu_Q^3} \min_{g \in \partial \delta_{\Theta}(\theta)} \left\|\nabla l(\theta) + g \right\|^2_2.
    \end{aligned}
\end{equation*}
The equality $(a)$ comes from (\ref{eq: deterministic policy gradient}) and the fact that
\begin{equation*}
    \begin{aligned}
        \nabla l(\theta) = 
        \begin{bmatrix}
            & \nabla_{\theta_1} \E_{s_1 \sim \rho} \left[ Q^{\pi_\theta}_1 \left(s_1, \pi_1(s_1|\theta_1) \right) \right]\\[8pt]
            & \nabla_{\theta_2} \E_{s_2 \sim \rho_2(\cdot|\pi_\theta)} \left[ Q^{\pi_\theta}_2 \left(s_2, \pi_2(s_2|\theta_2) \right) \right]\\[8pt]
            & \vdots\\[8pt]
            & \nabla_{\theta_T} \E_{s_T \sim \rho_T(\cdot|\pi_\theta)} \left[ Q^{\pi_\theta}_T \left(s_T, \pi_T(s_T|\theta_T) \right) \right]
        \end{bmatrix},\quad
        g = 
        \begin{bmatrix}
            g_1 \\[8pt]
            g_2 \\[8pt]
            \vdots\\[8pt]
            g_T 
        \end{bmatrix}.
    \end{aligned}
\end{equation*}
Therefore, $l(\theta)$ satisfies the P{\L}K condition. This completes the proof. \Halmos

\subsection{P{\L}K Constant under Weaker Assumptions}
\label{subsection: KL with weak assumptions}
In section~\ref{subsection: KL in PG}, we provide a standard approach that can establish a slightly weaker condition than sequential decomposition inequalities in Theorem~\ref{theorem: main result}:
\begin{equation}
\label{weak sequential decomposition inequality}
    \begin{aligned}
        & \ \left\| \nabla_{\theta_t} l(\theta_1, \dots, \theta_k, \theta_{k+1}^*, \dots, \theta_T^*) - \nabla_{\theta_t} l \left( \theta_1, \dots, \theta_{k-1}, \theta_k^*, \dots, \theta_T^* \right) \right\|_2\\
        \le & \ M_g \sqrt{\mathbb{E}_{s_k \sim \rho_k(\cdot|\pi_\theta)} \left[ Q_k^{\pi_{\theta^*}} \left( s_k, \pi_k(s_k|\theta_k) \right)\right] - \mathbb{E}_{s_k \sim \rho_k(\cdot|\pi_\theta)} \left[Q_k^{\pi_{\theta^*}} \left( s_k, \pi_k(s_k|\theta^*_k) \right) \right]}.
    \end{aligned}
\end{equation}

In what follows, we show that if relying on the weaker condition (\ref{weak sequential decomposition inequality}), our analysis leads to a suboptimal characterization of the P{\L}K constant, resulting in an exponential dependence on $T$.

We proceed with the same proof sketch in Theorem~\ref{theorem: main result}. In step 1, we aim to establish the gradient mismatch inequality (\ref{bounded gradient mismatch}). Using the definition of $X_t$ and $Y_t$ and the weaker condition (\ref{weak sequential decomposition inequality}), we have
\begin{equation*}
    \begin{aligned}
        X_t - Y_t \overset{(a)}{\le} & \ \sum_{k=t+1}^T \max_{g_t \in \partial \delta_{\Theta_t}(\theta_t)} \left\| \nabla_{\theta_t} l(\theta_1, \dots, \theta_{k-1}, \theta_k, \theta_{k+1}^*, \dots, \theta_T^*) - \nabla_{\theta_t} l \left( \theta_1, \dots, \theta_{k-1}, \theta_k^*, \theta_{k+1}^*, \dots, \theta_T^* \right) \right\|_2\\
        \overset{(b)}{\le} & \ M_g \sum_{k=t+1}^T \sqrt{\mathbb{E}_{s_k \sim \rho_k(\cdot|\pi_\theta)} \left[ Q_k^{\pi_{\theta^*}} \left( s_k, \pi_k(s_k|\theta_k) \right)\right] - \mathbb{E}_{s_k \sim \rho_k(\cdot|\pi_\theta)} \left[Q_k^{\pi_{\theta^*}} \left( s_k, \pi_k(s_k|\theta^*_k) \right) \right]}\\
        \overset{(c)}{\le} & \ \frac{M_g}{\sqrt{2\mu_Q}} \sum_{k=t+1}^T \min_{g_k \in \partial \delta_{\Theta_k}(\theta_k)} \left\| \nabla_{\theta_k} \E_{s_k \sim \rho_k(\cdot|\pi_\theta)} \left[ Q^{\pi_{\theta^*}}_k \left(s_k, \pi_k(s_k|\theta_k) \right) \right] + g_k \right\|_2\\
        = & \ \frac{M_g}{\sqrt{2\mu_Q}} \sum_{k=t+1}^T X_k.
    \end{aligned}
\end{equation*}
Inequality (a) uses the first three steps in (\ref{multi-stage KL: step 1}). (b) comes from the assumption (\ref{weak sequential decomposition inequality}). (c) uses the P{\L}K condition of expected optimal Q-value functions. 

To proceed, we establish a hard instance when $M_g > 1$ and $G > 1$.
\begin{lemma}
    \label{lemma: weak sequence}
    Assume that the nonnegative sequences $\{X_t\}_{t=1}^T$ and $\{Y_t\}_{t=1}^T$ satisfy
    \begin{equation}
        \label{ineq: weak bounded sequential gradient difference}
        |X_t - Y_t| \le M_g \sum_{k=t+1}^T X_k, \quad X_T = Y_T, \quad X_t, Y_t \le G, \quad \forall t\in[T],
    \end{equation}
    with constants $M_g > 1$ and $G > 1$. Then the best $M$ with $\sum_{t=1}^TX_t^2 \le M \sum_{t=1}^T Y_t^2$ for all sequences $\{X_t\}_{t=1}^T$ and $\{Y_t\}_{t=1}^T$ satisfying (\ref{ineq: weak bounded sequential gradient difference}) cannot be smaller than $M_g^{2(T-1)}$.
\end{lemma}

\proof{Proof of Lemma~\ref{lemma: weak sequence}} The analysis is through construction of a hard instance. Let $X_t = GM_g^{1-t}$ for any $t\in[T]$. Set $Y_T = GM_g^{1-T}$ and $Y_t = 0$ for $1\le t < T$. The two sequences satisfy (\ref{ineq: weak bounded sequential gradient difference}) because 
\begin{equation*}
    |X_t - Y_t| = |X_t| = GM_g^{1-t} = M_g \times GM_g^{-t} = M_gX_{t+1} \le M_g \sum_{k=t+1}^T X_k, \quad \forall t<T,
\end{equation*}
and $|X_T - Y_T| = 0$. Therefore, we have
\begin{equation*}
    \frac{\sum_{t=1}^TX_t^2}{\sum_{t=1}^TY_t^2} = \frac{\sum_{t=1}^TX_t^2}{Y_T^2} \ge \frac{X_1^2}{Y_T^2} = M_g^{2(T-1)}.
\end{equation*}
This concludes the proof. \Halmos

Lemma~\ref{lemma: weak sequence} implies that the constant of the gradient mismatch inequality (\ref{bounded gradient mismatch}) depends at least exponentially on $T$. Following the same steps in the proof of Theorem~\ref{theorem: main result}, the P{\L}K constant admits an exponential dependence on $T$. It remains unclear whether the dependence can be improved under (\ref{weak sequential decomposition inequality}) through alternative analysis. To remove the exponential dependence, we apply the stronger sequential decomposition inequality in Theorem~\ref{theorem: main result} that holds in various applications.

\subsection{P{\L}K Condition vs Gradient Dominance Condition}\label{subsection: PLK vs GD}
We work with the P{\L}K condition because it aligns better with the proof structure in our main result compared to the gradient dominance condition. Consequently, the resulting convergence rate for policy gradient methods exhibits a sharper dependence on the planning horizon $T$.
    
Indeed, we also attempted to carry out the analysis directly using the gradient dominance condition, but the resulting convergence rate for policy gradient methods exhibits a worse dependence on $T$. The key issue comes from the gradient mismatch. To explain this point, recall our approach to analyzing the landscape of $l(\theta)$. If we assume the P{\L}K condition of expected optimal $Q$-value functions with constant $\mu_{\text{P{\L}K}}$,
\begin{equation*}
    \begin{aligned}
        l(\theta) - l(\theta^*) &= \sum_{t=1}^T \left( \mathbb{E}_{s_t \sim \rho_t(\cdot|\pi_\theta)} \left[ Q_t^{\pi_{\theta^*}} \left( s_t, \pi_t(s_t|\theta_t) \right)\right] - \mathbb{E}_{s_t \sim \rho_t(\cdot|\pi_\theta)} \left[Q_t^{\pi_{\theta^*}} \left( s_t, \pi_t(s_t|\theta^*_t) \right) \right] \right)\\
        &\le \sum_{t=1}^T \frac{1}{2\mu_{\text{P{\L}K}}} \min_{g_t \in \partial \delta_{\Theta_t}(\theta_t)} \left\| \nabla_{\theta_t} \mathbb{E}_{s_t \sim \rho_t(\cdot|\pi_\theta)} \left[ Q_t^{\pi_{\theta^*}} \left( s_t, \pi_t(s_t|\theta_t) \right)\right] + g_t \right\|_2^2.
    \end{aligned}
\end{equation*}
To establish the P{\L}K condition of $l(\theta)$, we need to prove a gradient mismatch inequality:
\begin{equation*}
    \begin{aligned}
        & \ \sum_{t=1}^T  \min_{g_t \in \partial \delta_{\Theta_t}(\theta_t)} \left\| \nabla_{\theta_t} \E_{s_t \sim \rho_t(\cdot|\pi_\theta)} \left[ Q^{\pi_{\theta^*}}_t \left(s_t, \pi_t(s_t|\theta_t) \right) \right] + g_t \right\|_2^2\\
        \le & \ M_{\text{P{\L}K}} \sum_{t=1}^T\min_{g_t \in \partial \delta_{\Theta_t}(\theta_t)}\left\|\nabla_{\theta_t} \E_{s_t \sim \rho_t(\cdot|\pi_\theta)} \left[ Q^{\pi_\theta}_t \left(s_t, \pi_t(s_t|\theta_t) \right) \right] + g_t \right\|^2_2.
    \end{aligned}
\end{equation*}
If we choose the $(c, \mu_{\text{GD}})$-gradient dominance condition, a similar gradient mismatch inequality is needed:
\begin{equation*}
    \begin{aligned}
        & \ \sum_{t=1}^T  \max_{\theta_t'\in\Theta_t} \left\{ c \langle \nabla_{\theta_t} \E_{s_t \sim \rho_t(\cdot|\pi_\theta)} \left[ Q^{\pi_{\theta^*}}_t \left(s_t, \pi_t(s_t|\theta_t) \right) \right], \theta_t - \theta_t'\rangle - \frac{\mu_{\text{GD}}}{2}\|\theta_t - \theta_t'\|_2^2 \right\} \\
        \le & \ M_{\text{GD}} \sum_{t=1}^T \max_{\theta_t'\in\Theta_t} \left\{ c \langle \nabla_{\theta_t} \E_{s_t \sim \rho_t(\cdot|\pi_\theta)} \left[ Q^{\pi_\theta}_t \left(s_t, \pi_t(s_t|\theta_t) \right) \right], \theta_t - \theta_t'\rangle - \frac{\mu_{\text{GD}}}{2}\|\theta_t - \theta_t'\|_2^2 \right\}.
    \end{aligned}
\end{equation*}
In our analysis, $M_{\text{GD}}$ scales worse with $T$ than $M_{\text{P{\L}K}}$, leading to a worse $T$-dependence in the convergence rate for policy gradient methods under gradient dominance. For this reason, we adopt the P{\L}K condition.
    
We also note that this distinction disappears in the approach of \citet{bhandari2024global}, which relies on different structural conditions. Specifically, their analysis assumes that \textit{any} expected $Q$-value function satisfies the $(c,\mu_{\text{GD}})$-gradient dominance condition, which is stronger than requiring the P{\L}K condition only for the expected \textit{optimal} $Q$-value function in our paper. Under this stronger structural condition, 
\begin{equation*}
    \begin{aligned}
        l(\theta) - l(\theta^*) &\le \kappa \sum_{t=1}^T \bigl( \mathbb{E}_{s_t \sim \rho_t(\cdot|\pi_\theta)} \left[ Q_t^{\pi_{\theta}} \left( s_t, \pi_t(s_t|\theta_t) \right)\right] - \min_{\theta_t'\in\Theta_t} \mathbb{E}_{s_t \sim \rho_t(\cdot|\pi_\theta)} \left[Q_t^{\pi_{\theta}} \left( s_t, \pi_t(s_t|\theta_t') \right) \right] \bigr)\\
        &\le \kappa \sum_{t=1}^T \max_{\theta_t'\in\Theta_t} \left\{ c \langle \nabla_{\theta_t} \E_{s_t \sim \rho_t(\cdot|\pi_\theta)} \left[ Q^{\pi_\theta}_t \left(s_t, \pi_t(s_t|\theta_t) \right) \right], \theta_t - \theta_t'\rangle - \frac{\mu_{\text{GD}}}{2}\|\theta_t - \theta_t'\|_2^2 \right\}\\
        &= \kappa \sum_{t=1}^T \max_{\theta_t'\in\Theta_t} \left\{ c \langle \nabla_{\theta_t} l(\theta), \theta_t - \theta_t'\rangle - \frac{\mu_{\text{GD}}}{2}\|\theta_t - \theta_t'\|_2^2 \right\}.
    \end{aligned}
\end{equation*}
The first inequality follows from other structural conditions in \citet[Theorem 2]{bhandari2024global}. The last equation avoids the gradient mismatch encountered in our analysis. Consequently, their approach does not distinguish between working under gradient dominance or a P{\L}K condition. Both cases yield the same complexity results. We highlight that requiring \textit{all} expected $Q$-value functions to be gradient dominated is typically violated in operations models, whereas the structural conditions we impose continue to hold.

\section{Omitted Proofs in Section~\ref{section: tabular MDPs}}

\subsection{Feasible Region}\label{subsection: feasible region}
In this subsection, we examine properties of the optimal solution $\theta_t^*$ to restrict the feasible set to $\Theta_t = \{\theta_t \in\R^{m\times n}: \sum_{i=1}^n \theta_t(s_t, i) = 1, \theta_t(s_t,i) \ge \lambda/(n\bar{C}T+n\lambda), \forall s_t\in\cs, \forall i \in \cn\}$. From the Bellman equation (\ref{bellman optimality equation}), for any $s_t\in\cs$, $\theta_t^*(s_t,\cdot)$ minimizes the following function with the constraint $\sum_{i=1}^n \theta_t(s_t, i) = 1$:
\begin{equation*}
    \begin{aligned}
        f_t(\theta_t(s_t, \cdot)) &\coloneqq \lambda\mathcal{R} \left( \theta_t(s_t, \cdot) \right) + \sum_{i\in\cn} \theta_t(s_t, i) \Bigl( C_t(s_t, i) + \sum_{s_{t+1}\in\cs}P_t (s_{t+1}|s_t, i) V_{t+1}^{\pi_{\theta^*}}(s_{t+1}) \Bigr).
    \end{aligned}
\end{equation*}
From the KKT condition \citep{nocedal1999numerical}, there exists $\nu(s_t) \in \mathbb{R}$ such that
\begin{subequations}
    \begin{align}
        & \nabla_{\theta_t(s_t, i)} f_t(\theta_t^*(s_t, i)) + \nu(s_t) = 0, & \forall i \in \cn, \forall s_t \in \cs, \label{eq: complementary} \\
        & \sum_{i=1}^n \theta^*_t(s_t, i) = 1, & \forall s_t \in \cs. \label{eq: primal feasibility}
    \end{align}
\end{subequations}
For simplicity, let us define $q_t(s_t, i) \coloneqq C_t(s_t, i) + \sum_{s_{t+1}\in\cs}P_t (s_{t+1}|s_t, i) V_{t+1}^{\pi_{\theta^*}}(s_{t+1})$. Then, (\ref{eq: complementary}) implies
\begin{equation*}
    \nabla_{\theta_t(s_t, i)} f_t(\theta_t^*) = -\frac{\lambda}{n\theta^*_t(s_t, i)} + q_t(s_t, i) = -\nu(s_t).
\end{equation*}
Therefore,
\begin{equation}\label{eq: optimal theta}
    \theta^*_t(s_t, i) = \frac{\lambda}{n \left( q_t(s_t, i) + \nu(s_t) \right)}.
\end{equation}
Since $\theta_t^*(s_t, i) > 0$, we have $q_t(s_t, i) + \nu(s_t) > 0$ for any $s_t\in\cs$ and $i\in\cn$. By (\ref{eq: primal feasibility}),
\begin{equation}\label{ineq: v upper}
    1 = \sum_{i=1}^n \theta^*_t(s_t, i) = \sum_{i=1}^n \frac{\lambda}{n \left( q_t(s_t, i) + \nu(s_t) \right)} \le \frac{\lambda}{\min_{i\in\cn}q_t(s_t, i) + \nu(s_t)}.
\end{equation}
The inequality uses $q_t(s_t, i) \ge \min_{i\in\cn}q_t(s_t, i)$ for any $i\in\cn$. Plugging (\ref{ineq: v upper}) back to (\ref{eq: optimal theta}),
\begin{equation*}
    \theta^*_t(s_t, i) = \frac{\lambda}{n \left( q_t(s_t, i) + \nu(s_t) \right)} \ge \frac{\lambda}{n(\bar{C}T + \lambda - \min_{i\in\cn}q_t(s_t, i))} \ge \frac{\lambda}{n(\bar{C}T + \lambda)}.
\end{equation*}
The first inequality holds because the accumulated cost is bounded above by $T\bar{C}$. To see this, let $\theta^{U}_t(s_t, i) = 1/n$ for all $s_t\in\cs, i \in \cn$. Then $\mathcal{R} (\theta^{U}_t(s_t,\cdot))=0$, and $V_{t+1}^{\pi_{\theta^*}}(\cdot) \le  V_{t+1}^{\pi_{\theta^{U}}}(\cdot) \le  (T-t)\bar{C}$. The second inequality holds because per-period costs are non-negative. Thus, we can impose the additional constraint $\theta_t(s_t,i)\ge \lambda/(n\bar{C}T+n\lambda)$ in the feasible region.

\subsection{P{\L}K Condition of Expected Optimal Q-Value Function}
\proof{Proof of Lemma~\ref{lemma: tabular MDPs single stage KL}} The continuous differentiability has already been discussed in Section~\ref{subsection: tabular MDPs PLK}. Then, it is sufficient to verify the strong convexity. First, we construct a useful bound. For any $s_t \in \cs$, $t\in[T]$,
\begin{equation*}
    \rho_{t+1}(s_{t+1}|\pi_\theta) = \sum_{s_t\in\cs}\rho_t(s_t|\pi_\theta) \sum_{i\in\cn}\theta_t(s_t, i)P_t(s_{t+1}|s_t, i) \ge \underline{p}\sum_{s_t\in\cs}\rho_t(s_t|\pi_\theta) \sum_{i\in\cn}\theta_t(s_t,i) = \underline{p}.
\end{equation*}
The inequality uses $P_t(s_{t+1}|s_t, i) \ge \underline{p}$. For $t = 1$, we have $\rho_1(s_1) \ge \underline{\rho}$. From the Bellman equations (\ref{bellman optimality equation}),
\begin{equation*}
    \begin{aligned}
        & \ \E_{s_t \sim \rho_t(\cdot|\pi_\theta)} \left[ Q^{\pi_{\theta^*}}_t \left(s_t, \pi_t(s_t|\theta_t) \right) \right]\\
        = & \ \sum_{s_t\in\cs} \rho_t(s_t|\pi_\theta) \Bigl( \lambda\mathcal{R} \left( \theta_t(s_t, \cdot) \right) + \sum_{i\in\cn} \theta_t(s_t, i) \Bigl( C_t(s_t, i) + \sum_{s_{t+1}\in\cs}P_{t} (s_{t+1}|s_t, i) V_{t+1}^{\pi_{\theta^*}}(s_{t+1}) \Bigr) \Bigr).
    \end{aligned}
\end{equation*}
By definition, we can calculate
\begin{equation*}
    \nabla^2_{\theta_t(s_t,\cdot)} \E_{s_t \sim \rho_t(\cdot|\pi_\theta)} \left[ Q^{\pi_{\theta^*}}_t \left(s_t, \pi_t(s_t|\theta_t) \right) \right] = \frac{\lambda}{n} \rho_t(s_t|\pi_\theta) \text{Diag} \left( \frac{1}{\theta_t(s_t, i)^2} \right) \succeq \frac{\lambda}{n}\min\{\underline{\rho}, \underline{p}\} I, \quad \forall s_t\in\cs.
\end{equation*}
This implies the expected optimal Q-value function is $\lambda\min\{\underline{\rho}, \underline{p}\} / n$-strongly convex. By Corollary~\ref{corollary: strong convexity implies KL}, the expected optimal Q-value function satisfies the P{\L}K condition with constant $\lambda\min\{\underline{\rho}, \underline{p}\} / n$. \Halmos

\subsection{Bounded Gradient}
\proof{Proof of Lemma~\ref{lemma: tabular MDPs bounded gradient}} From the Policy Gradient Theorem \citep[Theorem 1]{sutton1999policy}, we have
\begin{equation*}
    \begin{aligned}
        \nabla_{\theta_t(s_t, i)} l(\theta) &= \nabla_{\theta_t(s_t, i)} \E_{s_t \sim \rho_t(\cdot|\pi_\theta)} \left[ Q^{\pi_\theta}_t \left(s_t, \pi_t(s_t | \theta_t) \right) \right]\\
        &= \rho_t(s_t|\pi_\theta) \Bigl( -\underbrace{\frac{\lambda}{n\theta_t(s_t, i)}}_{\text{(I)}} + \underbrace{C_t(s_t, i)}_{\text{(II)}} + \sum_{s_{t+1}\in\cs} P_t(s_{t+1}|s_t, i) \underbrace{V_{t+1}^{\pi_\theta}(s_{t+1})}_{\text{(III)}} \Bigr).
    \end{aligned}
\end{equation*}

Utilizing the assumption that $\theta_t(s_t, i) \ge \underline{\theta}$, we have $|\text{(I)}| \le \lambda / (n\underline{\theta})$. The second term admits $|\text{(II)}| \le \bar{C}$. Lastly, we have the recursive form for (III) using the Bellman equation (\ref{bellman equation}):
\begin{equation*}
    \begin{aligned}
        V_t^{\pi_\theta}(s_t) = C_t^r \left( s_t, \pi_t(s_t|\theta_t) \right) + \sum_{s_{t+1}}P_t \left( s_{t+1}| s_t, \pi_t(s_t|\theta_t) \right) V_{t+1}^{\pi_\theta}(s_{t+1}).
    \end{aligned}
\end{equation*}
By definition, we have
\begin{equation*}
    \begin{aligned}
        \left| C_t^r \left( s_t, \pi_t(s_t|\theta_t) \right) \right| &= \left| C_t \left( s_t, \pi_t(s_t|\theta_t) \right) + \lambda \mathcal{R} \left( \pi_t(s_t|\theta_t) \right) \right|\\
        &\le \sum_{i\in\cn} \theta_t(s_t, i) \left| C_t(s_t, i) \right| + \left| \lambda \mathcal{R} \left( \pi_t(s_t|\theta_t) \right) \right|\\
        &\le \bar{C} + \lambda \log \left( 1 / (n\underline{\theta}) \right).
    \end{aligned}
\end{equation*}
The last inequality uses the assumption that $\pi_t(a_t | s_t) \ge \underline{\theta}$. Therefore, by mathematical induction, we have
\begin{equation*}
    \left| V_t^{\pi_\theta}(s_t) \right| \le (T-t+1) \bar{C} + \lambda (T-t+1) \log \left( 1 / (n\underline{\theta}) \right), \quad \forall s_t \in \cs.
\end{equation*}
Combining all the results, we have
\begin{equation*}
    \begin{aligned}
        \| \nabla_{\theta_t} l(\theta) \|_F^2 &= \sum_{s_t \in \cs, i\in\cn} \left| \nabla_{\theta_t(s_t, i)} l(\theta) \right|^2 \\
        &\le \sum_{s_t \in \cs, i \in \cn} \left[ \rho_t(s_t|\pi_\theta) \Bigl( \frac{\lambda}{n\underline{\theta}} + \bar{C} + \sum_{s_{t+1}} P_t(s_{t+1}| s_t, i) \bigl( (T-t) \bar{C} + \lambda (T-t) \log ( \frac{1}{n\underline{\theta}} ) \bigr) \Bigr) \right]^2 \\
        &\le n \bigl[ 2T \bar{C} + \lambda + \lambda T \log((T\bar{C}+\lambda) / \lambda) \bigr]^2.
    \end{aligned}
\end{equation*}
The last inequality uses $\underline{\theta} = \lambda / (n(T\bar{C} + \lambda))$. This completes the proof. \Halmos

\subsection{Sequential Decomposition Inequality}\label{subsection: tabular MDPs sequential decomposition inequality}

\proof{Proof of Lemma~\ref{lemma: tabular MDPs sequential decomposition}} Let $\theta_\alpha = (\theta_1, \dots, \theta_{k-1}, \theta_k, \theta_{k+1}^*, \dots, \theta_T^*)$ and $\theta_\beta = (\theta_1, \dots, \theta_{k-1}, \theta_k^*, \theta_{k+1}^*, \dots, \theta_T^*)$, then we have the following inequalities:
\begin{equation*}
    \begin{aligned}
        & \ \left\| \nabla_{\theta_t} l(\theta_1, \dots, \theta_{k-1}, \theta_k, \theta_{k+1}^*, \dots, \theta_T^*) - \nabla_{\theta_t} l ( \theta_1, \dots, \theta_{k-1}, \theta_k^*, \theta_{k+1}^*, \dots, \theta_T^* ) \right\|_F\\
        \le & \ \sum_{s_t \in \cs, i \in \cn} \left| \nabla_{\theta_t(s_t, i)} l(\theta_1, \dots, \theta_{k-1}, \theta_k, \theta_{k+1}^*, \dots, \theta_T^*) - \nabla_{\theta_t(s_t, i)} l(\theta_1, \dots, \theta_{k-1}, \theta_k^*, \theta_{k+1}^*, \dots, \theta_T^*) \right|\\
        = & \ \sum_{s_t \in \cs, i \in \cn} \biggl| \rho_t(s_t | \pi_\theta) \sum_{s_{t+1} \in \cs} P_t(s_{t+1} | s_t, i) \left( Q_{t+1}^{\pi_\alpha} \left(s_{t+1}, \pi_{t+1}(s_{t+1}|\theta_{t+1}) \right) - Q_{t+1}^{\pi_\beta} \left(s_{t+1}, \pi_{t+1}(s_{t+1}|\theta_{t+1}) \right) \right) \biggr|.
    \end{aligned}
\end{equation*}

Similarly, we have
\begin{equation}
\label{ineq: tabular MDP recursive}
    \begin{aligned}
        & \ Q_{t+1}^{\pi_\alpha} \left(s_{t+1}, \pi_{t+1}(s_{t+1}|\theta_{t+1}) \right) - Q_{t+1}^{\pi_\beta} \left(s_{t+1}, \pi_{t+1}(s_{t+1}|\theta_{t+1}) \right)\\
        = & \ \sum_{i_{t+1} \in \cn} \theta_{t+1}(s_{t+1}, i_{t+1}) \sum_{s_{t+2} \in \cs} P_{t+1}(s_{t+2} | s_{t+1}, i_{t+1})\\
        & \quad\quad\quad\quad\quad\quad\quad\quad\quad\quad\quad\quad \left( Q_{t+2}^{\pi_\alpha} \left( s_{t+2}, \pi_{t+2}(s_{t+2}|\theta_{t+2}) \right) - Q_{t+2}^{\pi_\beta} \left( s_{t+2}, \pi_{t+2}(s_{t+2}|\theta_{t+2}) \right) \right).
    \end{aligned}
\end{equation}

Applying (\ref{ineq: tabular MDP recursive}) recursively, we conclude that
\begin{equation*}
    \begin{aligned}
        & \ \left\| \nabla_{\theta_t} l(\theta_1, \dots, \theta_{k-1}, \theta_k, \theta_{k+1}^*, \dots, \theta_T^*) - \nabla_{\theta_t} l ( \theta_1, \dots, \theta_{k-1}, \theta_k^*, \theta_{k+1}^*, \dots, \theta_T^* ) \right\|_F\\
        \le & \ \sum_{s_t \in \cs, i\in\cn} \Biggl| \rho_t(s_t | \pi_\theta) \sum_{s_{t+1} \in \cs} P_t(s_{t+1} | s_t, i) \sum_{i_{t+1}\in\cn} \theta_{t+1}(s_{t+1}, i_{t+1}) \sum_{s_{t+2} \in \cs} P_{t+1}(s_{t+2} | s_{t+1}, i_{t+1}) \dots\\
        & \quad \quad \quad  \sum_{i_{k-1}\in\cn} \theta_{k-1}(s_{k-1}, i_{k-1}) \sum_{s_k\in\cs} P_{k-1}(s_k|s_{k-1}, i_{k-1}) \left( Q_k^{\pi_{\theta^*}} \left( s_k, \pi_k(s_k|\theta_k) \right) - Q_k^{\pi_{\theta^*}} \left( s_k, \pi_k(s_k|\theta_k^*) \right) \right) \Biggr|.
    \end{aligned}
\end{equation*}

Since $Q_k^* \left( s_k, \pi_k(s_k|\theta_k) \right) \ge Q_k^* \left( s_k, \pi_k(s_k|\theta_k^*) \right)$, the absolute function can be removed. Then we have
\begin{equation*}
    \begin{aligned}
        & \ \left\| \nabla_{\theta_t} l(\theta_1, \dots, \theta_{k-1}, \theta_k, \theta_{k+1}^*, \dots, \theta_T^*) - \nabla_{\theta_t} l ( \theta_1, \dots, \theta_{k-1}, \theta_k^*, \theta_{k+1}^*, \dots, \theta_T^* ) \right\|_F\\
        \le & \ \sum_{s_t \in \cs, i\in\cn} \rho_t(s_t | \pi_\theta) \sum_{s_{t+1}\in\cs} P_t(s_{t+1} | s_t, i) \sum_{i_{t+1}\in\cn} \theta_{t+1}(s_{t+1}, i_{t+1}) \dots \left( Q_k^{\pi_{\theta^*}} \left( s_k, \pi_k(s_k|\theta_k) \right) - Q_k^{\pi_{\theta^*}} \left( s_k, \pi_k(s_k|\theta_k^*) \right) \right)\\
        \overset{(a)}{\le} & \ \frac{1}{\underline{\theta}} \sum_{s_t \in \cs} \rho_t(s_t | \pi_\theta) \sum_{i\in\cn} \theta_t(s_t, i) \sum_{s_{t+1}\in\cs} P_t(s_{t+1} | s_t, i) \sum_{i_{t+1}\in\cn} \theta_{t+1}(s_{t+1}, i_{t+1}) \sum_{s_{t+2} \in \cs} P_{t+1}(s_{t+2} | s_{t+1}, i_{t+1}) \dots \\
        & \quad \quad \quad  \sum_{i_{k-1}\in\cn} \theta_{k-1}(s_{k-1}, i_{k-1}) \sum_{s_k\in\cs} P_{k-1}(s_k|s_{k-1}, i_{k-1}) \left( Q_k^{\pi_{\theta^*}} \left( s_k, \pi_k(s_k|\theta_k) \right) - Q_k^{\pi_{\theta^*}} \left( s_k, \pi_k(s_k|\theta_k^*) \right) \right)\\
        = & \ \frac{1}{\underline{\theta}} \left( \E_{s_k \sim \rho_k(\cdot|\pi_\theta)} \left[ Q_k^{\pi_{\theta^*}} \left( s_k, \pi_k(s_k|\theta_k) \right) \right] - \E_{s_k \sim \rho_k(\cdot|\pi_\theta)} \left[ Q_k^{\pi_{\theta^*}} \left( s_k, \pi_k(s_k|\theta_k^*) \right) \right] \right).
    \end{aligned}
\end{equation*}

Inequality (a) uses the assumption that $\theta_t(s_t, i) \ge \underline{\theta}$ for any $i \in \cn$. We multiply and divide by $\theta_t(s_t,i_t)$, and since $\theta_t(s_t,i_t)\ge \underline{\theta}$, we have $1 \le \theta_t(s_t,i_t) / \underline{\theta}$. This completes the proof. \Halmos

\section{Omitted Proofs in Section~\ref{section: LQR}}

\subsection{P{\L}K Condition of Expected Optimal Q-Value Function}
\proof{Proof of Lemma~\ref{lemma: LQR single stage KL}} From (\ref{bellman equation}), the Q-value function has the following expression:
\begin{equation*}
    \begin{aligned}
        V_t^{\pi_\theta}(s_t) = Q_t^{\pi_\theta} \left( s_t, \pi_t(s_t|\theta_t) \right) = \underbrace{s_t^\top Q_t s_t + s_t^\top \theta_t^\top R_t \theta_t s_t}_{(\text{I})} + \underbrace{\E_{w_t} \left[ V_{t+1}^{\pi_\theta} \left( (A + B\theta_t)s_t + w_t \right) \right]}_{(\text{II})}.
    \end{aligned}
\end{equation*}

Suppose that $V_{t+1}^{\pi_\theta}(s_{t+1})$ is continuously differentiable in $s_{t+1}$. Term (I) is a quadratic function of $\theta_t$ and, therefore, is continuously differentiable. Term (II) is continuously differentiable since the composition of a continuously differentiable function and a linear function is continuously differentiable. From the induction base $V_T^{\pi_\theta}(s_T) = s_T^\top Q_T s_T$ that is continuously differentiable, we can prove that the Q-value function is continuously differentiable by mathematical induction, which implies the continuous differentiability of the expected Q-value function.

If we plug $\pi^*$ into (II), we get the explicit expression of the optimal Q-value function. \citet{bertsekas1995dynamic} demonstrated the convexity of $V_t^*$ by mathematical induction. Therefore, its composition with linear function $As_t + B(\theta_t s_t) + w_t$ is still convex, which implies the convexity of the term (II). From Assumption~\ref{assumption: LQR} and Lemma~\ref{lemma: positive definite state covariance matrix}, we know that matrices $R_t$ and $\E_{s_t\sim\rho_t(\cdot|\pi_\theta)}[s_t s_t^\top]$ are positive definite. Taking the expectation on the term (I) gives a $2\underline{\sigma}_X\underline{\sigma}_R$-strongly convex function in $\theta_t$ \citep{bhandari2024global}. Combining these, we conclude that the expected optimal Q-value function is $2\underline{\sigma}_X\underline{\sigma}_R$-strongly convex. Leveraging Corollary~\ref{corollary: strong convexity implies KL}, we establish the P{\L}K condition of the expected optimal Q-value function with P{\L}K constant $2\underline{\sigma}_X\underline{\sigma}_R$ for all $t=0,\dots, T-1$. \Halmos

\subsection{Bounded Gradient}
\proof{Proof of Lemma~\ref{lemma: LQR bounded gradient}} For any $\theta\in\Theta$, we derive the following inequality using (\ref{eq: LQR P update}):
\begin{equation*}
    \begin{aligned}
        \|P_t\|_2 &\le \|Q_t\|_2 + \|\theta_t^\top R_t \theta_t\|_2 + \|(A+B\theta_t)^\top P_{t+1} (A+B\theta_t)\|_2\\
        & \overset{(a)}{\le} \bar{\sigma}_Q + \|\theta_t\|_2^2 \|R_t\|_2 + \|A+B\theta_t\|_2^2 \|P_{t+1}\|_2\\
        & \overset{(b)}{\le} \bar{\sigma}_Q + \bar{\sigma}_\Theta^2 \bar{\sigma}_R + \|P_{t+1}\|_2.
    \end{aligned}
\end{equation*}
Here inequality (a) uses the assumption that $\sigma_{\text{max}}(Q_t) \le \bar{\sigma}_Q$ for all $0\le t \le T$, and inequality (b) uses the assumption that $\|\theta_t\|_2 \le \bar{\sigma}_\Theta$ and $\|A + B\theta_t\|_2 \le 1$ for all $0\le t \le T-1$. Since $P_T = Q_T$, we conclude that
\begin{equation}
\label{ineq: LQR P bound}
    \|P_t\|_2 \le (T-t+1)\bar{\sigma}_Q + (T-t)\bar{\sigma}_\Theta^2 \bar{\sigma}_R, \quad \forall t = 0, \dots, T.
\end{equation}
Next, we have the following inequality from (\ref{eq: LQR E update}) for all $t=0,\dots, T-1$:
\begin{equation}
\label{ineq: LQR E bound}
    \begin{aligned}
        \|E_t\|_2 &= \|R_t \theta_t + B^\top P_{t+1} (A + B\theta_t)\|_2\\
        &\le \|R_t \theta_t\|_2 + \|B^\top P_{t+1} (A + B\theta_t)\|_2\\
        &\overset{(a)}{\le} \bar{\sigma}_\Theta \bar{\sigma}_R + \|B\|_2 \|P_{t+1}\|_2 \|A+B\theta_t\|_2\\
        &\overset{(b)}{\le} \bar{\sigma}_\Theta \bar{\sigma}_R + (T-t)\bar{\sigma}_Q\|B\|_2 + (T - t - 1) \bar{\sigma}_\Theta^2 \bar{\sigma}_R\|B\|_2,
    \end{aligned}
\end{equation}
where inequality (a) comes from $\sigma_{\text{max}}(R_t) \le \bar{\sigma}_R$ for any $0\le t \le T-1$ and $\|\theta_t\|_2 \le \bar{\sigma}_\Theta$ for any $\theta_t \in \Theta_t, 0\le t \le T-1$. Inequality (b) uses (\ref{ineq: LQR P bound}) and $\|A + B\theta_t\|_2 \le 1$ for any $\theta_t\in\Theta_t, 0\le t \le T-1$. Recall the linear dynamic function $s_{t+1} = As_t + Ba_t + w_t$, we have the following result:
\begin{equation*}
    \begin{aligned}
        \E[s_{t+1} s_{t+1}^\top] = (A+B\theta_t)\E[s_t s_t^\top](A+B\theta_t)^\top + \E[w_t w_t^\top],
    \end{aligned}
\end{equation*}
since cross terms vanish by independence and zero mean. Therefore, for all $t=0,\dots, T-1$:
\begin{equation*}
    \begin{aligned}
        \|\E[s_{t+1} s_{t+1}^\top]\|_2 &\le \|A+B\theta_t\|_2^2 \|\E[s_t s_t^\top]\|_2 + \|\E[w_t w_t^\top]\|_2\\
        &\overset{(a)}{\le} \|\E[s_t s_t^\top]\|_2 + \|\E[w_t w_t^\top]\|_2\\
        &\overset{(b)}{\le} \|\E[s_t s_t^\top]\|_2 + \bar{\sigma}_W.
    \end{aligned}
\end{equation*}
Here inequality (a) uses $\|A + B\theta_t\|_2 \le 1$ for any $\theta_t \in \Theta_t, 0 \le t \le T-1$ and inequality (b) comes from $\sigma_{\text{max}}(\E[w_t w_t^\top]) \le \bar{\sigma}_W$. Taking the telescoping sum, for any $t=0,\dots,T-1,$, we have
\begin{equation}
\label{ineq: LQR Sigma bound}
    \begin{aligned}
        \|\E[s_{t+1} s_{t+1}^\top]\|_2 &\le \sigma_{\text{max}}(\E[s_0 s_0^\top]) + (t+1)\bar{\sigma}_W\\
        &\le \bar{\sigma}_X + (t+1)\bar{\sigma}_W.
    \end{aligned}
\end{equation}
Thus, combining (\ref{ineq: LQR E bound}) and (\ref{ineq: LQR Sigma bound}), we conclude that
\begin{equation*}
    \begin{aligned}
        \|\nabla_{\theta_t} l(\theta)\|_F &\overset{(a)}{\le} \sqrt{\min\{m, n\}} \|\nabla_{\theta_t} l(\theta)\|_2\\
        &\le 2\sqrt{\min\{m, n\}} \|E_t\|_2 \|\E[s_t s_t^\top]\|_2\\
        &\le 2\sqrt{\min\{m, n\}} (\bar{\sigma}_\Theta \bar{\sigma}_R + (T-t)\bar{\sigma}_Q\|B\|_2 + (T-t-1) \bar{\sigma}_\Theta^2 \bar{\sigma}_R\|B\|_2) (\bar{\sigma}_X + t\bar{\sigma}_W),
    \end{aligned}
\end{equation*}
where (a) uses $\|A\|_F \le \sqrt{r} \|A\|_2$ with $r = \text{rank}(A)$ \citep{golub2013matrix}. The right-hand side of the inequality is polynomial in the model parameters $(m, n, T, \bar{\sigma}_Q, \bar{\sigma}_R, \bar{\sigma}_\Theta, \bar{\sigma}_X, \bar{\sigma}_W, \|B\|_2)$. \Halmos

\subsection{Sequential Decomposition Inequality}
\proof{Proof of Lemma~\ref{lemma: LQR sequential decomposition}} First define
\begin{equation*}
    \Pi_{[j_1:j_2]} \coloneqq (A + B\theta_{j_2})(A + B\theta_{j_2 - 1}) \dots (A + B\theta_{j_1 + 1})(A + B\theta_{j_1}),
\end{equation*}
for $j_1 \le j_2$. Therefore, we have $\|\Pi_{[j_1:j_2]}\|_2 \le 1$ for any $0 \le j_1 \le j_2 \le T-1$ since $\|A+B\theta_t\|_2 \le 1$ for any $\theta_t \in \Theta_t, 0\le t\le T-1$. Recall the gradient formulation in Proposition~\ref{proposition: LQR policy gradient}, we can derive
\begin{equation*}
    \begin{aligned}
        & \nabla_{\theta_t} l(\theta_{[0:k-1]}, \theta_k, \theta^*_{[k+1:T-1]}) - \nabla_{\theta_t} l(\theta_{[0:k-1]}, \theta^*_k, \theta^*_{[k+1:T-1]})\\
        = \ & 2B^\top \left( P_{t+1}(\theta_{[0:k-1]}, \theta_k, \theta^*_{[k+1:T-1]}) - P_{t+1}(\theta_{[0:k-1]}, \theta^*_k, \theta^*_{[k+1:T-1]}) \right) (A + B\theta_t) \E_{s_t\sim\rho_t(\cdot|\pi_\theta)}[s_ts_t^\top]\\
        \overset{(a)}{=} \ & 2B^\top (A + B\theta_{t+1})^\top \bigl( P_{t+2}(\theta_{[0:k-1]}, \theta_k, \theta^*_{[k+1:T-1]})\\
        & \quad \quad \quad \quad \quad \quad \quad - P_{t+2}(\theta_{[0:k-1]}, \theta^*_k, \theta^*_{[k+1:T-1]}) \bigr) (A + B\theta_{t+1}) (A + B\theta_t) \E_{s_t\sim\rho_t(\cdot|\pi_\theta)}[s_ts_t^\top]\\
        &\dots\\
        = \ & 2B^\top \Pi_{[t+1:k-1]}^\top \left( P_k(\theta_{[0:k-1]}, \theta_k, \theta^*_{[k+1:T-1]}) - P_k(\theta_{[0:k-1]}, \theta^*_k, \theta^*_{[k+1:T-1]}) \right) \Pi_{[t:k-1]} \E_{s_t\sim\rho_t(\cdot|\pi_\theta)}[s_ts_t^\top].
    \end{aligned}
\end{equation*}
Equation (a) uses the update (\ref{eq: LQR P update}). Utilizing the explicit expression of the Q-value function in Proposition~\ref{proposition: LQR policy gradient}, we conclude that
\begin{equation*}
    \begin{aligned}
        & \E_{s_k\sim\rho_k(\cdot|\pi_\theta)} \left[ Q_k^{\pi_{\theta^*}} \left( s_k, \pi_k(s_k|\theta_k) \right) \right] - \E_{s_k\sim\rho_k(\cdot|\pi_\theta)} \left[ Q_k^{\pi_{\theta^*}} \left( s_k, \pi_k(s_k|\theta_k^*) \right) \right]\\
        = \ & \E_{s_k\sim\rho_k(\cdot|\pi_\theta)} \left[ s_k^\top \left( P_k(\theta_{[0:k-1]}, \theta_k, \theta^*_{[k+1:T-1]}) - P_k(\theta_{[0:k-1]}, \theta^*_k, \theta^*_{[k+1:T-1]}) \right) s_k \right]\\
        = \ & \text{Tr} \left( \left(P_k(\theta_{[0:k-1]}, \theta_k, \theta^*_{[k+1:T-1]}) - P_k(\theta_{[0:k-1]}, \theta^*_k, \theta^*_{[k+1:T-1]}) \right) \E_{s_k\sim\rho_k(\cdot|\pi_\theta)}[s_ks_k^\top] \right).
    \end{aligned}
\end{equation*}
Let us define $\Sigma_k \coloneqq (\E_{s_k\sim\rho_k(\cdot|\pi_\theta)}[s_ks_k^\top])^{1/2}$. Therefore,
\begin{equation}\label{ineq: LQR SDI1}
    \begin{aligned}
        & \left\| \nabla_{\theta_t} l(\theta_{[0:k-1]}, \theta_k, \theta^*_{[k+1:T-1]}) - \nabla_{\theta_t} l(\theta_{[0:k-1]}, \theta^*_k, \theta^*_{[k+1:T-1]}) \right\|_F\\
        = \ & 2 \left\| B^\top \Pi_{[t+1:k-1]}^\top \left( P_k(\theta_{[0:k-1]}, \theta_k, \theta^*_{[k+1:T-1]}) - P_k(\theta_{[0:k-1]}, \theta^*_k, \theta^*_{[k+1:T-1]}) \right) \Pi_{[t:k-1]} \E_{s_t\sim\rho_t(\cdot|\pi_\theta)}[s_ts_t^\top] \right\|_F\\
        \overset{(a)}{\le} \ & 2 \|B\|_2 \left\| \left( P_k(\theta_{[0:k-1]}, \theta_k, \theta^*_{[k+1:T-1]}) - P_k(\theta_{[0:k-1]}, \theta^*_k, \theta^*_{[k+1:T-1]}) \right) \Sigma_k \right\|_F  \left\| \Sigma_k^{-1} \Pi_{[t:k-1]} \E_{s_t\sim\rho_t(\cdot|\pi_\theta)}[s_ts_t^\top] \right\|_2.\\
        \overset{(b)}{\le} \ & \frac{2\|B\|_2\left\| \E_{s_t\sim\rho_t(\cdot|\pi_\theta)}[s_ts_t^\top] \right\|_2}{\underline{\sigma}_X} \left\| \Sigma_k \left( P_k(\theta_{[0:k-1]}, \theta_k, \theta^*_{[k+1:T-1]}) - P_k(\theta_{[0:k-1]}, \theta^*_k, \theta^*_{[k+1:T-1]}) \right) \Sigma_k \right\|_F.
    \end{aligned}
\end{equation}
Inequality (a) uses $\|\Pi_{[j_1:j_2]}\|_2 \le 1$ for any $0 \le j_1 \le j_2 \le T-1$. Inequality (b) uses $\|\Sigma_k^{-1}\|_2 \le 1/\sqrt{\underline{\sigma}_X}$. Then,
\begin{equation}\label{ineq: LQR SDI2}
    \begin{aligned}
        & \ \left\| \Sigma_k \left( P_k(\theta_{[0:k-1]}, \theta_k, \theta^*_{[k+1:T-1]}) - P_k(\theta_{[0:k-1]}, \theta^*_k, \theta^*_{[k+1:T-1]}) \right) \Sigma_k \right\|_F \\
        \le & \ \text{Tr} \left( \Sigma_k \left( P_k(\theta_{[0:k-1]}, \theta_k, \theta^*_{[k+1:T-1]}) - P_k(\theta_{[0:k-1]}, \theta^*_k, \theta^*_{[k+1:T-1]}) \right) \Sigma_k \right)\\
        = & \ \text{Tr} \left( \left(P_k(\theta_{[0:k-1]}, \theta_k, \theta^*_{[k+1:T-1]}) - P_k(\theta_{[0:k-1]}, \theta^*_k, \theta^*_{[k+1:T-1]}) \right) \E_{s_k\sim\rho_k(\cdot|\pi_\theta)}[s_ks_k^\top] \right)\\
        = & \ \E_{s_k\sim\rho_k(\cdot|\pi_\theta)} \left[ Q_k^{\pi_{\theta^*}} \left( s_k, \pi_k(s_k|\theta_k) \right) \right] - \E_{s_k\sim\rho_k(\cdot|\pi_\theta)} \left[ Q_k^{\pi_{\theta^*}} \left( s_k, \pi_k(s_k|\theta_k^*) \right) \right].
    \end{aligned}
\end{equation}
The inequality holds because $\Sigma_k ( P_k(\theta_{[0:k-1]}, \theta_k, \theta^*_{[k+1:T-1]}) - P_k(\theta_{[0:k-1]}, \theta^*_k, \theta^*_{[k+1:T-1]}) ) \Sigma_k \succeq 0$. Combine (\ref{ineq: LQR SDI1}) and (\ref{ineq: LQR SDI2}), we conclude that
\begin{equation*}
    \begin{aligned}
        & \| \nabla_{\theta_t} l(\theta_{[0:k-1]}, \theta_k, \theta^*_{[k+1:T-1]}) - \nabla_{\theta_t} l(\theta_{[0:k-1]}, \theta^*_k, \theta^*_{[k+1:T-1]}) \|_F\\
        \le \ & \frac{2\|B\|_2 \left\| \E_{s_t\sim\rho_t(\cdot|\pi_\theta)}[s_ts_t^\top] \right\|_2}{\underline{\sigma}_X} \left( \E_{s_k\sim\rho_k(\cdot|\pi_\theta)} \left[ Q_k^{\pi_{\theta^*}} \left( s_k, \pi_k(s_k|\theta_k) \right) \right] - \E_{s_k\sim\rho_k(\cdot|\pi_\theta)} \left[ Q_k^{\pi_{\theta^*}} \left( s_k, \pi_k(s_k|\theta_k^*) \right) \right] \right).
    \end{aligned}
\end{equation*}
By Lemma~\ref{lemma: LQR bounded gradient}, $\left\| \E_{s_t\sim\rho_t(\cdot|\pi_\theta)}[s_ts_t^\top] \right\|_2$ is polynomial in model parameters. This concludes the proof. \Halmos

\section{Omitted Proofs in Section~\ref{section: inventory system}}\label{appendix: inventory models}
In Appendix~\ref{appendix: inventory models} and \ref{appendix: cash balance}, we frequently use the following argument. Let $f:\cx \to\R$ be a continuously differentiable function. Assume that $\xi$ is a random variable whose cumulative distribution function $\P$ is Lipschitz continuous. Consider a function $F(x) \coloneqq \E_{\xi\sim\P}[f(x \wedge \xi)]$. \citet{chen2024efficient} proved that:
\begin{enumerate}
    \item $\frac{\partial}{\partial x} f(x \wedge \xi) \overset{\text{a.s.}}{=} \mathbf{1}(x \le \xi) \times f'(x \wedge \xi) = \mathbf{1}(x \le \xi) \times f'(x)$.
    \item $F'(x) = \frac{\partial}{\partial x}\E_{\xi\sim\P}[f(x \wedge \xi)] = \E_{\xi\sim\P}[\frac{\partial}{\partial x} f(x \wedge \xi)] = \E_{\xi\sim\P}[\mathbf{1}(x \le \xi) \times f'(x)] = \P(x \le \xi) \times f'(x)$.
\end{enumerate}
Similar arguments hold when $F(x) \coloneqq \E_{\xi\sim\P}[f(x \vee \xi)]$. In the following sections, we directly use the results.

\subsection{P{\L}K Condition of Expected Optimal Q-Value Function}

\proof{Proof of Lemma~\ref{lemma: markov modulated demand single stage KL}} We use three parts to complete the proof. First, we demonstrate the relationships between suboptimality gaps $F_t(\theta_t) - F_t(\theta_t^*)$ and $f_t(\theta_{t,i} | i) - f_t(\theta^*_{t,i} | i)$. Next, we show the connections of their gradients. Finally, we prove the P{\L}K property of $F_t$.

\textbf{Step 1: Relationship between suboptimality gaps.} Applying the Bellman equation (\ref{bellman optimality equation}), it holds that
\begin{equation*}
    \begin{aligned}
        & \ \E_{(x_t, i_t) \sim \rho_t(\cdot|\pi_\theta)} \left[ Q^{\pi_{\theta^*}}_t \left(x_t, i_t, \pi_t(x_t, i_t | \theta_t) \right) \right]\\
        = & \ \E_{(x_t, i_t) \sim \rho_t(\cdot|\pi_\theta)} \left[ L_t \left(x_t \vee \theta_{t, i_t} | i_t \right) + \E_{i_{t+1}\sim p(\cdot|i_t), D_t \sim P_D(\cdot | i_t)} \left[ V^{\pi_{\theta^*}}_{t+1} \left( x_t \vee \theta_{t,i_t} - D_t, i_{t+1} \right) \right] \right].
    \end{aligned}
\end{equation*}

Recalling the definition (\ref{equation: markov-modulated function f}), we express expected optimal Q-value functions as
\begin{equation*}
    F_t(\theta_t) = \E_{(x_t, i_t) \sim \rho_t(\cdot|\pi_\theta)} \left[ Q^{\pi_{\theta^*}}_t \left(x_t, i_t, \pi_t(x_t, i_t|\theta_t) \right) \right] = \E_{(x_t, i_t) \sim \rho_t(\cdot|\pi_\theta)} \left[ f_t \left( x_t \vee \theta_{t,i_t} | i_t \right) \right].
\end{equation*}

By the law of total expectation, we rewrite the suboptimality gap:
\begin{equation*}
    \begin{aligned}
        F_t(\theta_t) - F_t(\theta_t^*) &= \E_{(x_t, i_t) \sim \rho_t(\cdot|\pi_\theta)} \left[ f_t \left( x_t \vee \theta_{t,i_t} | i_t \right) - f_t \left( x_t \vee \theta^*_{t,i_t} | i_t \right) \right]\\
        &= \E_{i_t \sim \nu} \left[ \E_{x_t} \left[ f_t \left( x_t \vee \theta_{t, i_t} | i_t \right) - f_t \left( x_t \vee \theta^*_{t,i_t} | i_t \right) \bigl| i_t \right] \right]\\
        &= \sum_{i\in\ci} \nu_i \times \E_{x_t} \left[ f_t \left( x_t \vee \theta_{t,i} | i \right) - f_t \left( x_t \vee \theta^*_{t,i} | i \right) \bigl| i_t = i \right].
    \end{aligned}
\end{equation*}

Without loss of generality, we assume that $\theta_{t, i} \le \theta^*_{t, i}$. For any random variable $\xi$ and its corresponding cumulative distribution function $P(\xi)$, we have
\begin{equation*}
    \begin{aligned}
        & \ \E_{\xi\sim P(\xi)} \left[ f_t \left( \xi \vee \theta_{t,i} | i \right) - f_t \left( \xi \vee \theta^*_{t,i} | i \right) \right]\\
        = & \ \int_{-\infty}^{\theta_{t,i}} \left( f_t \left( \theta_{t,i} | i \right) - f_t \left( \theta^*_{t,i} | i \right) \right) d P(\xi) + \int_{\theta_{t,i}}^{\theta^*_{t,i}} \left( f_t ( \xi | i ) - f_t \left( \theta^*_{t,i} | i \right) \right) d P(\xi)\\
        \overset{(a)}{\le} & \ \int_{-\infty}^{\theta_{t,i}} \left( f_t \left( \theta_{t,i} | i \right) - f_t \left( \theta^*_{t,i} | i \right) \right) d P(\xi) + \int_{\theta_{t,i}}^{\theta^*_{t,i}} \left( f_t \left( \theta_{t,i} | i \right) - f_t \left( \theta^*_{t,i} | i \right) \right) d P(\xi)\\
        \overset{(b)}{\le} & \ f_t \left( \theta_{t,i} | i \right) - f_t \left( \theta^*_{t,i} | i \right).
    \end{aligned}
\end{equation*}

Inequality (a) holds as $\theta^*_{t,i}$ is a minimizer of $f_t(\cdot | i)$, which implies that $f_t(\cdot | i)$ is non-increasing on the interval $[\theta_{t,i}, \theta^*_{t,i}]$ due to its convexity. Inequality (b) comes from the fact that $f_t \left( \theta_{t,i} | i \right) - f_t \left( \theta^*_{t,i} | i \right) \ge 0$. The same result holds when $\theta_{t,i} > \theta^*_{t,i}$. Therefore, we have
\begin{equation*}
    \begin{aligned}
        F_t(\theta_t) - F_t(\theta_t^*) &= \sum_{i\in\ci} \nu_i \times \E_{x_t} \left[ f_t \left( x_t \vee \theta_{t,i} | i \right) - f_t \left( x_t \vee \theta^*_{t,i} | i \right) \bigl| i_t = i \right]\\
        &\le \sum_{i\in\ci} \nu_i \times \left[ f_t \left( \theta_{t,i} | i \right) - f_t \left( \theta^*_{t,i} | i \right) \right].
    \end{aligned}
\end{equation*}

\textbf{Step 2: Relationship between gradients.} By definition, we calculate the partial gradient of $F_t$:
\begin{equation*}
    \begin{aligned}
        \nabla_{\theta_{t,i}} F_t(\theta_t) &= \nabla_{\theta_{t,i}} \E_{(x_t, i_t) \sim \rho_t(\cdot|\pi_\theta)} \left[ f_t \left( x_t \vee \theta_{t,i_t} | i_t \right) \right]\\
        &\overset{(a)}{=} \E_{(x_t, i_t) \sim \rho_t(\cdot|\pi_\theta)} \left[ \mathbf{1}(i_t = i) \times \mathbf{1}(\theta_{t,i} \ge x_t) \times f'_t \left( \theta_{t,i} | i \right) \right]\\
        &= \P \left( i_t = i, \theta_{t,i} \ge x_t \right) \times f'_t \left( \theta_{t,i} | i \right).
    \end{aligned}
\end{equation*}

Here equation (a) utilizes the chain rule. Next, we analyze the property of $\P(i_t = i, \theta_{t,i} \ge x_t)$:
\begin{equation*}
    \begin{aligned}
        \P \left( i_t = i, \theta_{t,i} \ge x_t \right) &\overset{(a)}{=} \E_{i_{t-1}} \left[ \P \left( i_t = i, \theta_{t,i} \ge x_t | i_{t-1} \right) \right]\\
        &\overset{(b)}{=} \E_{i_{t-1}} \left[ \P \left( i_t = i| i_{t-1} \right) \times \P \left( \theta_{t,i} \ge x_t | i_{t-1} \right) \right]\\
        &= \E_{i_{t-1}} \left[ p(i|i_{t-1}) \times \P \left( \theta_{t,i} \ge x_t | i_{t-1} \right) \right].
    \end{aligned}
\end{equation*}

Equation (a) uses the law of total expectation. Equation (b) holds because $x_t$ and $i_t$ are independent conditioned on $i_{t-1}$. For simplicity, we use $\rho^x_{t-1}(\cdot| i_{t-1})$ to denote the CDF of $x_{t-1}$ conditioned on $i_{t-1}$. From the transition kernel $x_t = x_{t-1} \vee \theta_{t-1, i_{t-1}} - D_{t-1}$, we derive the following inequalities:
\begin{equation*}
    \begin{aligned}
        \P \left( \theta_{t,i} \ge x_t | i_{t-1} \right) &= \P \left( \theta_{t,i} \ge x_{t-1} \vee \theta_{t-1, i_{t-1}} - D_{t-1} | i_{t-1} \right)\\
        &= \P \left( \theta_{t,i} \ge x_{t-1} - D_{t-1}, \theta_{t,i} \ge \theta_{t-1, i_{t-1}} - D_{t-1} | i_{t-1} \right)\\
        &= \int_0^{\theta_{t-1,i_{t-1}} - \theta_{t,i}} \P \left( \theta_{t,i} \ge x_{t-1} - D, \theta_{t,i} \ge \theta_{t-1, i_{t-1}} - D | i_{t-1} \right) d P_D(D | i_{t-1})\\
        & \quad\quad\quad\quad + \int_{\theta_{t-1,i_{t-1}} - \theta_{t,i}}^{\infty} \P \left( \theta_{t,i} \ge x_{t-1} - D, \theta_{t,i} \ge \theta_{t-1, i_{t-1}} - D | i_{t-1} \right) d P_D(D | i_{t-1})\\
        &= \int_{\theta_{t-1,i_{t-1}} - \theta_{t,i}}^\infty \rho^x_{t-1} \left( D + \theta_{t,i} | i_{t-1} \right) d P_D(D | i_{t-1}).
    \end{aligned}
\end{equation*}

The last equation holds as $\mathbf{1}( \theta_{t,i} \ge \theta_{t-1, i_{t-1}} - D) = 0$ for any $D \in [0, \theta_{t-1,i_{t-1}} - \theta_{t,i})$. If $ \theta_{t-1, i_{t-1}} < \theta_{t,i}$, the first integral is empty and the second starts at $0$. Since $\rho^x_{t-1} \left( \cdot | i_{t-1} \right)$ is non-negative, we have
\begin{equation*}
    \int_{\theta_{t-1,i_{t-1}} - \theta_{t,i}}^\infty \rho^x_{t-1} \left( D + \theta_{t,i} | i_{t-1} \right) d P_D(D | i_{t-1}) \ge \int_{B - \theta_{t,i}}^\infty \rho^x_{t-1} \left( D + \theta_{t,i} | i_{t-1} \right) d P_D(D | i_{t-1}).
\end{equation*}

From the transition kernel, we know that $x_{t+1} = x_t \vee \theta_{t, i_t} - D_t$. Since $x_1 \in (-\infty, B]$, $\theta_{t, i} \in [0, B]$, and $D_t \in [0, +\infty)$, we have $x_t \le B$ for any $t \in [T]$. Therefore, $\rho^x_{t-1} \left( B | i_{t-1} \right) = 1$, which implies that
\begin{equation*}
    \begin{aligned}
        \P \left( \theta_{t,i} \ge x_t | i_{t-1} \right) &\ge \int_{B - \theta_{t,i}}^\infty \rho^x_{t-1} \left( D + \theta_{t,i} | i_{t-1} \right) d P_D(D | i_{t-1}) \ge \int_{B - \theta_{t,i}}^\infty d P_D(D | i_{t-1}).
    \end{aligned}
\end{equation*}
Therefore, we have
\begin{equation*}
    \P \left( \theta_{t, i} \ge x_t | i_{t-1} \right) \ge \int_{B - \theta_{t,i}}^\infty d P_D(D | i_{t-1}) \overset{(a)}{\ge} 1-P_D(B|i_{t-1}) \overset{(b)}{\ge} \alpha > 0.
\end{equation*}
Inequality (a) holds because $\theta_{t,i} \ge 0$ and inequality (b) comes from Assumption~\ref{assumption: inventory}.\ref{assumption: inventory KL}. For $t = 1$, we have $\P (\theta_{1, i} \ge x_1) \ge \rho(0) \ge \alpha$. Thus,
\begin{equation*}
    \P \left( i_t = i, \theta_{t,i} \ge x_t \right) = \E_{i_{t-1}} \left[ p(i|i_{t-1}) \times \P \left( \theta_{t,i} \ge x_t | i_{t-1} \right) \right] \ge \alpha \E_{i_{t-1}} \left[ p(i|i_{t-1}) \right] = \alpha \nu_i.
\end{equation*}

\textbf{Step 3: P{\L}K condition of $F_t$.} First from step 1, we have
\begin{equation*}
    F_t(\theta_t) - F_t(\theta_t^*) \le \sum_{i\in\ci} \nu_i \times \left[ f_t \left( \theta_{t,i} | i \right) - f_t \left( \theta^*_{t,i} | i \right) \right].
\end{equation*}

Based on Assumption~\ref{assumption: inventory}.\ref{assumption: inventory SC}, we know that the per-period cost $L_t(\cdot|i)$ exhibits $\min_{t\in[T]}\{h_t+b_t\}\mu_D$-strong convexity. Combining this with the convexity of the cost-to-go function $V^*_{t+1}(\cdot, i)$, we conclude that $f_t(\cdot|i)$ also possesses $\min_{t\in[T]}\{h_t+b_t\}\mu_D$-strong convexity for any $i\in\ci$. By Corollary~\ref{corollary: strong convexity implies KL}, we know that strong convexity implies P{\L}K condition. Therefore,
\begin{equation*}
    \begin{aligned}
        F_t(\theta_t) - F_t(\theta_t^*) &\overset{(a)}{\le} \sum_{i\in\ci} \frac{\nu_i}{2\min_{t\in[T]}\{h_t+b_t\}\mu_D} \times \min_{g_{t,i} \in \partial \delta_{[0,B]} (\theta_{t,i}) } \left| f'_t \left( \theta_{t,i} | i \right) + g_{t,i} \right|^2\\
        &\overset{(b)}{=}\sum_{i\in\ci} \frac{\nu_i}{2\min_{t\in[T]}\{h_t+b_t\}\mu_D} \min_{g_{t,i} \in \partial \delta_{[0,B]} (\theta_{t,i}) } \left| \frac{\nabla_{\theta_{t,i}} F_t(\theta_t)}{\P(i_t = i, \theta_{t,i} \ge x_t)} + g_{t,i} \right|^2\\
        &\overset{(c)}{=} \sum_{i\in\ci} \frac{\nu_i}{2\min_{t\in[T]}\{h_t+b_t\}\mu_D\P(i_t = i, \theta_{t,i} \ge x_t)^2} \min_{g_{t,i} \in \partial \delta_{[0,B]} (\theta_{t,i}) } \left| \nabla_{\theta_{t,i}} F_t(\theta_t) + g_{t,i} \right|^2\\
        &\overset{(d)}{\le} \sum_{i\in\ci} \frac{1}{2\min_{t\in[T]}\{h_t+b_t\}\mu_D\nu_i\alpha^2} \min_{g_{t,i} \in \partial \delta_{[0,B]} (\theta_{t,i}) } \left| \nabla_{\theta_{t,i}} F_t(\theta_t) + g_{t,i} \right|^2\\
        &\overset{(e)}{\le} \frac{1}{2\min_{t\in[T]}\{h_t+b_t\}\mu_D \alpha^2 \min_{i\in\ci}\{\nu_i\}} \min_{g_t \in \partial \delta_{\Theta_t}(\theta_t)} \left\| \nabla F_t(\theta_t) + g_t \right\|_2^2.
    \end{aligned}
\end{equation*}

Here inequality (a) utilizes the P{\L}K condition of $f_t(\cdot|i)$. Equation (b) uses the relationship between different gradients in step 2. Equation (c) holds because $\partial \delta_{[0,B]}(\theta_{t,i})$ is a cone. Inequality (d) holds because $\P \left( i_t = i, \theta_{t,i} \ge x_t \right) \ge \alpha \nu_i$. Equation (e) is true since $\nu_i^{-1} \le \min_{i\in\ci}\{\nu_i\}^{-1}$. This completes the proof. \Halmos

\subsection{Gradient Formulation}

\proof{Proof of Proposition~\ref{proposition: Markov modulated demand gradient formulation}}
We first prove the recursive form for partial derivatives of value functions. From the Bellman equation (\ref{bellman equation}), we have:
\begin{equation*}
    \begin{aligned}
        \nabla_x V_t^{\pi_\theta} (x_t, i_t) &= \frac{\partial}{\partial x_t} Q_t^{\pi_\theta} \left( x_t, i_t, \pi_t(x_t, i_t | \theta_t) \right)\\
        &= \frac{\partial}{\partial x_t} \Bigl( L_t \left( x_t \vee \theta_{t,i_t} | i_t \right) + \sum_{i_{t+1}\in\ci} p(i_{t+1} | i_t) \E_{D_t \sim P_D(\cdot | i_t)} \left[ V^{\pi_\theta}_{t+1} \left( x_t \vee \theta_{t,i_t} - D_t, i_{t+1} \right) \right] \Bigr)\\
        &= \mathbf{1} \left( x_t \ge \theta_{t,i_t} \right) \times \Bigl( L'_t( x_t | i_t ) + \sum_{i_{t+1}\in\ci} p(i_{t+1} | i_t) \E_{D_t \sim P_D(\cdot | i_t)} \left[ \nabla_x V_{t+1}^{\pi_\theta}(x_t - D_t, i_{t+1}) \right] \Bigr),
    \end{aligned}
\end{equation*}
where $\nabla_x V_{T+1}^{\pi_\theta}(\cdot, \cdot) = 0$. The last equation uses the chain rule. Then for the policy gradient objective function $l(\theta)$, we calculate its partial derivative by:
\begin{equation*}
    \begin{aligned}
        \frac{\partial}{\partial \theta_{t,i}} l(\theta) \overset{(a)}{=} & \ \E_{(x_t,i_t) \sim \rho_t(\cdot|\pi_\theta)} \left[ \frac{\partial}{\partial \theta_{t,i}} Q^{\pi_\theta}_t \left( x_t, i_t, \pi_t(x_t, i_t | \theta_t) \right) \right]\\
        \overset{(b)}{=} & \ \E_{(x_t, i_t) \sim \rho_t(\cdot|\pi_\theta)} \left[ \frac{\partial}{\partial \theta_{t,i}} \pi_t(x_t, i_t | \theta_t) \times  \frac{\partial}{\partial a_t} Q^{\pi_\theta}_t (x_t, i_t, a_t) \biggl|_{a_t = \pi_t(x_t, i_t | \theta_t) } \right]\\
        \overset{(c)}{=} & \ \E_{(x_t,i_t) \sim \rho_t(\cdot|\pi_\theta)} \Bigl[ \mathbf{1} \left( i_t = i, \theta_{t,i} \ge x_t \right) \times \Bigl( L'_t \left( \theta_{t,i} | i \right) + \sum_{i_{t+1} \in \ci} p(i_{t+1} | i) \E_{D_t \sim P_D(\cdot|i)} \left[ \nabla_x V_{t+1}^{\pi_\theta} \left( \theta_{t,i} - D_t, i_{t+1} \right) \right] \Bigr) \Bigr].
    \end{aligned}
\end{equation*}
Here equation (a) utilizes the Deterministic Policy Gradient Theorem \citep{silver2014deterministic}. Equation (b) applies the chain rule. Equation (c) uses the explicit expression of $\pi_t(x_t,i_t|\theta_t) = \theta_{t,i_t} \vee x_t$ and the Bellman equation (\ref{equation: inventory expected Q value function}). This concludes the proof. \Halmos

\subsection{Bounded Gradient}

\proof{Proof of Lemma~\ref{lemma: markov modulated demand bounded gradient}}
From Proposition~\ref{proposition: Markov modulated demand gradient formulation}, we bound the partial derivative as follows:
\begin{equation*}
    \begin{aligned}
        \left| \frac{\partial}{\partial \theta_{t,i}} l(\theta) \right| &= \Bigl| \P ( i_t = i, \theta_{t,i} \ge x_t ) \times \Bigl( L'_{t,i} ( \theta_{t,i} ) + \sum_{i_{t+1}\in\ci} p(i_{t+1}|i) \E_{D_t\sim P_D(\cdot|i)} \left[ \nabla_x V_{t+1}^{\pi_\theta}(\theta_{t,i} - D_t, i_{t+1}) \right] \Bigr) \Bigr|\\
        &\le \P (i_t = i) \times \Bigl| \Bigl( L'_{t,i} ( \theta_{t,i} ) + \sum_{i_{t+1}\in\ci} p(i_{t+1}|i) \E_{D_t\sim P_D(\cdot|i)} \left[ \nabla_x V_{t+1}^{\pi_\theta}(\theta_{t,i} - D_t, i_{t+1}) \right] \Bigr) \Bigr|\\
        &\overset{(a)}{\le} \nu_i \Bigl( \left| L'_{t,i} ( \theta_{t,i} ) \right| + \sum_{i_{t+1}\in\ci} p(i_{t+1}|i) \left| \E_{D_t\sim P_D(\cdot|i)} \left[ \nabla_x V_{t+1}^{\pi_\theta}(\theta_{t,i} - D_t, i_{t+1}) \right] \right| \Bigr)\\
        &\overset{(b)}{\le} \nu_i \Bigl( \max_{t\in[T]} \left\{ \max\{h_t, b_t\} \right\} + \sum_{i_{t+1}\in\ci} p(i_{t+1}|i) \left| \E_{D_t\sim P_D(\cdot|i)} \left[ \nabla_x V_{t+1}^{\pi_\theta}(\theta_{t,i} - D_t, i_{t+1}) \right] \right| \Bigr).
    \end{aligned}
\end{equation*}

Here inequality (a) employs the triangle inequality and inequality (b) holds because $|L_{t,i}'(\theta_{t,i})| \le \max_{t\in[T]} \left\{ \max\{h_t, b_t\} \right\}$ for any $\theta\in\Theta$, $t\in[T]$ and $i\in\ci$. From (\ref{markov modulated demand: recursive gradient}), we have
\begin{equation}
\label{markov modulated demand: recursive gradient bound}
    \begin{aligned}
        \left| \nabla_x V_t^{\pi_\theta}(x_t, i_t) \right| &= \Bigl| \mathbf{1}(x_t \ge \theta_{t,i} ) \times \Bigl( L'_{t,i}(x_t) + \sum_{i_{t+1}\in\ci} p(i_{t+1}|i_t) \E_{D_t\sim P_D(\cdot|i_t)} \left[ \nabla_x V_{t+1}^{\pi_\theta}(x_t - D_t, i_{t+1}) \right] \Bigr) \Bigr| \\
        &\le \left| L'_{t,i}(x_t) \right| + \sum_{i_{t+1}\in\ci} p(i_{t+1}|i_t) \left| \E_{D_t\sim P_D(\cdot|i_t)} \left[ \nabla_x V_{t+1}^{\pi_\theta}(x_t - D_t, i_{t+1}) \right] \right|\\
        &\le \max_{t\in[T]} \left\{ \max\{h_t, b_t\} \right\} + \sum_{i_{t+1}\in\ci} p(i_{t+1}|i_t) \left| \E_{D_t\sim P_D(\cdot|i_t)} \left[ \nabla_x V_{t+1}^{\pi_\theta}(x_t - D_t, i_{t+1}) \right] \right|.
    \end{aligned}
\end{equation}

We use mathematical induction to prove $| \nabla_x V_t^{\pi_\theta}(x_t, i_t) | \le (T-t+1) \max_{t\in[T]} \{ \max\{h_t, b_t\} \}$ for any $\theta\in\Theta$, $t\in[T]$, $x_t\in(-\infty, B]$, and $i_t\in\ci$.

\textbf{Induction Base:} As $\nabla_x V_{T+1}^{\pi_\theta}(\cdot, \cdot) = 0$, it is obvious that $\left| \nabla_x V_T^{\pi_\theta}(x_t, i_t) \right| \le \max_{t\in[T]} \left\{ \max\{h_t, b_t\} \right\}$.

\textbf{Induction Step:} Suppose we have $| \nabla_x V_{t+1}^{\pi_\theta}(x_{t+1}, i_{t+1}) | \le (T-t) \max_{t\in[T]} \{ \max\{h_t, b_t\} \}$ for any $\theta\in\Theta$, $t\in[T]$, $x_{t+1}\in(-\infty,B]$, and $i_{t+1} \in \ci$, applying (\ref{markov modulated demand: recursive gradient bound}) recursively, we have
\begin{equation*}
    \begin{aligned}
        \left| \nabla_x V_t^{\pi_\theta}(x_t, i_t) \right| &\le \max_{t\in[T]} \left\{ \max\{h_t, b_t\} \right\} + \sum_{i_{t+1}\in\ci} p(i_{t+1}|i_t) \left| \E_{D_t\sim P_D(\cdot|i_t)} \left[ \nabla_x V_{t+1}^{\pi_\theta}(x_t - D_t, i_{t+1}) \right] \right|\\
        &\le \max_{t\in[T]} \left\{ \max\{h_t, b_t\} \right\} + \sum_{i_{t+1}\in\ci} p(i_{t+1}|i_t) (T-t) \max_{t\in[T]} \left\{ \max\{h_t, b_t\} \right\}\\
        &= (T- t + 1) \max_{t\in[T]} \left\{ \max\{h_t, b_t\} \right\}.
    \end{aligned}
\end{equation*}

By mathematical induction, we have $| \nabla_x V_t^{\pi_\theta}(x_t, i_t) | \le (T-t+1) \max_{t\in[T]} \{ \max\{h_t, b_t\} \}$ for any $\theta\in\Theta$, $t\in[T]$, $x_t\in(-\infty, B]$, and $i_t\in\ci$. Therefore, we conclude that
\begin{equation*}
    \begin{aligned}
        \left| \frac{\partial}{\partial \theta_{t,i}} l(\theta) \right| &\le \nu_i \Bigl( \max_{t\in[T]} \left\{ \max\{h_t, b_t\} \right\} + \sum_{i_{t+1}\in\ci} p(i_{t+1}|i) \left| \E_{D_t\sim P_D(\cdot|i)} \left[ \nabla_x V_{t+1}^{\pi_\theta}(\theta_{t,i} - D_t, i_{t+1}) \right] \right| \Bigr)\\
        &\le \nu_i (T-t+1)\max_{t\in[T]} \left\{ \max\{h_t, b_t\} \right\}\\
        &\le \nu_i T\max_{t\in[T]} \left\{ \max\{h_t, b_t\} \right\}
    \end{aligned}
\end{equation*}

Thus, we obtain
\begin{equation*}
    \left\| \nabla_{\theta_t} l(\theta) \right\|_2 \le \left\| \nabla_{\theta_t} l(\theta) \right\|_1 \le \sum_{i\in\ci} \nu_i T\max_{t\in[T]} \left\{ \max\{h_t, b_t\} \right\} = T\max_{t\in[T]} \left\{ \max\{h_t, b_t\} \right\}.
\end{equation*}

This completes the proof. \Halmos

\subsection{Sequential Decomposition Inequality}

\proof{Proof of Lemma~\ref{lemma: markov modulated demand bounded sequence gradient difference}} For simplicity, we define $\theta_\alpha = (\theta_{[1:k]}, \theta^*_{[k+1:T]})$ and $\theta_\beta = (\theta_{[1:k-1]}, \theta^*_{[k:T]})$. Furthermore, we denote $\pi_\alpha$ and $\pi_\beta$ as the policies deploying parameters $\theta_\alpha$ and $\theta_\beta$, respectively. Similar to the previous discussion, let $\pi_\theta$ denote the policy using parameters $\theta = (\theta_1, \dots, \theta_T)$. Then
\begin{equation*}
    \begin{aligned}
        \left\| \nabla_{\theta_t} l(\theta_\alpha) - \nabla_{\theta_t} l(\theta_\beta) \right\|_2 \le \left\| \nabla_{\theta_t} l(\theta_\alpha) - \nabla_{\theta_t} l(\theta_\beta) \right\|_1 = \sum_{i\in\ci} \left| \frac{\partial}{\partial \theta_{t,i}} l(\theta_\alpha) - \frac{\partial}{\partial \theta_{t,i}} l(\theta_\beta) \right|.
    \end{aligned}
\end{equation*}

For any $i\in\ci$, we can derive the following inequalities by Proposition~\ref{proposition: Markov modulated demand gradient formulation}:
\begin{equation*}
    \begin{aligned}
        & \ \left| \frac{\partial}{\partial \theta_{t,i}} l(\theta_\alpha) - \frac{\partial}{\partial \theta_{t,i}} l(\theta_\beta) \right|\\
        = & \ \Biggl| \E_{(x_t, i_t)\sim\rho_t(\cdot|\pi_\theta)} \biggl[ \mathbf{1}(i_t = i, \theta_{t, i} \ge x_t)\\
        & \quad \times \sum_{i_{t+1}\in\ci}p(i_{t+1}|i) \left( \E_{D_t\sim P_D(\cdot|i)} \left[\nabla_x V_{t+1}^{\pi_\alpha}(\theta_{t,i} - D_t, i_{t+1}) \right] - \E_{D_t\sim P_D(\cdot|i)} \left[ \nabla_x V_{t+1}^{\pi_\beta}(\theta_{t,i} - D_t, i_{t+1}) \right] \right) \biggr] \Biggr|\\
        \le & \ \nu_i \sum_{i_{t+1} \in \ci} p(i_{t+1}|i) \left| \E_{D_t\sim P_D(\cdot|i)} \left[\nabla_x V_{t+1}^{\pi_\alpha}(\theta_{t,i} - D_t, i_{t+1}) \right] - \E_{D_t\sim P_D(\cdot|i)} \left[ \nabla_x V_{t+1}^{\pi_\beta}(\theta_{t,i} - D_t, i_{t+1}) \right] \right|\\
        \overset{(a)}{\le} & \ \nu_i \sum_{i_{t+1} \in \ci} p(i_{t+1}|i) \sum_{i_{t+2}\in\ci} p(i_{t+2}|i_{t+1})\\
        & \quad \times \left| \E_{D_{[t:t+1]}} \left[\nabla_x V_{t+2}^{\pi_\alpha}(\theta_{t,i} - D_{[t:t+1]}, i_{t+2}) \right] - \E_{D_{[t:t+1]}} \left[ \nabla_x V_{t+2}^{\pi_\beta}(\theta_{t,i} - D_{[t:t+1]}, i_{t+2}) \right] \right|\\
        \le & \ \nu_i \sum_{i_{t+1} \in \ci} p(i_{t+1}|i) \sum_{i_{t+2}\in\ci} p(i_{t+2}|i_{t+1}) \sum_{i_{t+3}\in\ci} p(i_{t+3}|i_{t+2})\\
        & \quad \times \left| \E_{D_{[t:t+2]}} \left[\nabla_x V_{t+3}^{\pi_\alpha}(\theta_{t,i} - D_{[t:t+2]}, i_{t+3}) \right] - \E_{D_{[t:t+2]}} \left[ \nabla_x V_{t+3}^{\pi_\beta}(\theta_{t,i} - D_{[t:t+2]}, i_{t+3}) \right] \right|\\
        & \dots\\
        \overset{(b)}{\le} & \ \nu_i \sum_{i_{t+1} \in \ci} p(i_{t+1}|i) \sum_{i_{t+2}\in\ci} p(i_{t+2}|i_{t+1})\dots\sum_{i_k\in\ci}p(i_k|i_{k-1})\\
        & \quad \times \left| \E_{D_{[t:k-1]}} \left[\nabla_x V_{k}^{\pi_\alpha}(\theta_{t,i} - D_{[t:k-1]}, i_k) \right] - \E_{D_{[t:k-1]}} \left[ \nabla_x V_k^{\pi_\beta}(\theta_{t,i} - D_{[t:k-1]}, i_k) \right] \right|.
    \end{aligned}
\end{equation*}

Here inequality (a) applies (\ref{markov modulated demand: recursive gradient}) and utilizes $\mathbf{1}(\theta_{t,i} - D_t \ge \theta_{t+1, i_{t+1}}) \le 1$. Inequality (b) holds by applying (\ref{markov modulated demand: recursive gradient}) recursively. Recall the definition of $\theta_\alpha = (\theta_{[1:k]}, \theta^*_{[k+1:T]})$, $\theta_\beta = (\theta_{[1:k-1]}, \theta^*_{[k:T]})$, and $f_t(\cdot|i)$. For any sample path $(i, i_{t+1}, \dots, i_k)$, we have
\begin{equation*}
    \left\{
    \begin{aligned}
        \nabla_x V_k^{\pi_\alpha}(\theta_{t,i} - D_{[t:k-1]}, i_k) &= \mathbf{1}(\theta_{t,i} - D_{[t:k-1]} \ge \theta_{k,i_k}) \times f_{k}'(\theta_{t,i} - D_{[t:k-1]} | i_k),\\
        \nabla_x V_k^{\pi_\beta}(\theta_{t,i} - D_{[t:k-1]}, i_k) &= \mathbf{1}(\theta_{t,i} - D_{[t:k-1]} \ge \theta^*_{k,i_k}) \times f_{k}'(\theta_{t,i} - D_{[t:k-1]} | i_k).
    \end{aligned}
    \right.
\end{equation*}

Without loss of generality, we assume that $\theta_{k,i_k} \le \theta^*_{k,i_k}$. Then for any sample path $(i, i_{t+1}, \dots, i_k)$, $D_t\sim F_D(\cdot|i)$, and $D_j \sim F_D(\cdot|i_j), \forall t + 1\le j \le k-1$, we have
\begin{equation*}
    \begin{aligned}
        & \ \left| \E_{D_{[t:k-1]}} \left[\nabla_x V_{k}^{\pi_\alpha}(\theta_{t,i} - D_{[t:k-1]}, i_k) \right] - \E_{D_{[t:k-1]}} \left[ \nabla_x V_k^{\pi_\beta}(\theta_{t,i} - D_{[t:k-1]}, i_k) \right] \right|\\
        \le & \ \E_{D_{[t:k-1]}} \left[ \left| \mathbf{1}(\theta_{k, i_k} \le \theta_{t,i} - D_{[t:k-1]} \le \theta^*_{k,i_k}) \times f_{k}'(\theta_{t,i} - D_{[t:k-1]} | i_k) \right| \right]\\
        \overset{(a)}{=} & \ \int_{\theta_{k, i_k}}^{\theta_{k, i_k}^*} - f_{k}' (x | i_k) \psi(x) dx,
    \end{aligned}
\end{equation*}
where $\psi$ is the probability density function of $x \coloneqq \theta_{t,i} - D_{[t:k-1]}$. Equation (a) holds because $f_{k}(\cdot | i_k)$ is convex, therefore $f_{k}' (\cdot | i_k) \le 0$ within the interval $[\theta_{k, i_k}, \theta_{k, i_k}^*]$. From Assumption~\ref{assumption: inventory}.\ref{assumption: inventory Lipschitz of rv}, the cumulative distribution function of the random demand $D_t$ is $L_{D}$-Lipschitz continuous. Then the probability density function of $D_t$ is upper bounded by $L_D$. Suppose $\psi_{D_1}(\cdot)$ and $\psi_{D_2}(\cdot)$ are probability density functions of $D_1$ and $D_2$, we have the following inequalities:
\begin{equation*}
    \psi_{D_1+D_2}(\omega) = \int_0^\omega \psi_{D_1}(\nu)\psi_{D_2}(\omega - \nu) d\nu \le L_D \int_0^\omega \psi_{D_1}(\nu) d \nu \le L_D,
\end{equation*}
which implies that the probability density function of cumulative demands is upper bounded by $L_D$ and thus $\psi(\cdot) \le L_D$. Hence,
\begin{equation*}
    \begin{aligned}
        \int_{\theta_{k, i_k}}^{\theta_{k, i_k}^*} - f_{k}' (x | i_k) \psi(x) dx \le \ & L_D \int_{\theta_{k,i_k}}^{\theta_{k,i_k}^*} - f_{k}' (x | i_k) dx = L_D \left( f_{k} (\theta_{k,i_k} | i_k) - f_{k} (\theta_{k,i_k}^* | i_k) \right).
    \end{aligned}
\end{equation*}

The same result holds when $\theta_k > \theta_k^*$ following a similar derivation. Therefore,
\begin{equation*}
    \begin{aligned}
        \left\| \nabla_{\theta_t} l(\theta_\alpha) - \nabla_{\theta_t} l(\theta_\beta) \right\|_2 &\le \sum_{i\in\ci} \left| \frac{\partial}{\partial \theta_{t,i}} l(\theta_\alpha) - \frac{\partial}{\partial \theta_{t,i}} l(\theta_\beta) \right|\\
        &\le L_D \sum_{i\in\ci}\nu_i \sum_{i_{t+1}\in\ci} p(i_{t+1}|i) \dots \sum_{i_k\in\ci} p(i_k|i_{k-1}) \left( f_{k} (\theta_{k,i_k} | i_k) - f_{k} (\theta_{k,i_k}^* | i_k) \right)\\
        &= L_D\sum_{i\in\ci} \nu_i \left( f_{k} (\theta_{k,i} | i) - f_{k} (\theta_{k,i}^* | i) \right).
    \end{aligned}
\end{equation*}
The last equation is true as $\nu$ is a stationary distribution of the exogenous Markov chain. Recalling the definition of $F_k(\theta_k)$, we have
\begin{equation*}
    \begin{aligned}
        F_k(\theta_k) - F_k(\theta_k^*) &= \sum_{i\in\ci} \nu_i \times \E_{x_k} \left[ f_k \left( x_k \vee \theta_{k,i} | i \right) - f_k \left( x_k \vee \theta^*_{k,i} | i \right) \bigl| i_k = i \right].
    \end{aligned}
\end{equation*}

Without loss of generality, we assume that $\theta_{k,i} \le \theta_{k, i}^*$. For any random variable $\xi$ and its corresponding cumulative distribution function $P(\xi)$, we have
\begin{equation*}
    \begin{aligned}
        \E_{\xi\sim P(\xi)} \left[ f_k \left( \xi \vee \theta_{k,i} | i \right) - f_k \left( \xi \vee \theta^*_{k,i} | i \right) \right] &= \int_{-\infty}^{\theta_{k,i}} f_k \left( \theta_{k,i} | i \right) - f_k \left( \theta^*_{k,i} | i \right) d P(\xi) + \int_{\theta_{k,i}}^{\theta^*_{k,i}} f_k ( \xi | i ) - f_k \left( \theta^*_{k,i} | i \right) d P(\xi)\\
        & \ge \int_{-\infty}^{\theta_{k,i}} f_k \left( \theta_{k,i} | i \right) - f_k \left( \theta^*_{k,i} | i \right) d P(\xi).
    \end{aligned}
\end{equation*}
The last inequality holds as $f_k(\xi | i) \ge f_k(\theta_{k,i}^* | i)$ for any $\xi \in [\theta_{k,i}, \theta_{k, i}^*]$. Therefore, we have
\begin{equation*}
    \begin{aligned}
        F_k(\theta_k) - F_k(\theta_k^*) &\ge \sum_{i\in\ci} \nu_i \int_{-\infty}^{\theta_{k,i}} d P(x_k|i_k = i) \left( f_k \left( \theta_{k,i} | i \right) - f_k \left( \theta^*_{k,i} | i \right) \right)\\
        &= \sum_{i \in \ci} \P(x_k \le \theta_{k,i}, i_k = i) \left( f_k \left( \theta_{k,i} | i \right) - f_k \left( \theta^*_{k,i} | i \right) \right).
    \end{aligned}
\end{equation*}

Similar results hold when $\theta_{k,i} > \theta_{k, i}^*$. Following the same procedure in the proof of Lemma~\ref{lemma: markov modulated demand single stage KL} step 2, we have $\P(x_k \le \theta_{k,i}, i_k = i) \ge \alpha \nu_i$. Thus, we conclude that
\begin{equation*}
    \begin{aligned}
        \left\| \nabla_{\theta_t} l(\theta_\alpha) - \nabla_{\theta_t} l(\theta_\beta) \right\|_2 &\le L_D\sum_{i\in\ci} \nu_i \left( f_{k} (\theta_{k,i} | i) - f_{k} (\theta_{k,i}^* | i) \right)\\
        &\le \frac{L_D}{\alpha} \sum_{i \in \ci} \P(x_k \le \theta_{k,i}, i_k = i) \left( f_k \left( \theta_{k,i} | i \right) - f_k \left( \theta^*_{k,i} | i \right) \right)\\
        &\le \frac{L_D}{\alpha} \left( F_k(\theta_k) - F_k(\theta_k^*) \right).
    \end{aligned}
\end{equation*}

This completes the proof. \Halmos

\subsection{Smoothness}\label{subsection: sample complexity}
In this section, we verify the smoothness of $l(\theta)$. Let $\Gamma = (x_1, i_1, D_1, i_2, D_2, \dots,i_T, D_T)$ denote the sample path generated by the exogenous Markov chain and demand process.

\begin{lemma}[Smoothness]
\label{lemma: lipschitz_pg invt control}
    Suppose that Assumption \ref{assumption: inventory} holds. The policy gradient $\nabla l(\theta)$ is $S_l$-lipschitz continuous, where $S_l = 2\max_{t\in[T]}\{h_t + b_t\}\max \{L_D, L_\rho\}T^2\sqrt{|\ci|}$.
\end{lemma}

The remaining part of this subsection is to prove Lemma~\ref{lemma: lipschitz_pg invt control}, which relies on several technical lemmas. We first define a random variable $\Lambda_t(\theta, \Gamma)$ representing the last time (after period $t$) that the on-hand inventory level is greater than the base-stock level in all subsequent periods. More specifically,
\begin{equation*}
    \Lambda_t(\theta, \Gamma) \coloneqq \max \left\{ t' | t \le t' \le T, \theta_{t, i_t} - D_{[t:l-1]} \ge \theta_{l, i_l}, \forall t \le l \le t' \right\}.
\end{equation*}
By definition, the event $j \le \Lambda_t(\theta, \Gamma)$ is equivalent to the event $\theta_{t, i_t} - D_{[t:l-1]} \ge \theta_{l, i_l}$, $\forall t \le l \le j$. For any $\bar{\theta} \in \Theta$, let $\bar{x}_t$ denote the inventory level at the beginning of period $t$ when using the policy $\pi_{\bar{\theta}}$. For simplicity of notation, define $\hat{L}_t(x_t|i_t) \coloneqq h_t\mathbf{1}(x_t \ge D_t) - b_t\mathbf{1}(x_t < D_t)$. By Proposition~\ref{proposition: Markov modulated demand gradient formulation}, we obtain the following:
\begin{equation*}
    \begin{aligned}
        & \ \left\| \nabla l(\theta) - \nabla l(\bar{\theta}) \right\|_2\\
        \le & \ \sum_{t\in[T], i\in\ci} \left| \frac{\partial}{\partial \theta_{t,i}}l(\theta) - \frac{\partial}{\partial \theta_{t,i}}l(\bar{\theta}) \right|\\
        = & \ \sum_{t\in[T], i\in\ci} \Biggl| \E_\Gamma \Biggl[ \sum_{j=t}^T \biggl( \hat{L}_j(\theta_{t,i} - D_{[t:j-1]}|i_j) \mathbf{1} \left( i_t = i, \theta_{t, i} \ge x_t, j \le \Lambda_t(\theta, \Gamma) \right)\\
        &\quad\quad\quad\quad\quad\quad\quad - \hat{L}_j(\bar{\theta}_{t,i} - D_{[t:j-1]}|i_j) \mathbf{1} \left( i_t = i, \bar{\theta}_{t, i} \ge \bar{x}_t, j \le \Lambda_t(\bar{\theta}, \Gamma) \right) \biggr) \Biggr] \Biggr|\\
        \le & \ \sum_{t\in[T], i\in\ci} \sum_{j=t}^T \Biggl( \underbrace{\E_\Gamma \left[ \left| \left( \hat{L}_j(\theta_{t,i} - D_{[t:j-1]}|i_j) - \hat{L}_j(\bar{\theta}_{t,i_t} - D_{[t:j-1]}|i_j) \right) \mathbf{1} \left( i_t = i, \theta_{t, i} \ge x_t, j \le \Lambda_t(\theta, \Gamma) \right) \right| \right]}_{\text{(I)}}\\
        & \quad\quad\quad\quad\quad + \underbrace{\E_\Gamma \left[ \left| \hat{L}_j(\bar{\theta}_{t,i} - D_{[t:j-1]}|i_j) \mathbf{1} \left( i_t = i, j \le \Lambda_t(\theta, \Gamma) \right) \left( \mathbf{1}(\theta_{t, i} \ge x_t) - \mathbf{1} (\bar{\theta}_{t, i} \ge \bar{x}_t) \right) \right| \right]}_{\text{(II)}}\\
        & \quad\quad\quad\quad\quad + \underbrace{\E_\Gamma \left[ \left| \hat{L}_j(\bar{\theta}_{t,i} - D_{[t:j-1]}|i_j) \mathbf{1} ( i_t = i, \bar{\theta}_{t, i} \ge \bar{x}_t ) \left( \mathbf{1} \left( j \le \Lambda_t(\theta, \Gamma) \right) - \mathbf{1} \left( j \le \Lambda_t(\bar{\theta}, \Gamma) \right) \right) \right| \right]}_{\text{(III)}} \Biggr).
    \end{aligned}
\end{equation*}
The last inequality applies the triangle inequality. Further upper bounding the right hand side uses the following lemma.

\begin{lemma}[Lipschitz Continuity under Convolution]
    \label{lemma:lip_preserve}
    Suppose $F_{\cz}$ is the cumulative distribution function of the random variable $Z = X + Y$, where $X$ and $Y$ are two independent random variables from the distributions $F_{\cx}$ and $F_{\cy}$. Assume that $F_{\cx}$ and $F_{\cy}$ are $L_{\cx}$ and $L_{\cy}$-Lipschitz continuous respectively. Then $F_{\cz}$ is $\min\{L_{\cx}, L_{\cy}\}$-Lipschitz continuous.
\end{lemma}
\proof{Proof of Lemma \ref{lemma:lip_preserve}}
Denote $\cx, \cy, \cz$ as the support of $X, Y$, and $Z$ respectively. For any $z_1, z_2\in\cz$, we have
\begin{equation*}
    \begin{aligned}
        \left|F_{\cz}(z_1) - F_{\cz}(z_2)\right| &= \left|\int_{\cx} F_{\cy}(z_1 - X) dF_{\cx}(X) - \int_{\cx} F_{\cy}(z_2 - X) dF_{\cx}(X)\right|\\
        & \le \int_{\cx} \left| F_{\cy}(z_1 - X) - F_{\cy}(z_2 - X) \right| dF_{\cx}(X)\\
        & \overset{(a)}{\le} L_{\cy}|z_1 - z_2|,
    \end{aligned}
\end{equation*}
where (a) comes from the $L_{\cy}$-Lipschitz continuity of $F_{\cy}$. Similarly, we can derive
\begin{equation*}
    \begin{aligned}
        \left|F_{\cz}(z_1) - F_{\cz}(z_2)\right| = \left|\int_{\cy} F_{\cx}(z_1 - Y) dF_{\cy}(Y) - \int_{\cy} F_{\cx}(z_2 - Y) dF_{\cy}(Y)\right|
        \le L_{\cx}|z_1 - z_2|.
    \end{aligned}
\end{equation*}
Inequality applies $L_{\cx}$-Lipschitz continuity of $F_{\cx}$. Hence, $F_{\cz}$ is $\min\{L_{\cx}, L_{\cy}\}$-Lipschitz continuous.
\Halmos

Then for part (I), we have:
\begin{equation*}
    \begin{aligned}
        \text{(I)} &\le \E_\Gamma \left[ \left| \hat{L}_j(\theta_{t,i} - D_{[t:j-1]}|i_j) - \hat{L}_j(\bar{\theta}_{t,i} - D_{[t:j-1]}|i_j) \right| \right]\\
        & \overset{(a)}{=} \E_{i_{[t:j]}} \left[ \E_{\Gamma|i_{[t:j]}} \left[ \left| \hat{L}_j(\theta_{t,i} - D_{[t:j-1]}|i_j) - \hat{L}_j(\bar{\theta}_{t,i} - D_{[t:j-1]}|i_j) \right| \right] \right]\\
        & = \E_{i_{[t:j]}} \left[ \E_{\Gamma|i_{[t:j]}} \left[ \left| (h_j + b_j) \mathbf{1}(\theta_{t,i} - D_{[t:j-1]} \ge D_j) -  (h_j + b_j) \mathbf{1}(\bar{\theta}_{t,i} - D_{[t:j-1]} \ge D_j)\right| \right]\right]\\
        & = \E_{i_{[t:j]}} \left[ (h_j + b_j) \int_{\min\{\theta_{t,i}, \bar{\theta}_{t,i}\}}^{\max\{\theta_{t,i}, \bar{\theta}_{t,i}\}} dF_{[t:j]}(D | i_{[t:j]}) \right]\\
        & \overset{(b)}{\le} (h_j + b_j) L_D |\theta_{t,i} - \bar{\theta}_{t,i}|.
    \end{aligned}
\end{equation*}
Here $F_{[t:j]}(D | i_{[t:j]})$ denote the cumulative distribution function of $D_{[t:j]}$ conditioned on $i_t, \dots, i_j$. Equation (a) applies the law of total expectation, and inequality (b) follows from Assumption~\ref{assumption: inventory}.\ref{assumption: inventory Lipschitz of rv} and Lemma~\ref{lemma:lip_preserve}, which implies that $F_{[t:j]}(D | i_{[t:j]})$ is $L_D$-Lipschitz. 

The analysis of part (II) uses the following lemma.

\begin{lemma}
\label{lemma: lipschitz s_t}
    Suppose that Assumption~\ref{assumption: inventory} holds. Then for any $c \in \R$, we have
    \begin{equation}
    \label{ineq: lipschitz s_t 1}
        \E_\Gamma \left[ \left| \mathbf{1}(c \ge x_t) - \mathbf{1}(c \ge \bar{x}_t) \right| \right] \le \max \{L_D, L_\rho\} \sum_{l=1}^{t-1} \max_{i\in\ci}|\theta_{l,i} - \bar{\theta}_{l,i}|.
    \end{equation}
    Consequently, the following inequality holds.
    \begin{equation}
    \label{ineq: lipschitz s_t 2}
        \E_\Gamma \left[ \left| \mathbf{1}(\theta_{t, i} \ge x_t) - \mathbf{1} (\bar{\theta}_{t, i} \ge \bar{x}_t) \right| \right] \le t \max \{L_D, L_\rho\} |\theta_{t,i} - \bar{\theta}_{t,i}| + \max \{L_D, L_\rho\} \sum_{l=1}^{t-1} \max_{i\in\ci}|\theta_{l,i} - \bar{\theta}_{l,i}|.
    \end{equation}
\end{lemma}
\proof{Proof of Lemma~\ref{lemma: lipschitz s_t}} We prove the first argument by mathematical induction.

\textbf{Induction Base}: When $t = 1$, $\bar{x}_1 = x_1$ by definition. Then for any $c\in\R$, we have
\begin{equation*}
    \begin{aligned}
        \E_\Gamma \left[ \left| \mathbf{1}(c \ge x_1) - \mathbf{1} (c \ge \bar{x}_1) \right| \right] = 0.
    \end{aligned}
\end{equation*}

\textbf{Induction Step}: Suppose that for any $c\in\R$ and $k < t$, we obtain
\begin{equation*}
    \E_\Gamma \left[ \left| \mathbf{1}(c \ge x_k) - \mathbf{1} (c \ge \bar{x}_k) \right| \right] \le \max\{L_D, L_\rho\} \sum_{l=1}^{k-1} \max_{i\in\ci}|\theta_{l,i} - \bar{\theta}_{l,i}|.
\end{equation*}
Then by definition, for any $c\in\R$, we have
\begin{equation*}
    \begin{aligned}
        & \ \E_\Gamma \left[ \left| \mathbf{1}(c \ge x_{k+1}) - \mathbf{1}(c \ge \bar{x}_{k+1}) \right| \right]\\
        = & \ \E_\Gamma \left[ \left| \mathbf{1}(c \ge x_{k}\vee\theta_{k, i_k} - D_{k}) - \mathbf{1}(c \ge \bar{x}_{k} \vee \bar{\theta}_{k, i_k} - D_{k}) \right| \right]\\
        \overset{(a)}{=} & \ \E_\Gamma \left[ \left| \mathbf{1}(c \ge x_{k} - D_{k}, c \ge \theta_{k, i_k} - D_{k}) - \mathbf{1}(c \ge \bar{x}_{k} - D_{k}, c \ge \bar{\theta}_{k, i_k} - D_{k}) \right| \right].
    \end{aligned}
\end{equation*}
Here equality $(a)$ holds because $\mathbf{1}(t\ge x\vee y)$ is equivalent to $\mathbf{1}(t\ge x, t \ge y)$. Then from the triangle inequality, we obtain the following:
\begin{equation*}
    \begin{aligned}
        & \ \E_\Gamma \left[ \left| \mathbf{1}(c \ge x_{k} - D_{k}, c \ge \theta_{k, i_k} - D_{k}) - \mathbf{1}(c \ge \bar{x}_{k} - D_{k}, c \ge \bar{\theta}_{k, i_k} - D_{k}) \right| \right]\\
        \le & \ \E_\Gamma \left[ \left| \mathbf{1}(c \ge x_{k} - D_{k}) \times (\mathbf{1}(c \ge \theta_{k, i_k} - D_{k}) - \mathbf{1}(c \ge \bar{\theta}_{k, i_k} - D_{k})) \right| \right]\\
        & \ + \E_\Gamma \left[ \left| \mathbf{1}(c \ge \bar{\theta}_{k, i_k} - D_{k}) \times (\mathbf{1}(c \ge x_{k} - D_{k}) - \mathbf{1}(c \ge \bar{x}_{k} - D_{k})) \right| \right]\\
        \le & \ \underbrace{\E_\Gamma \left[ \left| \mathbf{1}(c \ge \theta_{k, i_k} - D_{k}) - \mathbf{1}(c \ge \bar{\theta}_{k, i_k} - D_{k}) \right| \right]}_{\text{(A)}} + \underbrace{\E_\Gamma \left[ \left| \mathbf{1}(c \ge x_{k} - D_{k}) - \mathbf{1}(c \ge \bar{x}_{k} - D_{k}) \right| \right]}_{\text{(B)}}.
    \end{aligned}
\end{equation*}
For the first part, we have
\begin{equation*}
    \begin{aligned}
        \text{(A)} &\overset{(a)}{=} \E_{i_k} \left[ \E_{\Gamma|i_k} \left[ \left| \mathbf{1}(c \ge \theta_{k, i_k} - D_{k}) - \mathbf{1}(c \ge \bar{\theta}_{k, i_k} - D_{k}) \right| \right] \right]\\
        &= \E_{i_k} \left[\int_{\min\{\theta_k, \bar{\theta}_k\} - c}^{\max\{\theta_k, \bar{\theta}_k\} - c}dP_D(D_k|i_k) \right]\\
        &\overset{(b)}{\le} L_D \E_{i_k} \left[ | \theta_{k,i_k} - \bar{\theta}_{k,i_k} | \right]\\
        &\le L_D \max_{i\in\ci} | \theta_{k,i} - \bar{\theta}_{k,i} |.
    \end{aligned}
\end{equation*}
where equation (a) applies the law of total expectation, and inequality (b) follows from the $L_D$-Lipschitz continuity of $P_D(D_k|i_k)$. For the second part, we have
\begin{equation*}
    \begin{aligned}
        \text{(B)} & \overset{(a)}{=} \E_{D_k} \left[ \E \left[ \left| \mathbf{1}(c \ge x_{k} - D_{k}) - \mathbf{1}(c \ge \bar{x}_{k} - D_{k}) \right| \Bigl| D_k \right] \right]\\
        & \overset{(b)}{\le} \E_{D_k} \left[ \max \{L_D, L_\rho\} \sum_{l=1}^{k-1} \max_{i\in\ci}|\theta_{l,i} - \bar{\theta}_{l,i}| \right]\\
        & = \max \{L_D, L_\rho\} \sum_{l=1}^{k-1} \max_{i\in\ci}|\theta_{l,i} - \bar{\theta}_{l,i}|.
    \end{aligned}
\end{equation*}
Here equality $(a)$ applies the law of total expectation, and inequality $(b)$ comes from the induction assumption. Therefore, we can finish the induction by combining the results of (A) and (B),
\begin{equation*}
    \E_\Gamma \left[ \left| \mathbf{1}(c \ge x_{k+1}) - \mathbf{1}(c \ge \bar{x}_{k+1}) \right| \right] \le \text{(A)} + \text{(B)} \le \max \{L_D, L_\rho\} \sum_{l=1}^{k} \max_{i\in\ci}|\theta_{l,i} - \bar{\theta}_{l,i}|.
\end{equation*}

To show the second argument, from the triangle inequality, we have
\begin{equation*}
    \begin{aligned}
        \E_\Gamma \left[ \left| \mathbf{1}(\theta_{t, i} \ge x_t) - \mathbf{1} (\bar{\theta}_{t, i} \ge \bar{x}_t) \right| \right] \le \E_\Gamma \left[ \left| \mathbf{1}(\theta_{t, i} \ge x_t) - \mathbf{1} (\bar{\theta}_{t, i} \ge x_t) \right| \right] + \E_\Gamma \left[ \left| \mathbf{1}(\bar{\theta}_{t, i} \ge x_t) - \mathbf{1} (\bar{\theta}_{t, i} \ge \bar{x}_t) \right| \right].
    \end{aligned}
\end{equation*}
For the first part, we obtain the following:
\begin{equation*}
    \begin{aligned}
        & \ \E_\Gamma \left[ \left| \mathbf{1}(\theta_{t, i} \ge x_t) - \mathbf{1} (\bar{\theta}_{t, i} \ge x_t) \right| \right]\\
        = & \ \E_\Gamma \left[ \left| \mathbf{1}(\theta_{t, i} \ge x_{t-1}\vee\theta_{t-1, i_{t-1}} - D_{t-1}) - \mathbf{1}(\bar{\theta}_{t, i} \ge x_{t-1} \vee \theta_{t-1, i_{t-1}} - D_{t-1}) \right| \right]\\
        = & \ \E_\Gamma \left[ \left| \mathbf{1}(\theta_{t, i} \ge x_{t-1} - D_{t-1}, \theta_{t, i} \ge \theta_{t-1, i_{t-1}} - D_{t-1}) - \mathbf{1}(\bar{\theta}_{t, i} \ge x_{t-1} - D_{t-1}, \bar{\theta}_{t, i} \ge \theta_{t-1, i_{t-1}} - D_{t-1}) \right| \right]\\
        \le & \ \E_\Gamma \left[ \left| \mathbf{1}(\theta_{t, i} \ge x_{t-1} - D_{t-1}) - \mathbf{1}(\bar{\theta}_{t, i} \ge x_{t-1} - D_{t-1}) \right| \right]\\
        &\quad\quad\quad\quad\quad\quad+ \E_\Gamma \left[ \left| \mathbf{1}(\theta_{t, i} \ge \theta_{t-1, i_{t-1}} - D_{t-1}) - \mathbf{1}(\bar{\theta}_{t, i} \ge \theta_{t-1, i_{t-1}} - D_{t-1}) \right| \right]\\
    \end{aligned}
\end{equation*}
Applying mathematical induction in a similar way, we prove that 
\begin{equation*}
    \E_\Gamma \left[ \left| \mathbf{1}(\theta_{t, i} \ge x_t) - \mathbf{1} (\bar{\theta}_{t, i} \ge x_t) \right| \right] \le t \max \{L_D, L_\rho\} |\theta_{t,i} - \bar{\theta}_{t,i}|.
\end{equation*}
In addition, we have
\begin{equation*}
    \begin{aligned}
        & \ \E_\Gamma \left[ \left| \mathbf{1}(\theta_{t, i} \ge x_t) - \mathbf{1} (\bar{\theta}_{t, i} \ge \bar{x}_t) \right| \right]\\
        \le & \  \E_\Gamma \left[ \left| \mathbf{1}(\theta_{t, i} \ge x_t) - \mathbf{1} (\bar{\theta}_{t, i} \ge x_t) \right| \right] + \E_\Gamma \left[ \left| \mathbf{1}(\bar{\theta}_{t, i} \ge x_t) - \mathbf{1} (\bar{\theta}_{t, i} \ge \bar{x}_t) \right| \right]\\
        \le & \ t \max \{L_D, L_\rho\} |\theta_{t,i} - \bar{\theta}_{t,i}| + \max \{L_D, L_\rho\} \sum_{l=1}^{t-1} \max_{i\in\ci}|\theta_{l,i} - \bar{\theta}_{l,i}|.
    \end{aligned}
\end{equation*}
This concludes the proof. \Halmos

For part (II), applying $\left| \hat{L}_j(x_j | i_j) \right| \le \max\{h_j, b_j\}$ and Lemma~\ref{lemma: lipschitz s_t} returns
\begin{equation*}
    \text{(II)} \le \max\{h_j, b_j\} \max \{L_D, L_\rho\}  \left(t |\theta_{t,i} - \bar{\theta}_{t,i}| + \sum_{l=1}^{t-1} \max_{i\in\ci}|\theta_{l,i} - \bar{\theta}_{l,i}| \right).
\end{equation*}

The analysis of part (III) uses the following lemma. 
\begin{lemma}
\label{lemma: lipschitz lambda}
    Suppose that Assumption~\ref{assumption: inventory} holds. Then for any $\theta, \bar{\theta} \in \Theta$ and $j\ge t$, we have
    \begin{equation*}
        \E_\Gamma \left[ \left| \mathbf{1} \left( j \le \Lambda_t(\theta, \Gamma) \right) - \mathbf{1} \left( j \le \Lambda_t(\bar{\theta}, \Gamma) \right) \right| \right] \le L_D \left( (j-t)|\theta_{t,i_t} - \bar{\theta}_{t,i_t}| + \sum_{k=t+1}^{j} \max_{i\in\ci}|\theta_{k, i} - \bar{\theta}_{k,i}| \right).
    \end{equation*}
\end{lemma}
\proof{Proof of Lemma~\ref{lemma: lipschitz lambda}} By definition, the event $j \le \Lambda_t(\theta, \Gamma)$ is equivalent to the event $\theta_{t,i_t} - D_{[t:l-1]} \ge \theta_{l, i_l}, \forall t\le l \le j$. Therefore we have
\begin{equation*}
    \mathbf{1} \left( j \le \Lambda_t(\theta, \Gamma) \right) = \mathbf{1}(D_t \le \theta_{t, i_t} - \theta_{t+1, i_{t+1}}, \dots, D_{[t:j-1]} \le \theta_{t,i_t} - \theta_{j,i_j}).
\end{equation*}
From the triangle inequality and telescoping sum,
\begin{equation*}
    \begin{aligned}
        & \ \E_\Gamma \left[ \left| \mathbf{1} \left( j \le \Lambda_t(\theta, \Gamma) \right) - \mathbf{1} \left( j \le \Lambda_t(\bar{\theta}, \Gamma) \right) \right| \right]\\
        = & \ \E_\Gamma \biggl[ \Bigl| \mathbf{1}(D_t \le \theta_{t,i_t} - \theta_{t+1, i_{t+1}}, \dots, D_{[t:j-1]} \le \theta_{t,i_t} - \theta_{j, i_j})\\
        &\quad- \mathbf{1}(D_t \le \bar{\theta}_{t,i_t} - \bar{\theta}_{t+1, i_{t+1}}, \dots, D_{[t:j-1]} \le \bar{\theta}_{t,i_t} - \bar{\theta}_{j, i_j}) \Bigr| \biggr]\\
        \le & \ \sum_{k=t+1}^{j} \E_\Gamma \biggl[ \Bigl| \mathbf{1}(D_{[t:t'-1]} \le \theta_{t,i_t} - \theta_{t', i_{t'}}, \forall t + 1\le t'\le k, D_{[t:t'-1]} \le \bar{\theta}_{t,i_t} - \bar{\theta}_{t', i_{t'}}, \forall k < t'\le j)\\
        & \quad \quad \ \ - \mathbf{1}(D_{[t:t'-1]} \le \theta_{t,i_t} - \theta_{t', i_{t'}}, \forall t+1 \le t'< k, D_{[t:t'-1]} \le \bar{\theta}_{t,i_t} - \bar{\theta}_{t', i_{t'}}, \forall k \le t'\le j) \Bigr| \biggr]\\
        \le & \ \sum_{k=t+1}^{j} \E_\Gamma \left[ \left| \mathbf{1} (D_{[t:k-1]} \le \theta_{t, i_t} - \theta_{k,i_k}) - \mathbf{1} (D_{[t:k-1]} \le \bar{\theta}_{t,i_t} - \bar{\theta}_{k, i_k}) \right| \right]\\
        \overset{(a)}{\le} & \ L_D \sum_{k=t+1}^{j} \bigl( |\theta_{t,i_t} - \bar{\theta}_{t,i_t}| + \max_{i\in\ci}|\theta_{k, i} - \bar{\theta}_{k,i}| \bigr)\\
        = & \ L_D \Bigl( (j-t)|\theta_{t,i_t} - \bar{\theta}_{t,i_t}| + \sum_{k=t+1}^{j} \max_{i\in\ci}|\theta_{k, i} - \bar{\theta}_{k,i}| \Bigr).
    \end{aligned}
\end{equation*}
Here inequality (a) applies Lemma~\ref{lemma:lip_preserve} that the cumulative distribution function of random demand $D_{[t:k-1]}$ is $L_D$-Lipschitz continuous. This completes the proof. \Halmos

Applying Lemma~\ref{lemma: lipschitz lambda}, we have
\begin{equation*}
    \text{(III)} \le \max\{h_j, b_j\} L_D \Bigl( (j-t)|\theta_{t,i} - \bar{\theta}_{t,i}| + \sum_{k=t+1}^{j} \max_{i\in\ci}|\theta_{k, i} - \bar{\theta}_{k,i}| \Bigr).
\end{equation*}
Combining the results for part (I), (II), and (III), we bound the partial derivative:
\begin{equation*}
    \begin{aligned}
        \left| \frac{\partial}{\partial \theta_{t,i}}l(\theta) - \frac{\partial}{\partial \theta_{t,i}}l(\bar{\theta}) \right| &\le \max_{t\in[T]}\{h_t + b_t\} \max \{L_D, L_\rho\}\sum_{j=t}^T \Bigl(j |\theta_{t,i} - \bar{\theta}_{t,i}| + \sum_{l=1}^{j} \max_{i\in\ci}|\theta_{l,i} - \bar{\theta}_{l,i}| \Bigr)\\
        &\le \max_{t\in[T]}\{h_t + b_t\} \max \{L_D, L_\rho\} \Bigl(T^2 \max_{i\in\ci}|\theta_{t,i} - \bar{\theta}_{t,i}| + T\sum_{l=1}^T \max_{i\in\ci}|\theta_{l,i} - \bar{\theta}_{l,i}| \Bigr).
    \end{aligned}
\end{equation*}
By definition, we obtain the following:
\begin{equation*}
    \begin{aligned}
        \|\nabla l(\theta) - \nabla l(\bar{\theta})\|_2^2 &= \sum_{t\in[T], i\in\ci} \left| \frac{\partial}{\partial \theta_{t,i}}l(\theta) - \frac{\partial}{\partial \theta_{t,i}}l(\bar{\theta}) \right|^2\\
        \le &\ \max_{t\in[T]}\{h_t + b_t\}^2\max \{L_D, L_\rho\}^2\sum_{t\in[T], i\in\ci} \Bigl( 2T^4 \max_{i\in\ci}|\theta_{t,i} - \bar{\theta}_{t,i}|^2 + 2T^2 \Bigl( \sum_{l\in[T]} \max_{i\in\ci}|\theta_{l,i} - \bar{\theta}_{l,i}| \Bigr)^2 \Bigr)\\
        \le &\ \max_{t\in[T]}\{h_t + b_t\}^2\max \{L_D, L_\rho\}^2\sum_{t\in[T], i\in\ci} \Bigl( 2T^4 \sum_{i\in\ci}|\theta_{t,i} - \bar{\theta}_{t,i}|^2 + 2T^3 \sum_{l\in[T]} \max_{i\in\ci}|\theta_{l,i} - \bar{\theta}_{l,i}|^2 \Bigr)\\
        \le &\ \max_{t\in[T]}\{h_t + b_t\}^2\max \{L_D, L_\rho\}^2\sum_{t\in[T], i\in\ci} \Bigl( 2T^4 \sum_{i\in\ci}|\theta_{t,i} - \bar{\theta}_{t,i}|^2 + 2T^3 \sum_{l\in[T]} \sum_{i\in\ci}|\theta_{l,i} - \bar{\theta}_{l,i}|^2 \Bigr)\\
        = &\ \max_{t\in[T]}\{h_t + b_t\}^2\max \{L_D, L_\rho\}^2 4T^4|\ci| \|\theta - \bar{\theta}\|_2^2.
    \end{aligned}
\end{equation*}
Therefore, $l(\theta)$ is $(2\max_{t\in[T]}\{h_t + b_t\}\max \{L_D, L_\rho\}T^2\sqrt{|\ci|})$-smooth, implying that Assumption~\ref{assumption: smoothness} holds. Given that the smoothness parameter and the P{\L}K constant all admit a polynomial dependence on the planning horizon, an $\epsilon$-optimal policy can be obtained using sample size in $\tilde{\co}(\epsilon^{-1})$ and polynomial in terms of the planning horizon by stochastic policy gradient methods.

\section{Omitted Proofs in Section~\ref{section: cash balance}} \label{appendix: cash balance}
The proof for the stochastic cash balance problem in Section~\ref{section: cash balance} shares some similarities with the inventory system in Section~\ref{section: inventory system} yet the per-period decision is a two-dimensional vector. We demonstrate the full proof for completeness.

\subsection{P{\L}K Condition of Optimal Q-value Function}
\proof{Proof of Lemma~\ref{lemma: cash balance single stage KL}} We divide the proof into three parts. First, we demonstrate the relationship between three suboptimality gaps $F_t(\theta_t) - F_t(\theta_t^*)$, $\underline{f}_t(\underline{\theta}_t) - \underline{f}_t(\underline{\theta}_t^*)$, and $\bar{f}_t(\bar{\theta}_t) - \bar{f}_t(\bar{\theta}_t^*)$. Next, we show how their gradients relate to each other. Finally, we prove the P{\L}K property of $F_t$.

\textbf{Step 1: Relationship between suboptimality gaps.} In the feasible region $\Theta_t$, we have $\underline{\theta}_t \le \bar{\theta}_t$. For function $f_t$, the following equations hold:
\begin{equation}
\label{cash balance decomposition}
    f_t \left( (s_t \vee \underline{\theta}_t) \wedge \bar{\theta}_t \right) = f_t (s_t \vee \underline{\theta}_t ) + f_t (s_t \wedge \bar{\theta}_t ) - f_t(s_t).
\end{equation}
It further holds that
\begin{equation*}
    \begin{aligned}
        F_t(\theta_t) - F_t(\theta_t^*) = & \ \E_{s_t\sim\rho_t(\cdot|\pi_\theta)} \left[ Q_t^* \left(s_t, \pi_t(s_t|\theta_t) \right) - Q_t^* \left(s_t, \pi_t(s_t|\theta_t^*) \right) \right]\\
        = & \ \E_{s_t\sim\rho_t(\cdot|\pi_\theta)} \left[ c\left( (s_t\vee\underline{\theta}_t) \wedge \bar{\theta}_t, s_t \right) + f_t\left( (s_t\vee\underline{\theta}_t) \wedge \bar{\theta}_t \right) - c\left( (s_t\vee\underline{\theta}_t^*) \wedge \bar{\theta}_t^*, s_t \right) - f_t\left( (s_t\vee\underline{\theta}_t^*) \wedge \bar{\theta}_t^* \right) \right]\\
        = & \ \underbrace{\E_{s_t\sim\rho_t(\cdot|\pi_\theta)} \left[ c(s_t\vee\underline{\theta}_t , s_t ) + f_t (s_t\vee\underline{\theta}_t ) - c(s_t\vee\underline{\theta}_t^*, s_t ) - f_t(s_t\vee\underline{\theta}_t^*) \right]}_{(\text{I})}\\
        & \quad\quad\quad\quad + \underbrace{\E_{s_t\sim\rho_t(\cdot|\pi_\theta)} \left[ c(s_t\wedge \bar{\theta}_t, s_t) + f_t(s_t\wedge\bar{\theta}_t) - c(s_t\wedge \bar{\theta}_t^*, s_t) - f_t(s_t\wedge\bar{\theta}_t^*) \right]}_{(\text{II})},
    \end{aligned}
\end{equation*}
where the last equation comes from (\ref{cash balance decomposition}). We analyze the first term (I). Without loss of generality, we assume that $\underline{\theta}_t \le \underline{\theta}_t^*$. With the expression of $c$, it holds that
\begin{equation*}
    \begin{aligned}
        (\text{I}) &= \int_{-\infty}^{\underline{\theta}_t} \left( \underline{f}_t(\underline{\theta}_t) - ks_t \right) d \rho_t(s_t|\pi_\theta) + \int_{\underline{\theta}_t}^{+\infty} f_t(s_t) d \rho_t(s_t|\pi_\theta) \\
        & \quad\quad - \int_{-\infty}^{\underline{\theta}_t^*} \left( \underline{f}_t(\underline{\theta}_t^*) - ks_t \right) d \rho_t(s_t|\pi_\theta) - \int_{\underline{\theta}_t^*}^{+\infty} f_t(s_t) d \rho_t(s_t|\pi_\theta) \\
        & = \int_{-\infty}^{\underline{\theta}_t} \left( \underline{f}_t(\underline{\theta}_t) - \underline{f}_t(\underline{\theta}_t^*) \right) d \rho_t(s_t|\pi_\theta) + \int_{\underline{\theta}_t}^{\underline{\theta}_t^*} \left( f_t(s_t) + ks_t - \underline{f}_t(\underline{\theta}_t^*) \right) d \rho_t(s_t|\pi_\theta).
    \end{aligned}
\end{equation*}
For the right-hand-side, we have the following inequalities:
\begin{equation*}
    \begin{aligned}
        \int_{\underline{\theta}_t}^{\underline{\theta}_t^*} \left( f_t(s_t) + ks_t - \underline{f}_t(\underline{\theta}_t^*) \right) d \rho_t(s_t|\pi_\theta) &\overset{(a)}{=} \int_{\underline{\theta}_t}^{\underline{\theta}_t^*} \left( \underline{f}_t(s_t) - \underline{f}_t(\underline{\theta}_t^*) \right) d \rho_t(s_t|\pi_\theta)\\
        &\overset{(b)}{\le} \int_{\underline{\theta}_t}^{\underline{\theta}_t^*} \left( \underline{f}_t(\underline{\theta}_t) - \underline{f}_t(\underline{\theta}_t^*) \right) d \rho_t(s_t|\pi_\theta),
    \end{aligned}
\end{equation*}
where Equation (a) uses the definition of $\underline{f}_t$ and $f_t$, and inequality (b) holds because $\underline{f}_t$ is non-increasing on the interval $[\underline{\theta}_t, \underline{\theta}_t^*]$. Therefore, we conclude that $(\text{I}) \le \underline{f}_t(\underline{\theta}_t) - \underline{f}_t(\underline{\theta}_t^*)$. The same result holds when $\underline{\theta}_t > \underline{\theta}_t^*$.

For the second term (II), we apply the same technique. Without loss of generality, we assume that $\bar{\theta}_t \le \bar{\theta}_t^*$:
\begin{equation*}
    \begin{aligned}
        (\text{II}) &= \int_{-\infty}^{\bar{\theta}_t} f_t(s_t) d \rho_t(s_t|\pi_\theta) + \int_{\bar{\theta}_t}^{+\infty} \left( \bar{f}_t(\bar{\theta}_t) + qs_t \right) d \rho_t(s_t|\pi_\theta) \\
        & \quad\quad - \int_{-\infty}^{\bar{\theta}_t^*} f_t(s_t) d \rho_t(s_t|\pi_\theta) - \int_{\bar{\theta}_t^*}^{+\infty} \left( \bar{f}_t(\bar{\theta}_t^*) + qs_t \right) d \rho_t(s_t|\pi_\theta) \\
        & = \int_{\bar{\theta}_t}^{\bar{\theta}_t^*} \left( \bar{f}_t(\bar{\theta}_t) + qs_t - f_t(s_t) \right) d \rho_t(s_t|\pi_\theta) + \int_{\bar{\theta}_t^*}^{+\infty} \left( \bar{f}_t(\bar{\theta}_t) - \bar{f}_t(\bar{\theta}_t^*) \right) d \rho_t(s_t|\pi_\theta).
    \end{aligned}
\end{equation*}
Similarly, it holds that
\begin{equation*}
    \begin{aligned}
        \int_{\bar{\theta}_t}^{\bar{\theta}_t^*} \left( \bar{f}_t(\bar{\theta}_t) + qs_t - f_t(s_t) \right) d \rho_t(s_t|\pi_\theta) &\overset{(a)}{=} \int_{\bar{\theta}_t}^{\bar{\theta}_t^*} \left( \bar{f}_t(\bar{\theta}_t) - \bar{f}_t(s_t) \right) d \rho_t(s_t|\pi_\theta)\\
        &\overset{(b)}{\le} \int_{\bar{\theta}_t}^{\bar{\theta}_t^*} \left( \bar{f}_t(\bar{\theta}_t) - \bar{f}_t(\bar{\theta}_t^*) \right) d \rho_t(s_t|\pi_\theta),
    \end{aligned}
\end{equation*}
where Equation (a) uses the definition of $\bar{f}_t$ and $f_t$, and Inequality (b) holds because $\bar{f}_t$ is non-increasing on the interval $[\bar{\theta}_t, \bar{\theta}_t^*]$. We thus have $(\text{II}) \le \bar{f}_t(\bar{\theta}_t) - \bar{f}_t(\bar{\theta}_t^*)$. The same result holds when $\bar{\theta}_t > \bar{\theta}_t^*$. Combining all the results, we conclude that
\begin{equation*}
    F_t(\theta_t) - F_t(\theta_t^*) \le (\text{I}) + (\text{II}) \le \underline{f}_t(\underline{\theta}_t) - \underline{f}_t(\underline{\theta}_t^*) + \bar{f}_t(\bar{\theta}_t) - \bar{f}_t(\bar{\theta}_t^*)
\end{equation*}

\textbf{Step 2: Relationship between gradients.} 
By definition, we calculate the gradient of $F_t$. We first show the partial derivative for $\underline{\theta}$ using the definition of $F_t$, $f_t$, $\underline{f}_t$, and $\bar{f}_t$.
\begin{equation*}
    \begin{aligned}
        \nabla_{\underline{\theta}_t} F_t(\theta_t) &= \nabla_{\underline{\theta}_t} \E_{s_t\sim\rho_t(\cdot|\pi_\theta)} \left[ Q_t^{\pi_{\theta^*}} \left( s_t, \pi_t(s_t | \theta_t) \right) \right]\\
        &= \nabla_{\underline{\theta}_t} \left[ \int_{-\infty}^{\underline{\theta}_t} \left( \underline{f}_t(\underline{\theta}_t) - ks_t \right) d \rho_t(s_t|\pi_\theta) + \int_{\underline{\theta}_t}^{\bar{\theta}_t} f_t(s_t) d \rho_t(s_t|\pi_\theta) + \int_{\bar{\theta}_t}^{\infty} \left( \bar{f}_t(\bar{\theta}_t) + qs_t \right) d\rho_t(s_t|\pi_\theta) \right]\\
        &\overset{(a)}{=} \underline{f}_t(\underline{\theta}_t) - k\underline{\theta}_t + \int_{-\infty}^{\underline{\theta}_t} \underline{f}_t'(\underline{\theta}_t) d\rho_t(s_t|\pi_\theta) - f_t(\underline{\theta}_t)\\[4pt]
        &= \P(s_t \le \underline{\theta}_t) \underline{f}_t'(\underline{\theta}_t).
    \end{aligned}
\end{equation*}
Equation (a) uses the Leibniz rule. Similarly, we derive $\nabla_{\bar{\theta}_t} F_t(\theta_t) = \P(s_t \ge \bar{\theta}_t) \bar{f}_t'(\bar{\theta}_t)$.

\textbf{Step 3: P{\L}K Condition of $F_t$}. By Assumption~\ref{assumption: cash balance}.\ref{assumption: cash balance SC}, the per-period holding or backlogging cost is $\min_{t\in[T]}\{h_t+b_t\}\mu_D$-strongly convex over $[\underline{B}, \bar{B}]$. By the convexity of cost-to-go functions, we have that $\underline{f}_t$ and $\bar{f}_t$ are both $\min_{t\in[T]}\{h_t+b_t\}\mu_D$-strongly convex over $[\underline{B}, \bar{B}]$. Therefore, we can derive
\begin{equation*}
    \begin{aligned}
        & \ F_t(\theta_t) - F_t(\theta_t^*)\\
        & \le \underline{f}_t(\underline{\theta}_t) - \underline{f}_t(\underline{\theta}_t^*) + \bar{f}_t(\bar{\theta}_t) - \bar{f}_t(\bar{\theta}_t^*)\\
        &\overset{(a)}{\le} \underline{f}_t'(\underline{\theta}_t)(\underline{\theta}_t - \underline{\theta}_t^*) - \frac{\min_{t\in[T]}\{h_t+b_t\}\mu_D}{2}\|\underline{\theta}_t - \underline{\theta}_t^*\|_2^2 + \bar{f}_t'(\bar{\theta}_t)(\bar{\theta}_t - \bar{\theta}_t^*) - \frac{\min_{t\in[T]}\{h_t+b_t\}\mu_D}{2}\|\bar{\theta}_t - \bar{\theta}_t^*\|_2^2\\
        &=\frac{\nabla_{\underline{\theta}_t}F_t(\theta_t)}{\P(s_t \le \underline{\theta}_t)} (\underline{\theta}_t - \underline{\theta}_t^*) - \frac{\min_{t\in[T]}\{h_t+b_t\}\mu_D}{2}\|\underline{\theta}_t - \underline{\theta}_t^*\|_2^2 + \frac{\nabla_{\bar{\theta}_t}F_t(\theta_t)}{\P(s_t \ge \bar{\theta}_t)} (\bar{\theta}_t - \bar{\theta}_t^*) - \frac{\min_{t\in[T]}\{h_t+b_t\}\mu_D}{2}\|\bar{\theta}_t - \bar{\theta}_t^*\|_2^2\\
        &\overset{(b)}{\le} \alpha^{-1} \left\langle \nabla_{\theta_t} F_t(\theta_t), \begin{bmatrix} \underline{\theta}_t \\[1pt] \bar{\theta}_t \end{bmatrix} - \begin{bmatrix} \underline{\theta}_t^* \\[1pt] \bar{\theta}_t^* \end{bmatrix} \right\rangle - \frac{\min_{t\in[T]}\{h_t+b_t\}\mu_D}{2} \left\| \begin{bmatrix} \underline{\theta}_t \\[1pt] \bar{\theta}_t \end{bmatrix} - \begin{bmatrix} \underline{\theta}_t^* \\[1pt] \bar{\theta}_t^* \end{bmatrix} \right\|_2^2\\
        &\overset{(c)}{\le} \max_{\theta'_t \in \Theta_t} \left\{ \alpha^{-1} \langle \nabla_{\theta_t} F_t(\theta_t), \theta_t - \theta'_t \rangle - \frac{\min_{t\in[T]}\{h_t+b_t\}\mu_D}{2} \| \theta_t - \theta'_t \|_2^2 \right\}.
    \end{aligned}
\end{equation*}
Inequality (a) uses the strong convexity of $\underline{f}_t$ and $\bar{f}_t$. The equality holds by the explicit expression derived in Step 2. Inequality (b) holds because $\P(s_t \le \underline{\theta}_t) \ge \alpha > 0$ and $\P(s_t \ge \bar{\theta}_t) \ge \alpha > 0$. Inequality (c) utilizes the fact that $\theta^* \in \Theta_t$. Therefore, $F_t(\theta_t)$ satisfies the $(\alpha^{-1}, \min_{t\in[T]}\{h_t+b_t\}\mu_D)$ gradient dominance condition, and thus the P{\L}K condition with constant $\min_{t\in[T]}\{h_t+b_t\}\mu_D\alpha^2$ by Lemma~\ref{lemma: gradient dominance implies KL}. This completes the proof. \Halmos

\subsection{Gradient Formulation}

\proof{Proof of Proposition~\ref{proposition: cash balance gradient formulation}}

By the Bellman equation (\ref{bellman equation}), we derive the recursive form of $(V_t^{\pi_\theta})'(s_t)$ for any $t \in [T]$:
\begin{equation*}
    \begin{aligned}
        (V_t^{\pi_\theta})' (s_t) &= \frac{\partial}{\partial s_t} Q_t^{\pi_\theta} \left( s_t, \pi_t(s_t | \theta_t) \right)\\
        &= \frac{\partial}{\partial s_t} \left( c \left( (s_t \vee \underline{\theta}_t) \wedge \bar{\theta}_t, s_t \right) + L_t \left( (s_t \vee \underline{\theta}_t) \wedge \bar{\theta}_t \right) + \E_{D_t} \left[ (V^{\pi_\theta}_{t+1})' \left( (s_t \vee \underline{\theta}_t) \wedge \bar{\theta}_t - D_t \right) \right] \right)\\
        &= -k \mathbf{1}(s_t \le \underline{\theta}_t) + q \mathbf{1}(s_t \ge \bar{\theta}_t) + \mathbf{1}(\underline{\theta}_t < s_t < \bar{\theta}_t) \times \left( L'_t(s_t) + \E_{D_t}\left[ (V^{\pi_\theta}_{t+1})'(s_t - D_t) \right] \right)
    \end{aligned}
\end{equation*}
with $(V_{T+1}^{\pi_\theta})'(\cdot) = 0$. For the policy gradient objective function $l(\theta)$, we calculate the partial derivative
\begin{equation*}
    \begin{aligned}
        \frac{\partial}{\partial \underline{\theta}_t} l(\theta) &\overset{(a)}{=} \E_{s_t \sim \rho_t(\cdot|\pi_\theta)} \left[ \frac{\partial}{\partial \underline{\theta}_t} Q^{\pi_\theta}_t \left( s_t, \pi_t(s_t | \theta_t) \right) \right]\\
        &\overset{(b)}{=} \E_{s_t \sim \rho_t(\cdot|\pi_\theta)} \left[ \frac{\partial}{\partial \underline{\theta}_t} \pi_t(s_t|\theta_t) \times  \frac{\partial}{\partial a_t} Q^{\pi_\theta}_t (s_t, a_t ) \biggl|_{a_t = \pi_t(s_t|\theta_t) } \right]\\
        &\overset{(c)}{=} \E_{s_t \sim \rho_t(\cdot|\pi_\theta)} \left[ \mathbf{1}(\underline{\theta}_t \ge s_t) \times \frac{\partial}{\partial a_t} \left( C_t(s_t, a_t) + \E_{D_t} \left[ V_{t+1}^{\pi_\theta} (s_t + a_t - D_t) \right] \right) \biggl|_{a_t = \pi_t(s_t|\theta_t)} \right]\\
        &\overset{(d)}{=} \E_{s_t \sim \rho_t(\cdot|\pi_\theta)} \left[ \mathbf{1}(\underline{\theta}_t \ge s_t) \times \left( k + L'_t(\underline{\theta}_t) + \E_{D_t} \left[ (V_{t+1}^{\pi_\theta})' (\underline{\theta}_t - D_t) \right] \right) \right].
    \end{aligned}
\end{equation*}
Equation (a) utilizes the Deterministic Policy Gradient Theorem \citep{silver2014deterministic}. Equation (b) applies the chain rule. Equation (c) uses the Bellman equation (\ref{bellman equation}). Lastly, equation (d) holds because $\mathbf{1}(\underline{\theta}_t \ge s_t)$. Similarly, we can calculate the partial derivative
\begin{equation*}
    \frac{\partial}{\partial \bar{\theta}_t} l(\theta) = \E_{s_t \sim \rho_t(\cdot|\pi_\theta)} \left[ \mathbf{1}( \bar{\theta}_t \le s_t) \times \left( -q + L'_t(\bar{\theta}_t) + \E_{D_t} \left[ (V_{t+1}^{\pi_\theta})' (\bar{\theta}_t - D_t) \right] \right) \right].
\end{equation*}
This concludes the proof. \Halmos

\subsection{Bounded Gradient}

\proof{Proof of Lemma~\ref{lemma: cash balance bounded gradient}}
Following Proposition~\ref{proposition: cash balance gradient formulation}, we can bound the partial derivative as follows:
\begin{equation*}
    \begin{aligned}
        \left| \frac{\partial}{\partial \underline{\theta}_t} l(\theta) \right| &= \left| \E_{s_t \sim \rho_t(\cdot|\pi_\theta)} \left[ \mathbf{1}(\underline{\theta}_t \ge s_t) \times \left( k + L'_t(\underline{\theta}_t) + \E_{D_t} \left[ (V_{t+1}^{\pi_\theta})' (\underline{\theta}_t - D_t) \right] \right) \right] \right|\\
        &\overset{(a)}{\le} k + |L'_t(\underline{\theta}_t)| + \left| \E_{D_t} \left[ (V_{t+1}^{\pi_\theta})' (\underline{\theta}_t - D_t) \right] \right|\\
        &\overset{(b)}{\le} k + \max_{t\in[T]} \left\{ \max\{h_t, b_t\} \right\} + \left| \E_{D_t} \left[ (V_{t+1}^{\pi_\theta})' (\underline{\theta}_t - D_t) \right] \right|.
    \end{aligned}
\end{equation*}
Inequality (a) employs the triangle inequality and utilizes the fact that $\mathbf{1}(\underline{\theta}_t \ge s_t) \le 1$. Inequality (b) holds because $|L_t'(\theta_t)| \le \max_{t\in[T]} \left\{ \max\{h_t, b_t\} \right\}$ for any $\theta_t$ and $t\in[T]$. From (\ref{cash balance: recursive gradient}), we have
\begin{equation}
\label{cash balance: recursive gradient bound}
    \begin{aligned}
        \left| (V_t^{\pi_\theta})' (s_t) \right| &= \left| -k \mathbf{1}(s_t \le \underline{\theta}_t) + q \mathbf{1}(s_t \ge \bar{\theta}_t) + \mathbf{1}(\underline{\theta}_t < s_t < \bar{\theta}_t) \times \left( L'_t(s_t) + \E_{D_t}\left[ (V^{\pi_\theta}_{t+1})'(s_t - D_t) \right] \right) \right| \\
        &\le k + |q| + |L'_t(s_t)| + \left| \E_{D_t} \left[ (V_{t+1}^{\pi_\theta})' (s_t - D_t) \right] \right|\\
        &\le k + |q| + \max_{t\in[T]} \left\{ \max\{h_t, b_t\} \right\} + \left| \E_{D_t} \left[ (V_{t+1}^{\pi_\theta})' (s_t - D_t) \right] \right|.
    \end{aligned}
\end{equation}

Applying (\ref{cash balance: recursive gradient bound}) recursively, we derive
\begin{equation*}
    \begin{aligned}
        \left| \frac{\partial}{\partial \underline{\theta}_t} l(\theta) \right| \le \left( k + |q| + \max_{t\in[T]} \left\{ \max\{h_t, b_t\} \right\} \right) T.
    \end{aligned}
\end{equation*}

Similarly, we have
\begin{equation*}
    \left| \frac{\partial}{\partial \bar{\theta}_t} l(\theta) \right| \le \left( k + |q| + \max_{t\in[T]} \left\{ \max\{h_t, b_t\} \right\} \right) T.
\end{equation*}

Thus, we obtain
\begin{equation*}
    \|\nabla_{\theta_t} l(\theta)\|_2 \le \| \nabla_{\theta_t} l(\theta) \|_1 \le 2\left( k + |q| + \max_{t\in[T]} \left\{ \max\{h_t, b_t\} \right\} \right) T.
\end{equation*}

This completes the proof. \Halmos

\subsection{Sequential Decomposition Inequality}

\proof{Proof of Lemma~\ref{lemma: cash balance bounded sequence gradient difference}} For simplicity, we define $\theta_\alpha = (\theta_{[1:k]}, \theta^*_{[k+1:T]})$ and $\theta_\beta = (\theta_{[1:k-1]}, \theta^*_{[k:T]})$. Furthermore, we denote $\pi_\alpha$ and $\pi_\beta$ as the policies deploying parameters $\theta_\alpha$ and $\theta_\beta$, respectively. Let $\pi_\theta$ denote the policy using parameters $\theta = (\theta_1, \dots, \theta_T)$. Then
\begin{equation*}
    \begin{aligned}
        \left\| \nabla_{\theta_t} l(\theta_\alpha) - \nabla_{\theta_t} l(\theta_\beta) \right\|_2 \le \left\| \nabla_{\theta_t} l(\theta_\alpha) - \nabla_{\theta_t} l(\theta_\beta) \right\|_1 = \underbrace{ \left| \frac{\partial}{\partial \underline{\theta}_t} l(\theta_\alpha) - \frac{\partial}{\partial \underline{\theta}_t} l(\theta_\beta) \right| }_{(\text{I})} + \underbrace{ \left| \frac{\partial}{\partial \bar{\theta}_t} l(\theta_\alpha) - \frac{\partial}{\partial \bar{\theta}_t} l(\theta_\beta) \right| }_{(\text{II})}.
    \end{aligned}
\end{equation*}

For the first part (I), we can derive the following inequality by Proposition~\ref{proposition: cash balance gradient formulation},
\begin{equation*}
    \begin{aligned}
        (\text{I}) &= \left| \E_{s_t \sim \rho_t(\cdot|\pi_\theta)} \left[ \mathbf{1}(\underline{\theta}_t \ge s_t) \times \left( \E_{D_t} \left[ (V_{t+1}^{\pi_{\alpha}})' (\underline{\theta}_t - D_t) \right] - \E_{D_t} \left[ (V_{t+1}^{\pi_{\beta}})' (\underline{\theta}_t - D_t) \right] \right) \right] \right|\\
        &\le \left| \E_{D_t} \left[ (V_{t+1}^{\pi_{\alpha}})' (\underline{\theta}_t - D_t) \right] - \E_{D_t} \left[ (V_{t+1}^{\pi_{\beta}})' (\underline{\theta}_t - D_t) \right] \right|.
    \end{aligned}
\end{equation*}
Applying (\ref{cash balance: recursive gradient}) recursively, we have
\begin{equation}
\label{cash balance: bounded sequential gradient difference 1}
    \begin{aligned}
        (\text{I}) &\le \left| \E_{D_t} \left[ (V_{t+1}^{\pi_{\alpha}})' (\underline{\theta}_t - D_t) \right] - \E_{D_t} \left[ (V_{t+1}^{\pi_{\beta}})' (\underline{\theta}_t - D_t) \right] \right|\\
        & = \biggl| \E_{D_t} \Bigl[ -k \mathbf{1}(\underline{\theta}_t - D_t \le \underline{\theta}_{t+1}) + q \mathbf{1}(\underline{\theta}_t - D_t \ge \bar{\theta}_{t+1})\\
        & \quad\quad\quad + \mathbf{1}(\underline{\theta}_{t+1} < \underline{\theta}_t - D_t < \bar{\theta}_{t+1}) \times \left( L'_{t+1}(\underline{\theta}_t - D_t) + \E_{D_{t+1}}\left[ (V^{\pi_\alpha}_{t+2})'(\underline{\theta}_t - D_{[t:t+1]}) \right] \right) \Bigr]\\
        & \quad - \E_{D_t} \Bigl[ -k \mathbf{1}(\underline{\theta}_t - D_t \le \underline{\theta}_{t+1}) + q \mathbf{1}(\underline{\theta}_t - D_t \ge \bar{\theta}_{t+1})\\
        & \quad\quad\quad + \mathbf{1}(\underline{\theta}_{t+1} < \underline{\theta}_t - D_t < \bar{\theta}_{t+1}) \times \left( L'_{t+1}(\underline{\theta}_t - D_t) + \E_{D_{t+1}}\left[ (V^{\pi_\beta}_{t+2})'(\underline{\theta}_t - D_{[t:t+1]}) \right] \right) \Bigr] \biggr|\\
        &\le \E_{D_{[t:t+1]}} \left[ \left| (V_{t+2}^{\pi_{\alpha}})' (\underline{\theta}_t - D_{[t:t+1]}) - (V_{t+2}^{\pi_{\beta}})' (\underline{\theta}_t - D_{[t:t+1]}) \right| \right]\\
        & \dots\\
        &\le \E_{D_{[t:k-1]}} \left[ \left| (V_{k}^{\pi_{\alpha}})' (\underline{\theta}_t - D_{[t:k-1]}) - (V_{k}^{\pi_{\beta}})' (\underline{\theta}_t - D_{[t:k-1]}) \right| \right].
    \end{aligned}
\end{equation}
Therefore, we can derive the following inequality using (\ref{cash balance: recursive gradient}):
\begin{equation*}
    \begin{aligned}
        (V_{k}^{\pi_{\alpha}})' (\underline{\theta}_t - D_{[t:k-1]}) &= -k \mathbf{1}(\underline{\theta}_t - D_{[t:k-1]} \le \underline{\theta}_k) + q \mathbf{1}(\underline{\theta}_t - D_{[t:k-1]} \ge \bar{\theta}_k)\\
        & \quad + \mathbf{1}(\underline{\theta}_k < \underline{\theta}_t - D_{[t:k-1]} < \bar{\theta}_k) \times \left( L'_k(\underline{\theta}_t - D_{[t:k-1]}) + \E_{D_k}\left[ (V^{\pi_\alpha}_{k+1})'(\underline{\theta}_t - D_{[t:k]}) \right] \right)\\
        &\overset{(a)}{=} -k \mathbf{1}(\underline{\theta}_t - D_{[t:k-1]} \le \underline{\theta}_k) + q \mathbf{1}(\underline{\theta}_t - D_{[t:k-1]} \ge \bar{\theta}_k)\\
        & \quad + \mathbf{1}(\underline{\theta}_k < \underline{\theta}_t - D_{[t:k-1]} < \bar{\theta}_k) \times \left( L'_k(\underline{\theta}_t - D_{[t:k-1]}) + \E_{D_k}\left[ (V^*_{k+1})'(\underline{\theta}_t - D_{[t:k]}) \right] \right)\\
        &\overset{(b)}{=} -k \mathbf{1}(\underline{\theta}_t - D_{[t:k-1]} \le \underline{\theta}_k) + q \mathbf{1}(\underline{\theta}_t - D_{[t:k-1]} \ge \bar{\theta}_k)\\
        & \quad + \mathbf{1}(\underline{\theta}_k < \underline{\theta}_t - D_{[t:k-1]} < \bar{\theta}_k) \times \left( f_k'(\underline{\theta}_t - D_{[t:k-1]}) \right).
    \end{aligned}
\end{equation*}
Here Equation (a) holds because $\pi_\alpha$ uses optimal $\theta_{[k+1:T]}^*$ starting from period $k+1$, and Equation (b) comes from the definition of $f_k$. Again, by definitions of $\underline{f}_k$ and $\bar{f}_k$, we have
\begin{equation}
\label{cash balance: bounded sequential gradient difference 2}
    \begin{aligned}
        (V_{k}^{\pi_{\alpha}})' (\underline{\theta}_t - D_{[t:k-1]}) & = -k \mathbf{1}(\underline{\theta}_t - D_{[t:k-1]} \le \underline{\theta}_k) + q \mathbf{1}(\underline{\theta}_t - D_{[t:k-1]} \ge \bar{\theta}_k)\\
        & \quad + \mathbf{1}(\underline{\theta}_k < \underline{\theta}_t - D_{[t:k-1]} < \bar{\theta}_k) \times \left( f_k'(\underline{\theta}_t - D_{[t:k-1]}) \right)\\
        &= f_k'(\underline{\theta}_t - D_{[t:k-1]}) - \mathbf{1}(\underline{\theta}_t - D_{[t:k-1]} \le \underline{\theta}_k) \times \underline{f}_k'(\underline{\theta}_t - D_{[t:k-1]})\\
        & \quad - \mathbf{1}(\underline{\theta}_t - D_{[t:k-1]} \ge \bar{\theta}_k) \times \bar{f}_k'(\underline{\theta}_t - D_{[t:k-1]}).
    \end{aligned}
\end{equation}

Similarly, we can derive that
\begin{equation}
\label{cash balance: bounded sequential gradient difference 3}
    \begin{aligned}
        (V_{k}^{\pi_{\beta}})' (\underline{\theta}_t - D_{[t:k-1]}) &= f_k'(\underline{\theta}_t - D_{[t:k-1]}) - \mathbf{1}(\underline{\theta}_t - D_{[t:k-1]} \le \underline{\theta}_k^*) \times \underline{f}_k'(\underline{\theta}_t - D_{[t:k-1]})\\
        & \quad - \mathbf{1}(\underline{\theta}_t - D_{[t:k-1]} \ge \bar{\theta}_k^*) \times \bar{f}_k'(\underline{\theta}_t - D_{[t:k-1]}).
    \end{aligned}
\end{equation}

Plugging the results of (\ref{cash balance: bounded sequential gradient difference 2}) and (\ref{cash balance: bounded sequential gradient difference 3}) into (\ref{cash balance: bounded sequential gradient difference 1}), we have
\begin{equation*}
    \begin{aligned}
        (\text{I}) &\le \E_{D_{[t:k-1]}} \left[ \left| (V_{k}^{\pi_{\alpha}})' (\underline{\theta}_t - D_{[t:k-1]}) - (V_{k}^{\pi_{\beta}})' (\underline{\theta}_t - D_{[t:k-1]}) \right| \right]\\
        &=\E_{D_{[t:k-1]}} \biggl[ \Bigl| \underline{f}_k'(\underline{\theta}_t - D_{[t:k-1]}) \times \left( \mathbf{1}(\underline{\theta}_t - D_{[t:k-1]} \le \underline{\theta}_k) - \mathbf{1}(\underline{\theta}_t - D_{[t:k-1]} \le \underline{\theta}_k^*) \right) \\
        & \quad \quad \quad \quad \quad + \bar{f}_k'(\underline{\theta}_t - D_{[t:k-1]}) \times \left( \mathbf{1}(\underline{\theta}_t - D_{[t:k-1]} \ge \bar{\theta}_k) - \mathbf{1}(\underline{\theta}_t - D_{[t:k-1]} \ge \bar{\theta}_k^*) \right) \Bigr| \biggr]\\
        & \le \underbrace{\E_{D_{[t:k-1]}} \left[ \left| \underline{f}_k'(\underline{\theta}_t - D_{[t:k-1]}) \times \left( \mathbf{1}(\underline{\theta}_t - D_{[t:k-1]} \le \underline{\theta}_k) - \mathbf{1}(\underline{\theta}_t - D_{[t:k-1]} \le \underline{\theta}_k^*) \right) \right| \right]}_{(\text{III})} \\
        & \quad + \underbrace{\E_{D_{[t:k-1]}} \left[ \left| \bar{f}_k'(\underline{\theta}_t - D_{[t:k-1]}) \times \left( \mathbf{1}(\underline{\theta}_t - D_{[t:k-1]} \ge \bar{\theta}_k) - \mathbf{1}(\underline{\theta}_t - D_{[t:k-1]} \ge \bar{\theta}_k^*) \right) \right| \right]}_{(\text{IV})}.
    \end{aligned}
\end{equation*}
For part (III), without loss of generality, we assume $\underline{\theta}_k \le \underline{\theta}_k^*$. Then
\begin{equation*}
    \begin{aligned}
        (\text{III}) = \E_{D_{[t:k-1]}} \left[ \left| \underline{f}_k'(\underline{\theta}_t - D_{[t:k-1]}) \times \mathbf{1}(\underline{\theta}_k < \underline{\theta}_t - D_{[t:k-1]} \le \underline{\theta}_k^*) \right| \right] \overset{(a)}{=} \int_{\underline{\theta}_k}^{\underline{\theta}_k^*} - \underline{f}_k' (x) \psi(x) dx,
    \end{aligned}
\end{equation*}
where $\psi$ is the probability density function of $\underline{\theta}_t - D_{[t:k-1]}$. Equation (a) holds because $\underline{f}_k(\cdot)$ is convex, $\underline{\theta}_k^*$ is its minimizer, and it uses the variable change $x = \underline{\theta}_t - D_{[t:k-1]}$. From Assumption~\ref{assumption: cash balance}.\ref{assumption: cash balance Lipschitz of rv}, the cumulative distribution function of the random demand $D_t$ is $L_D$-Lipschitz continuous. Then the probability density function of $D_t$ is upper bounded by $L_D$. Using a similar derivation in Section~\ref{section: inventory system}, the probability density function of cumulative demands is upper bounded by $L_D$ and thus $\psi(\cdot) \le L_D$. Hence,
\begin{equation*}
    \begin{aligned}
        (\text{III}) \le \ & L_D \int_{\underline{\theta}_k}^{\underline{\theta}_k^*} - \underline{f}_k' (x) dx = L_D \left( \underline{f}_k (\underline{\theta}_k) - \underline{f}_k (\underline{\theta}_k^*) \right).
    \end{aligned}
\end{equation*}
The same result holds when $\underline{\theta}_k > \underline{\theta}_k^*$. As for part (IV), we can derive a similar bound
\begin{equation*}
    (\text{IV}) \le L_D \left( \bar{f}_k(\bar{\theta}_k) - \bar{f}_k(\bar{\theta}_k^*) \right).
\end{equation*}
As a result, we have
\begin{equation*}
    (\text{I}) \le (\text{III}) + (\text{IV}) \le L_D \left( \underline{f}_k (\underline{\theta}_k) - \underline{f}_k (\underline{\theta}_k^*) \right) + L_D \left( \bar{f}_k(\bar{\theta}_k) - \bar{f}_k(\bar{\theta}_k^*) \right).
\end{equation*}

For the second part (II), we can derive a similar bound using the same technique.
\begin{equation*}
    (\text{II}) \le L_D \left( \underline{f}_k (\underline{\theta}_k) - \underline{f}_k (\underline{\theta}_k^*) \right) + L_D \left( \bar{f}_k(\bar{\theta}_k) - \bar{f}_k(\bar{\theta}_k^*) \right).
\end{equation*}
This implies that
\begin{equation}
\label{cash balance: sequential decomposition inequality 1}
    \left\| \nabla_{\theta_t} l(\theta_\alpha) - \nabla_{\theta_t} l(\theta_\beta) \right\|_2 \le (\text{I}) + (\text{II}) \le 2L_D \left( \underline{f}_k (\underline{\theta}_k) - \underline{f}_k (\underline{\theta}_k^*) + \bar{f}_k(\bar{\theta}_k) - \bar{f}_k(\bar{\theta}_k^*) \right).
\end{equation}
Recalling the definition of $F_k(\theta_k)$, we have
\begin{equation}
\label{cash balance: sequential decomposition inequality 2}
    \begin{aligned}
        F_k(\theta_k) - F_k(\theta_k^*) = & \ \underbrace{\E_{s_k\sim\rho_k(\cdot|\pi_\theta)} \left[ c(s_k\vee\underline{\theta}_k , s_k ) + f_k (s_k\vee\underline{\theta}_k ) - c(s_k\vee\underline{\theta}_k^*, s_k ) - f_k(s_k\vee\underline{\theta}_k^*) \right]}_{(\text{I})}\\
        & \quad\quad\quad\quad + \underbrace{\E_{s_k\sim\rho_k(\cdot|\pi_\theta)} \left[ c(s_k\wedge \bar{\theta}_k, s_k) + f_k(s_k\wedge\bar{\theta}_k) - c(s_k\wedge \bar{\theta}_k^*, s_k) - f_k(s_k\wedge\bar{\theta}_k^*) \right]}_{(\text{II})}.
    \end{aligned}
\end{equation}
Without loss of generality, we assume that $\underline{\theta}_k \le \underline{\theta}_k^*$. Then we have
\begin{equation*}
    \begin{aligned}
        (\text{I}) & = \int_{-\infty}^{\underline{\theta}_k} \left( \underline{f}_k(\underline{\theta}_k) - \underline{f}_k(\underline{\theta}_k^*) \right) d \rho_k(s_k|\pi_\theta) + \int_{\underline{\theta}_k}^{\underline{\theta}_k^*} \left( f_k(s_k) + ks_k - \underline{f}_k(\underline{\theta}_k^*) \right) d \rho_k(s_k|\pi_\theta)\\
        &\ge \int_{-\infty}^{\underline{\theta}_k} \left( \underline{f}_k(\underline{\theta}_k) - \underline{f}_k(\underline{\theta}_k^*) \right) d \rho_k(s_k|\pi_\theta).
    \end{aligned}
\end{equation*}
The last inequality holds as $\underline{f}_k(s_k) \ge \underline{f}_k(\underline{\theta}_k^*)$ for any $s_k \in [\underline{\theta}_k, \underline{\theta}_k^*]$. Therefore, we have
\begin{equation}
\label{cash balance: sequential decomposition inequality 3}
    \begin{aligned}
        (\text{I}) \ge \P(s_k \le \underline{\theta}_k) \left( \underline{f}_k(\underline{\theta}_k) - \underline{f}_k(\underline{\theta}_k^*) \right) \ge \alpha \left( \underline{f}_k(\underline{\theta}_k) - \underline{f}_k(\underline{\theta}_k^*) \right).
    \end{aligned}
\end{equation}
Similar results hold when $\underline{\theta}_k > \underline{\theta}_k^*$. For term (II), we apply the same technique and derive that
\begin{equation}
\label{cash balance: sequential decomposition inequality 4}
    (\text{II}) \ge \alpha \left( \bar{f}_k(\bar{\theta}_k) - \bar{f}_k(\bar{\theta}_k^*) \right).
\end{equation}
Combining (\ref{cash balance: sequential decomposition inequality 1}), (\ref{cash balance: sequential decomposition inequality 2}), (\ref{cash balance: sequential decomposition inequality 3}), and (\ref{cash balance: sequential decomposition inequality 4}), we conclude that
\begin{equation*}
    \begin{aligned}
        \left\| \nabla_{\theta_t} l(\theta_\alpha) - \nabla_{\theta_t} l(\theta_\beta) \right\|_2 &\le 2L_D \left( \underline{f}_k (\underline{\theta}_k) - \underline{f}_k (\underline{\theta}_k^*) + \bar{f}_k(\bar{\theta}_k) - \bar{f}_k(\bar{\theta}_k^*) \right)\\
        &\le \frac{2 L_D}{\alpha} \left( F_k(\theta_k) - F_k(\theta_k^*) \right).
    \end{aligned}
\end{equation*}
This concludes the proof. \Halmos

\section{Additional Experiments}\label{Appendix: Robustness Check}

This section presents two robustness checks that complement the numerical experiments in Section~\ref{sec: numerical exp} and assess the stability of the policy gradient methods’ performance. Specifically, the first check considers discrete demand, under which Assumption~\ref{assumption: inventory}.\ref{assumption: inventory Lipschitz of rv} does not hold. The second check considers settings in which Assumption~\ref{assumption: inventory}.\ref{assumption: inventory KL} is violated. In both cases, policy gradient methods continue to perform well and remain competitive relative to the benchmarks.

To demonstrate the robustness of policy gradient methods under nonsmooth objectives, we additionally test a Poisson demand family (PO). Because the Poisson distribution has discrete support, Assumption~\ref{assumption: inventory}.\ref{assumption: inventory Lipschitz of rv} does not hold; consequently, the policy gradient objective $l(\theta)$ need not be smooth. Specifically, for each period $t \in [T]$ and state $i \in \mathcal{I}$, we set $D_{t,i} \sim \mathrm{Poisson}(\lambda_{t,i})$, where $\lambda_{t,i}$ are drawn independently from $\mathrm{Unif}[1,10]$ for all $t \in [T]$ and $i \in \mathcal{I}$. All remaining parameters and implementation details follow those in Table~\ref{table:exp markov}. The results, reported in Table~\ref{table:exp markov nonsmooth}, show that policy gradient methods remain effective even when the smoothness assumption fails. Almost all the suboptimality gaps are within $0.1$, and the longest runtime is below $21$ seconds, including the largest horizon $T=100$.

\begin{table}[htbp]
    \small
    \centering
    \caption{Performance of PG for inventory models with Markov-modulated demand under nonsmooth objectives.}
    \renewcommand{\arraystretch}{1.2}
    \begin{tabular}{cc *{3}{ccc}}
    \hline\hline
    \multirow{2}{*}{Problem Setting}
      & \multirow{2}{*}{OP}
      & \multirow{2}{*}{DP Runtime}
      & \multicolumn{3}{c}{Gap}
      & \multicolumn{3}{c}{PG Runtime}\\
    \cline{4-9}
      & & 
      & MIN & AVG & MAX
      & MIN & AVG & MAX \\
    \hline
    $(20, 4, 0.1, \text{PO})$ & 9.1893 & 1068.976
    & 0.0070 & 0.0086 & 0.0112
    & 3.2482 & 3.2756 & 3.3163 \\
    $(20, 4, 0.25, \text{PO})$ & 17.7866 & 1049.034
    & 0.0125 & 0.0166 & 0.0211
    & 3.2539 & 3.2944 & 3.3287 \\
    $(20, 7, 0.1, \text{PO})$ & 9.4902 & 1854.250
    & 0.0119 & 0.0149 & 0.0168
    & 4.2280 & 4.2752 & 4.3201 \\
    $(20, 7, 0.25, \text{PO})$ & 18.3534 & 1881.284
    & 0.0198 & 0.0242 & 0.0278
    & 4.1963 & 4.2300 & 4.2766 \\
    $(50, 4, 0.1, \text{PO})$ & 22.4864 & 2817.266
    & 0.0195 & 0.0244 & 0.0266
    & 7.8664 & 7.9813 & 8.1345 \\
    $(50, 4, 0.25, \text{PO})$ & 42.6924 & 2673.256
    & 0.0397 & 0.0438 & 0.0474
    & 7.8684 & 7.9776 & 8.1850 \\
    $(50, 7, 0.1, \text{PO})$ & 22.8846 & 4716.023
    & 0.0327 & 0.0365 & 0.0407
    & 10.2765 & 10.4378 & 10.7380\\
    $(50, 7, 0.25, \text{PO})$ & 43.4531 & 4702.971
    & 0.0479 & 0.0533 & 0.0584
    & 10.2036 & 10.3700 & 10.6808\\
    $(100, 4, 0.1, \text{PO})$ & 45.3027 & 6404.459
    & 0.0418 & 0.0457 & 0.0506
    & 15.4568 & 15.5592 & 15.6562\\
    $(100, 4, 0.25, \text{PO})$ & 85.4369 & 6215.914
    & 0.0757 & 0.0838 & 0.0912
    & 15.9266 & 16.0616 & 16.3902\\
    $(100, 7, 0.1, \text{PO})$ & 45.6078 & 11128.036
    & 0.0689 & 0.0735 & 0.0774
    & 20.1604 & 20.3228 & 20.5808\\
    $(100, 7, 0.25, \text{PO})$ & 86.0153 & 11134.640
    & 0.0989 & 0.1098 & 0.1220
    & 20.3210 & 20.4842 & 20.7948\\
    
    \hline\hline
    \end{tabular}
    \label{table:exp markov nonsmooth}
\end{table}

We conduct additional experiments to evaluate the robustness of policy gradient (PG) when Assumption~\ref{assumption: inventory}.\ref{assumption: inventory KL} is violated. Table~\ref{table:exp markov assumption} reports results for the Markov-modulated demand model using the same notation and parameters as in Table~\ref{table:exp markov}. Here, BT denotes that demands follow a Beta distribution. Since Beta random variables have bounded support $[0,1]$, choosing a large bound (e.g., $B=20$) implies $P_D(B\mid i)=1$ for all $i\in\mathcal{I}$, and hence violates the second assumption.  Even without the second assumption, PG continues to perform well: in the worst case, it terminates within $22$ seconds and returns a solution with a suboptimality gap below $0.008$. These results suggest that PG is empirically robust to violations of the second assumption.

\begin{table}[htbp]
    \small
    \centering
    \caption{Performance of PG for inventory models with Markov-modulated demand when Assumption~\ref{assumption: inventory}.\ref{assumption: inventory KL} is violated.}
    \renewcommand{\arraystretch}{1.2}
    \begin{tabular}{cc *{3}{ccc}}
    \hline\hline
    \multirow{2}{*}{Problem Setting}
      & \multirow{2}{*}{OP}
      & \multirow{2}{*}{DP Runtime}
      & \multicolumn{3}{c}{Gap}
      & \multicolumn{3}{c}{PG Runtime}\\
    \cline{4-9}
      & & 
      & MIN & AVG & MAX
      & MIN & AVG & MAX \\
    \hline
    $(20, 4, 0.1, \text{BT})$ & 7.9685 & 1060.031
    & 0.0002 & 0.0003 & 0.0004
    & 3.2336 & 3.2838 & 3.3719 \\
    $(20, 4, 0.25, \text{BT})$ & 19.7695 & 1054.023
    & 0.0005 & 0.0008 & 0.0011
    & 3.2113 & 3.2477 & 3.3005 \\
    $(20, 7, 0.1, \text{BT})$ & 7.9174 & 1865.260
    & 0.0003 & 0.0003 & 0.0005
    & 4.1837 & 4.2472 & 4.2925 \\
    $(20, 7, 0.25, \text{BT})$ & 19.6442 & 1942.469
    & 0.0005 & 0.0007 & 0.0010
    & 4.1005 & 4.2155 & 4.2952 \\
    $(50, 4, 0.1, \text{BT})$ & 10.3318 & 2565.268
    & 0.0015 & 0.0022 & 0.0028
    & 7.8706 & 7.9472 & 8.0220 \\
    $(50, 4, 0.25, \text{BT})$ & 25.2732 & 2602.869
    & 0.0020 & 0.0026 & 0.0032
    & 7.8775 & 7.9869 & 8.2997 \\
    $(50, 7, 0.1, \text{BT})$ & 9.9724 & 4725.112
    & 0.0008 & 0.0011 & 0.0014
    & 10.4159 & 10.4715 & 10.5342 \\
    $(50, 7, 0.25, \text{BT})$ & 24.3839 & 4622.234
    & 0.0021 & 0.0026 & 0.0034
    & 10.4778 & 10.5351 & 10.6268 \\
    $(100, 4, 0.1, \text{BT})$ & 11.4952 & 6257.031
    & 0.0020 & 0.0024 & 0.0028
    & 15.6536 & 15.7692 & 16.1654\\
    $(100, 4, 0.25, \text{BT})$ & 27.4703 & 6285.321
    & 0.0053 & 0.0062 & 0.0073
    & 15.9749 & 16.0568 & 16.1602\\
    $(100, 7, 0.1, \text{BT})$ & 11.3371 & 11104.514
    & 0.0021 & 0.0024 & 0.0026
    & 20.8733 & 20.9689 & 21.2646\\
    $(100, 7, 0.25, \text{BT})$ & 27.0777 & 11117.824
    & 0.0054 & 0.0062 & 0.0068
    & 20.6542 & 20.8212 & 21.1322 \\
    \hline\hline
    \end{tabular}
    \label{table:exp markov assumption}
\end{table}

\end{APPENDICES}

\end{document}